\documentclass[10pt]{article}

\usepackage{amssymb,amsmath,amsthm,amsfonts}
\usepackage{graphicx,graphics}
\usepackage{epstopdf}
\usepackage{algorithm,algorithmic}
\graphicspath{{./pics/}}
\usepackage{verbatim}
\usepackage{mathtools}
\usepackage{threeparttable}
\usepackage{booktabs}
\usepackage{url}
\usepackage[nohead,margin=1.1in]{geometry}
\usepackage{color}
\usepackage{tikz}
\usetikzlibrary{calc}
\numberwithin{equation}{section}
\numberwithin{algorithm}{section}

\newtheorem{theorem}{Theorem}[section]
\newtheorem{lemma}{Lemma}[section]

\newtheorem{corollary}{Corollary}[section]

\newtheorem{remark}{Remark}[section]

\newtheorem{assumption}{Assumption}[section]

\allowdisplaybreaks
\newcommand{\E}{\mathbb{E}}
\newcommand{\nn}{\rm nn}
\newcommand{\bd}{\tilde d}
\def\Lip#1{L_{#1}}

\begin{document}
\title{On the Approximation of Bi-Lipschitz Maps by Invertible Neural Networks\thanks{The work of BJ is supported by a start-up
fund and Direct Grant for Research 2022/2023, both from The Chinese University of Hong Kong and Hong Kong RGC General Research Fund (Project 14306423), and that of JZ was
substantially supported by Hong Kong RGC General Research Fund (Projects 14306921 and 14308322).}}

\author{Bangti Jin\thanks{Department of Mathematics, The Chinese University of Hong Kong, Shatin, New Territories,
Hong Kong. (\texttt{b.jin@cuhk.edu.hk}, \texttt{bangti.jin@gmail.com}, \texttt{zou@math.cuhk.edu.hk})} \and Zehui Zhou\thanks{Department of Mathematics,
Rutgers University, 110 Frelinghuysen Road Piscataway, NJ 08854-8019 (\texttt{zz569@math.rutgers.edu})}
\and Jun Zou\footnotemark[2]}
\date{}
\maketitle

\begin{abstract}
Invertible neural networks (INNs) represent an important class of deep neural network architectures
that have been widely used in several applications. The universal approximation properties
of INNs have also been established recently. However, the approximation rate of INNs is largely
missing. In this work, we provide an analysis of the capacity of a class of coupling-based INNs
to approximate bi-Lipschitz continuous mappings on a compact domain, and the result
shows that it can well approximate both forward and inverse maps simultaneously. Furthermore,
we develop an approach for approximating bi-Lipschitz maps on infinite-dimensional spaces that simultaneously approximate the forward and
inverse maps, by combining model reduction with principal component analysis and INNs for approximating the reduced map, and we analyze the overall approximation error of the approach. Preliminary
numerical results show the feasibility of the approach for approximating the solution operator for
parameterized second-order elliptic problems.

\end{abstract}

\section{Introduction}

Invertible neural networks (INNs) are a class of neural network (NN) architectures with invertibility by design, via
 special invertible layers called the flow layers. INNs often enjoy tractable numerical algorithms to compute the inverse map and
Jacobian determinant, e.g., with explicit inversion formulas. These distinct features have made them very attractive
for a variety of machine learning tasks, e.g., generative modeling \cite{DinhKruegerBengio:2015,KingmaDhariwal:2018,KimLee:2019}, probabilistic modeling \cite{LouizosWelling:2017,DinhSohlBengio:2017,HuangKrueger:2018,BauerMnih:2019}, solving inverse problems \cite{ArdizzoneKruse:2019,Zabaras:2021,Arndt:2023}, modeling nonlinear dynamics \cite{Bevanda:2022} and point cloud generation \cite{Yang:2019}.

There are several different classes of INNs, including invertible residual networks (iResNet) \cite{BehrmannGrathwohlJacobsen:2019,YamazakiRathourLe:2021}, neural ordinary
differential equations (NODEs) \cite{ChangMeng:2018,ChenRubanovaDuvenaud:2018,DupontDoucet:2019} and
coupling-based neural networks \cite{DinhKruegerBengio:2015,DinhSohlBengio:2017,JacobsenSmeulders:2018,KingmaDhariwal:2018,
ArdizzoneKruse:2019}. For iResNet, Behrmann et al \cite{BehrmannGrathwohlJacobsen:2019} leveraged the viewpoint of ResNets as an Euler
discretization of ODEs and proved the standard ResNet architecture can be made invertible by adding a simple normalization step to control the Lipschitz constant of the NN during training. The inverse is not available in closed form but can be obtained through a fixed-point iteration. Chen et al \cite{ChenRubanovaDuvenaud:2018} proposed using black-box ODE solvers as a model component, and developed a class of new models, i.e., NODEs, for time-series modeling, supervised learning, and density estimation etc. NODEs indirectly models an invertible function by transforming an input vector through an ordinary differential equation (ODE). Dupont and Doucet \cite{DupontDoucet:2019} introduced
a class of more expressive and empirically stable models, augmented neural ODEs (ANODEs), which have a lower computational cost.
There is another family of popular INNs using coupling-based frameworks, including non-linear independent component estimation
(NICE) \cite{DinhKruegerBengio:2015}, real-valued non-volume preserving (RealNVPs) \cite{DinhSohlBengio:2017}, invertible
reversible residual network (i-RevNet) \cite{JacobsenSmeulders:2018}, and generative flow (Glow) \cite{KingmaDhariwal:2018}. In coupling-based NNs, one employs a highly restricted NN architecture where only some of the input variables undergo some transformations, and the rest of the input variables become the output as-is so that the splitting facilitates the inversion process.

The great empirical successes of INNs in diverse applications have sparked intensive interest in theoretical study of INNs, e.g., approximation theory and training. The construction of INNs imposes certain restrictions on the choice of admissible layers, and hence their approximation properties are more limited and far less well understood than that for standard fully connected NNs (see \cite{DeVoreHaninPetrova:2021} for an overview of the latter). Indeed, only very recently, Teschima et al \cite{TeshimaIshikawaSugiyama:2020} proved a first universal approximation property in $\mathcal{L}^p$ spaces
for coupling based INNs. Lyu et al \cite{LyuChen:2022} extended the results to $C^k$ universality for $C^k$ diffeomorphism and also established the parametric version. More recently, Ishikawa et al \cite{IshikawaTeshimaSugiyama:2022} developed a more general theoretical framework, covering also other types of INNs, e.g., NODEs. However, approximation rates of INNs appear unknown so far.

In this work, we continue the line of research on the approximation theory for INNs. We explicitly construct a coupling-based INN for simultaneously approximating the forward and inverse maps of a bi-Lipschitz map on a compact domain, by combining a shallow NN with a deep coupling-based NN, and provide an analysis of the approximation capacity of the constructed INNs. We show that quantitatively it can approximate the forward and inverse processes simultaneously to a certain tolerance; see Theorem \ref{thm:F_NN} for the precise statement. This result complements the existing results \cite{TeshimaIshikawaSugiyama:2020,LyuChen:2022,IshikawaTeshimaSugiyama:2022} by providing an explicit construction and also quantitative estimates. To the best of our knowledge, this result is the first of the kind in the literature. Methodologically, the analysis is inspired by the techniques for establishing approximation capabilities of INNs / NODEs \cite{AroraBasuMukherjee:2018,LiLinShen:2019,TeshimaIshikawaSugiyama:2020}.

To illustrate the use of the approximation result, we apply INNs to approximate bi-Lipschitz maps in infinite-dimensional spaces, which belongs to operator learning and represents a fundamental task in the study of many scientific problems \cite{Krasnoselski:1972}. We focus on the situation where the ambient dimensionality is high but the intrinsic dimensionality is low, so that model reduction can be applied effectively. There are many methods for feature projection to reduce the dimensionality, e.g., principal component analysis (PCA) \cite{Hotelling:1933,BerkoozHolmesLumley:1993}, nonnegative matrix factorization \cite{LeeSeung:1999} and Johnson-Lindenstrauss embedding \cite{JohnsonLindenstrauss:1984}. We propose to couple INN with PCA, i.e., first applying PCA encoder and decoder to reduced the bi-Lipschitz map in infinite-dimensional spaces to finite-dimensional subspaces, and then approximating the reduced map with INNs. We provide an analysis of the overall approximation error. Further, we provide preliminary numerical
illustrations with model problems in parametric PDEs.

The idea of coupling of DNNs  with model reduction has been explored in several works (see, e.g., \cite{WojcikKurdziel:2019,KarnikWangIwen:2022,PinedaPetersen:2022,BhattacharyaHosseiniStuart:2021,LiuYangLiao:2022,ChenWangYang:2023}).
W\'ojcik and Kurdziel \cite{WojcikKurdziel:2019} proposed using Johnson-Lindenstrauss embedding to reduce the input dimensionality, and Karnik et al \cite{KarnikWangIwen:2022} developed complexity bounds for approximating a H\"{o}lder (or uniformly) continuous function under a low-complexity structure assumption. Bhattacharya et al \cite{BhattacharyaHosseiniStuart:2021} proposed an approach using PCA for dimensionality reduction, and then employing a fully connected feedforward NN to approximate the reduced nonlinear mapping. Pineda and Petersen \cite{PinedaPetersen:2022} investigated restricting the domain of the forward operator to finite-dimensional spaces,
and proved the existence of a robust-to-noise NN approximation of the operator. Our work continues the active line of research on developing operator approximations, by extending the encoder-decoder type approach in \cite{BhattacharyaHosseiniStuart:2021} to simultaneously approximating the forward and inverse maps using INNs. Liu et al  \cite{LiuYangLiao:2022,ChenWangYang:2023} thoroughly analyzed general encoder-decoder based approaches in terms of approximation and generalization errors, and applied the results to canonical PDEs. In this work, we use the approximation result of INNs developed herein instead of more well-understood feedforward fully connected NNs, in order to derive the estimation error on the INN based operator approximations.

The rest of the paper is organized as follows. In Section \ref{sec:INN}, we develop an INN architecture for approximating the bi-Lipschitz map, and analyze the approximation error of the constructed INN. In Section \ref{sec:overall}, we develop a novel operator approximation based on model reduction (PCA) and INN, and discuss the error estimation for the overall scheme. Finally, in Section \ref{sec:numer} we present some preliminary numerical
experiments to illustrate the approach.

\section{Invertible neural networks}\label{sec:INN}

In this section, we construct a coupling-based invertible neural network (INN) to approximate the
bi-Lipschitz forward operator $F:\mathbb{R}^d\to \mathbb{R}^d$, $d\geq 2$.
\begin{assumption}\label{ass}
The operator $F$ is bi-Lipschitz continuous with Lipschitz constants $\Lip{F}$ and $\Lip{F^{-1}}$ for $F$ and its inverse $F^{-1}$,
respectively.
\end{assumption}

Coupling-based NNs split the input $x\in\mathbb{R}^d$ into two parts and then employ an affine
coupling layer whose inverse can be found directly. Specifically, fix $k\in \mathbb{N}$ with $k<d$, we split
$x\in\mathbb{R}^d$ into $x=(x_{\leq k},x_{>k})$, with the subvectors $x_{\leq k}:=(x_1,\ldots,x_k)\in \mathbb{R}^{k}$
and $x_{>k}:=(x_{k+1},\ldots,x_d)\in\mathbb{R}^{d-k}$, and set
\begin{align*}
 x'_{\leq k}=& x_{\leq k} \odot \exp(s_1(x_{>k}))+t_1(x_{>k}),\\
 x'_{>k}=& x_{>k} \odot \exp(s_2(x'_{\leq k}))+t_2(x'_{\leq k}),
\end{align*}
with operators $s_1,t_1 : \mathbb{R}^{d-k} \longrightarrow \mathbb{R}^{k}$ and $s_2,t_2 :
\mathbb{R}^{k} \longrightarrow \mathbb{R}^{d-k}$, which are all taken to be NNs. The
inverse of this process can be obtained directly as
\begin{align*}
 x_{>k}=& (x'_{>k}-t_2(x'_{\leq k}))\odot\exp(-s_2(x'_{\leq k})),\\
 x_{\leq k}=& (x'_{\leq k}-t_1(x_{>k}))\odot\exp(-s_1(x_{>k})).
\end{align*}
Intuitively, it reduces the spatial resolution of the coefficients and potentially increasing the layer size may improve the expressivity of the resulting NN. Numerically, it has been shown that the restriction to the affine
coupling layers will not noticeably reduce the expressivity of the NN \cite{ArdizzoneKruse:2019}.
The universal approximation properties of (affine) coupling based INNs were studied in
\cite{TeshimaIshikawaSugiyama:2020,LyuChen:2022,IshikawaTeshimaSugiyama:2022}.
We provide an explicit construction of INNs to approximate bi-Lipschitz maps,
which yields an explicit bound on the INN approximation.

\subsection{Construction}

Now we construct an INN, denoted by $F_{\nn}$, to approximate
the bi-Lipschitz map $F:\mathbb{R}^{d}\to \mathbb{R}^{d}$, over any compact set $K\subset \mathbb{R}^d$.
Throughout we assume that $K=[0,1]^{d}$ and $\{x^i\}_{i=1}^N$ are taken from a uniform grid in $K$.
For a general compact set $K\subset\mathbb{R}^d$, we can always find an affine mapping $A_f:\; x\to W_ax+b_a$
with $W_a\in \mathbb{R}^{d\times d}$ and $b_a\in\mathbb{R}^d$ such that $A_f(K)\subset [0,1]^{d}$.
Let $\{x^i\}_{i=1}^N=\{x^\alpha\}_{\alpha\in[n-1]^{d}}$, where $[n]$
denotes the set $[n] = \{0,\cdots,n\},$ $N=n^{d}$ and $x^\alpha=\frac{\alpha}{n}\in\mathbb{R}^d$. We also denote the
$j$th entry of a vector $x\in\mathbb{R}^d$ by $x_{j}$. For a given set $K$,
we denote by $|K|$ its Lebesgue measure. We denote the standard Euclidean norm of a vector (and also the spectral norm of a matrix)
by $\|\cdot\|_2$. For any $f:K\to\mathbb{R}^d$, we define
\begin{align*}
\|f\|_{\mathcal{L}^p(K)}=&\left\{\begin{aligned}
  \Big(\int_{ K}\|f(x)\|_2^p{\rm d}x\Big)^\frac1p,&\quad 1\leq p<\infty,\\
   \sup_{x\in K}\|f(x)\|_2, & \quad p =\infty.
\end{aligned}\right.
\end{align*}
Following the idea in \cite[Section 4.3]{LiLinShen:2019}, we first approximate the map $F$ by a composition
mapping without using NNs. This is achieved in the following lemma. The function
$h^r$ (depending on the choice of $n$ and $r$) serves a piecewise linear approximation
to the following piecewise constant function:
\begin{equation*}
  \lim_{r\to 1^-} h^r(x) = \left\{\begin{array}{ll}
  \frac{j}{n}, & x\in [\frac{j}{n},\frac{j+1}{n}), \ j \in [n-1],\\
    x, & x \in \mathbb{R}\setminus [0,1).
  \end{array}\right.
\end{equation*}
For any $r\in (0,1)$, $h^r$ is strictly increasing, and can be realized by a
rectified linear unit (ReLU) NN exactly; see Section \ref{subsubsec:H^r}. Then by choosing
$r\in(0,1)$ sufficiently close to 1, we obtain the final estimate.
\begin{lemma}\label{lem:H^r}
For any $r\in(0,1)$, let $H^r: \mathbb{R}^{d}\to\mathbb{R}^{d}$ be defined by $H^r(x)=(h^r(x_{1}),\cdots,
h^r(x_{d}))^t$, with the component function $h^r:\mathbb{R}\to \mathbb{R}$ given by
\begin{align*}
h^r(x)=\left\{\begin{array}{ll}
\frac{j}{n},   & x=\frac{j}{n}, \; j\in[n-1],\\
\frac{j+(1-r)}{n},   & x=\frac{j+r}{n}, \; j\in[n-1],\\
x,                  & x\in \mathbb{R}\setminus [0,1],
\end{array}\right.
\end{align*}
being continuous over $\mathbb{R}$ and linear in the intervals $[\frac{j}{n},\frac{j+r}{n}]$ and
$[\frac{j+r}{n},\frac{j+1}{n}]$ for $j\in[n-1]$. Then the
following estimates hold
\begin{align*}
\|F\circ H^r-F\|^2_{\mathcal{L}^2(K)}&\leq \Lip{F}^2(2r-1)^2dn^{-2},\\
\|H^r\circ F^{-1}- F^{-1}\|^2_{\mathcal{L}^2(F(K))}&\leq \Lip{F}^d (2r-1)^2 d n^{-2}.
\end{align*}
\end{lemma}
\begin{proof}
First, since $|K|=1$, we obtain the following estimate
\begin{align*}
\|F\circ H^r-F\|^2_{\mathcal{L}^2(K)} =& \|F\circ H^r-F\circ \mbox{id}\|^2_{\mathcal{L}^2(K)}= \|F\circ (H^r- \mbox{id})\|^2_{\mathcal{L}^2(K)}\\
\leq&|K|\|F\circ (H^r- \mbox{id})\|^2_{\mathcal{L}^\infty(K)}=\|F\circ (H^r- \mbox{id})\|^2_{\mathcal{L}^\infty(K)}.
\end{align*}
Then, by the definitions of the functions $H^r$ and $h^r$, we can derive
\begin{align*}
&\|F\circ (H^r- \mbox{id})\|_{\mathcal{L}^\infty(K)}
\leq \omega_{F} (\|H^r-\mbox{id}\|_{\mathcal{L}^\infty(K)})
\leq \omega_{F}(d^\frac12 \max_{x\in[0,n^{-1}]}(h^r(x)-x))=\omega_{F}(|2r-1|d^\frac12 n^{-1} ),
\end{align*}
where $\omega_{F}$ denotes the modulus of continuity of the mapping $F$, defined by
$\omega_{F}(s)=\sup_{\|x-\tilde{x}\|_2\leq s}\|F(x)-F(\tilde x)\|_2\leq \Lip{F}s$.
Thus, we obtain
\begin{align*}
\|F\circ H^r-F\|^2_{\mathcal{L}^2(K)}
\leq\Lip{F}^2(2r-1)^2dn^{-2}.
\end{align*}
Similarly, by the inequality $|F(K)|\leq \Lip{F}^{d}|K|=\Lip{F}^{d},$
we can bound $F^{-1}-H^r\circ F^{-1}$.
\end{proof}

\begin{remark}
The estimates in Lemma \ref{lem:H^r} are derived in the $\mathcal{L}^2$ norm, and the same results can
be extended to the $\mathcal{L}^p$ space, with $1\leq p\leq \infty$. By the construction of $h^r$,
$\lim_{r\to 1^-}h^r$ is a piecewise constant function. Since $F\in
\mathcal{L}^2(K)$, we can use piecewise constant mappings to approximate $F$ in $K$. The definition of
$H^r$ implies that $\lim_{r\to 1^-}F\circ H^r$ may serve as such a piecewise constant mapping.
\end{remark}

Next we construct two INNs  $H^r_{\rm nn}$ and $\tilde{F}_{\rm nn}$ to approximate
$H^r$ and $F$, respectively. The detailed procedures are given in the following two parts. Throughout, for any $s\in\mathbb{N}$, the notation $\boldsymbol{1}_s$ and $\boldsymbol{0}_s$ denote a length-$s$ column / row vector of ones and zeros, respectively. The notation $O_{d\times}$ denotes a zero matrix of size $d\times s$.

\subsubsection{INN for $H^r$}\label{subsubsec:H^r}
First we construct an INN $H^r_{\nn}$ to realize the map $H^r$. For any $r\in(0,1)$, the function
$h^r:\mathbb{R}\to\mathbb{R}$ is strictly increasing and hence bijective and continuous piecewise linear with
the non-differentiable points given by $(2n+1)$ tuples $\{(a_i,a_i)\}_{i=0}^{n}\cup\{(p_i,q_i)
\}_{i=1}^{n}$, with $a_i=\frac{i}{n}$, $p_i=\frac{i-1+r}{n}$ and $q_i=\frac{i-r}{n}$. It
can be exactly represented by a two-layer NN with the ReLU activation function $\sigma(t)$
\cite{AroraBasuMukherjee:2018} (see also \cite{HeLiXuZheng:2020}). Specifically,
let $\sigma(t)=\max(0,t)$. Then it follows from the elementary identities
\begin{align*}
\frac{a_i-q_i}{a_i-p_i} - \frac{q_i-a_{i-1}}{p_i-a_{i-1}} = \frac{\frac{i}{n}-\frac{i-r}{n}}{\frac{i}{n}-\frac{i-1+r}{n}} - \frac{\frac{i-r}{n}-\frac{i-1}{n}}{\frac{i-1+r}{n}-\frac{i-1}{n}} = \frac{r}{1-r}-\frac{1-r}{r} = \frac{2r-1}{r(1-r)},\\
\frac{q_{i+1}-a_i}{p_{i+1}-a_i} - \frac{a_i-q_{i}}{a_i-p_{i}} = \frac{\frac{i+1-r}{n}-\frac{i}{n}}{\frac{i+r}{n}-\frac{i}{n}} -\frac{\frac{i}{n}-\frac{i-r}{n}}{\frac{i}{n}-\frac{i-1+r}{n}} = \frac{1-r}{r} - \frac{r}{1-r}= -\frac{2r-1}{r(1-r)},
\end{align*}
for $i=1,\ldots,n-1$, that we can express $h^r$ as
\begin{align*}
h^r(x)=&-\sigma(-x)+\frac{q_1-a_0}{p_1-a_0}\sigma(x-a_0)+\Big(\frac{a_1-q_1}{a_1-p_1}-\frac{q_1-a_0}{p_1-a_0}\Big)\sigma(x-p_1)
+\Big(\frac{q_{2}-a_1}{p_{2}-a_1} - \frac{a_1-q_{1}}{a_1-p_{1}}\Big)\sigma(x-a_1)\\
&+\cdots+\Big(\frac{a_n-q_n}{a_n-p_n}-\frac{q_n-a_{n-1}}{p_n-a_{n-1}}\Big)\sigma(x-p_n)+\Big(1-\frac{a_n-q_n}{a_n-p_n}\Big)\sigma(x-a_n)\\
=&-\sigma(-x)+\frac{1-r}{r}\sigma(x-a_0)+\frac{2r-1}{r(1-r)}\sigma(x-p_1)-\frac{2r-1}{r(1-r)}\sigma(x-a_1)+\cdots\\
&-\frac{2r-1}{r(1-r)}\sigma(x-a_{n-1})+\frac{2r-1}{r(1-r)}\sigma(x-p_n)-\frac{2r-1}{1-r}\sigma(x-a_n)\\
=&-\sigma(-x)+\frac{1-r}{r}\sigma(x-a_0)-c_r\sum_{i=1}^{n-1}\sigma(x-a_{i})-\frac{2r-1}{1-r}\sigma(x-a_n)+c_r\sum_{i=1}^{n}\sigma(x-p_{i}),
\end{align*}
with the constant $c_r=\frac{2r-1}{r(1-r)}$. Now let
\begin{align*}
  W^1 &= (-1,\tfrac{1-r}{r},-c_r\boldsymbol{1}_{n-1},-\tfrac{2r-1}{1-r},c_r\boldsymbol{1}_{n})\in\mathbb{R}^{1\times(2n+2)},\\
  W^0 &= (-1,\boldsymbol{1}_{2n+1})^t\in \mathbb{R}^{(2n+2)\times1},\quad
  b^0 = -(0,a_0,\ldots,a_n,p_1,\ldots,p_n)^t\in \mathbb{R}^{2n+2}.
\end{align*}
Then we have
\begin{equation*}
  h^r(x) = W^1\sigma(W^0x+b^0):=h^r_{\rm nn}(x).
\end{equation*}
Finally, let $H^r_{\rm nn}(x):=(h^r_{\rm nn}(x_{1}),\cdots,h^r_{\rm nn}(x_{d}))$.
Then we have $H^r_{\rm nn}(x)=H^r(x)$, and obviously $H^r_{\rm nn}$ is bijective in $\mathbb{R}^{d}$.
Hence, $H^r_{\rm nn}$ is an INN with a two-layer NN $h^r_{\rm nn}$
with $2(n+1)$ neurons.

\begin{remark}
The preceding discussion focuses on representing $H^r(x)$ exactly using ReLU but the resulting NN
does not belong to the coupling-based family. This can be easily remedied as follows. Indeed, in order to approximate
$H^r(x)=(h^r(x_1),\cdots,h^r(x_{d}))$, we can define a one-dimensional control family $\mathcal{R}$
corresponding to $\mathcal{F}$, cf. \eqref{eq:control family} below, following the technique in \cite[Section 4.3]{LiLinShen:2019} as
$\mathcal{R}:=\{g(x)\;|\; g(x)=[f((x,0,\cdots,0))]_1, \;f\in\mathcal{F}\}$,
and let $\mathcal{A}_{\mathcal{R},1}:=\mathcal{A}_{\mathcal{R}}\times\{{\rm id}\}\times \cdots \times\{{\rm id}\}$.
By \cite[Proposition 4.4]{LiLinShen:2019}, and the definitions of $h^r$ and $\{x_i\}_{i=1}^{N}$, there exists
$\tilde{h}^r \in \mathcal{A}_\mathcal{R}$ such that
\begin{align*}
\|\tilde{h}^r-h^r\|_{\mathcal{F}^\infty([0,1])}\leq \omega_{h^r}(\Delta)=\omega_{h^r}(n^{-1})\leq 2n^{-1},
\end{align*}
where $\Delta$ is the maximum of the minimal distance of each data point with all other data points
and $\omega_{h^r}$ is the modulus of continuity, i.e., $\omega_{h^r}(\Delta)=\sup_{|x-\tilde{x}|\leq|\Delta| }(h^r(x)-h^r(\tilde x))$.
Then for $H^r_{\nn}(x):=(\tilde{h}^r(x_1),\cdots,\tilde{h}^r(x_{d}))$, we can derive
\begin{align}\label{eq:H_NN-Hr}
\|H^r_{\nn}-H^r\|_{\mathcal{F}^\infty(K)}\leq \sqrt{d}\|\tilde{h}^r-h^r\|_{\mathcal{F}^\infty([0,1])}\leq 2\sqrt{d}n^{-1}.
\end{align}
\end{remark}

\subsubsection{INN for $F$}

Now we explicitly construct an INN $\tilde F_{\rm nn}$ to approximate the bi-Lipschitz map $F$. The overall construction is quite lengthy and technical and is divided into several steps. Note that the construction is non-unique: there is an alternative construction by lifting the vector into $\mathbb{R}^{2d+2}$; see the appendix for details. Throughout, we assume that we have access to a finite collection of evaluations $\{(x^{\alpha} , y^{\alpha})\}_{\alpha\in[n-1]^{d}}$ on a uniform grid. The construction proceeds in the following four steps:
\begin{itemize}
  \item[(i)] make a small perturbation $\eta$ to $\{x^{\alpha}\}_{\alpha\in[n-1]^{d}}$ such that
  $\{\eta(x^{\alpha})_{d}\}_{\alpha\in[n-1]^{d}}$ are distinct, i.e., the last component of all
  the mapped points $\{\eta(x^{\alpha})_{d}\}_{\alpha\in[n-1]^{d}}$ are distinct;
  \item[(ii)] construct a mapping $\varphi^N$ that satisfies, for any $\alpha\in[n-1]^{d}$, that $(\varphi^N\circ\eta
  (x^{\alpha}))_{j}=y^\alpha_{j}$ for any $j=1,\cdots,d-1$ and $(\varphi^N\circ \eta(x^{\alpha}))_{d}
  =\eta(x^{\alpha})_{d}$, i.e., keeping the last entry of all the points unchanged;
  \item[(iii)] choose an index $j_0\in\{1,\cdots,d-1\}$, and then make a small perturbation $\tilde{\eta}$,
  which ensures that $\{[\tilde{\eta}\circ \varphi^N\circ\eta (x^{\alpha})]_{j_0}\}_{\alpha\in[n-1]^{d}}$ are
  all distinct;
  \item[(iv)] derive some $\tilde{\varphi}^N$ so that $[\tilde{\varphi}^N\circ\tilde{\eta}\circ \varphi^N\circ
  \eta (x^{\alpha})]_{d}=y^{\alpha}_{d}$ for any $\alpha\in[n-1]^{d}$ and $[\tilde{\varphi}^N\circ\tilde{\eta}\circ
  \varphi^N\circ\eta(x^{\alpha})]_{j}=[\tilde{\eta}\circ\varphi^N\circ\eta(x^{\alpha})]_{j}$ for any $j=1,\cdots,d-1$.
\end{itemize}

The key is to ensure that each of these steps can be realized by a coupling-based
INN, and Lipschitz constants of these INNs and their inverses are under precise
control. The explicit constructions of these invertible networks follow that of coupling based INNs, and the
details are given in the following lemmas. The construction represents the main technical
contribution of the work.

Now we introduce the main technical tool, i.e., control family, in the construction. This concept was
used in the dynamical system approach to universal approximation property of DNNs \cite{LiLinShen:2019}.
We employ a control family $\mathcal{F}$ of three-layer fully connected NN architecture with ReLU $\sigma$ defined by
\begin{equation}\label{eq:control family}
\mathcal{F}:=\big\{V\sigma\big(W^{(2)}\sigma(W^{(1)}\cdot+b^{(1)})+b^{(2)}\big)\;|\; V, W^{(1)}, W^{(2)}\in\mathbb{R}^{d\times d}. b^{(1)},b^{(2)}\in\mathbb{R}^{d}\big\}.
\end{equation}
Accordingly, we define the attainable set $\mathcal{A}_{\mathcal{F}}$ corresponding to $\mathcal{F}$ by
\begin{equation}\label{eq:attain}
\mathcal{A}_{\mathcal{F}}:=\bigcup_{T>0}\Big\{\psi_{\tau_k}^{f_k}\circ\psi_{\tau_{k-1}}^{f_{k-1}}\circ\cdots\circ\psi_{\tau_1}^{f_1}\;|\; k \geq 1, f_j\in\mathcal{F},\tau_j>0, \forall j=1,\cdots,k \;\mbox{and}\;\sum_{j=1}^{k}\tau_j=T\Big\},
\end{equation}
with $\psi_\tau^f=y(\tau)$, where $y(t)$ is the unique solution of the following initial value problem of the
systems of first-order autonomous ordinary differential equations (ODEs):
\begin{equation}\label{eq:ode}
\left\{\begin{aligned}
\frac{{\rm d}y}{{\rm d}t}&=f(y),\quad t> 0,\\
y(0)&=x.
\end{aligned}\right.
\end{equation}
Since $f(y)$ is Lipschitz continuous in $y$, by the classical
Picard-Lindel\"{o}f theorem for ODEs (see, e.g., \cite[Theorem 2.2.1]{Kolokoltsov:2019}),
problem \eqref{eq:ode} has a unique continuously differentiable solution $y(t)$.
Further,  we have
\begin{equation}\label{eq:ode_sol}
\psi_\tau^f(x)=y(\tau;x)=x+\int_{0}^\tau f(y(t;x)) {\rm d} t.
\end{equation}
The map $\psi_\tau^{-f}=z(\tau)$ given by the unique solution $z(t)$ of the following initial value problem
\begin{equation}\label{eq:ode_inv}
\left\{\begin{aligned}
\frac{{\rm d}z}{{\rm d}t}&=-f(z),\quad t> 0,\\
z(0)&=x,
\end{aligned}\right.
\end{equation}
is the inverse of the flow map $\psi_\tau^f$, since
$\psi_\tau^{-f}\circ \psi_\tau^f(x)=\psi_\tau^{f}\circ \psi_\tau^{-f}(x)=x$.
Thus, the functions in the attainable set $\mathcal{A}_\mathcal{F}$ are invertible,
which serve as the natural candidates for constructing
various INNs. Moreover, if $f(x)\in\mathcal{F}$, we have $D f(Ax+b)\in\mathcal{F}$, for any $ D,A\in
\mathbb{R}^{d\times d}$ and $b\in\mathbb{R}^{d}$. Indeed, the constructions below use these
elementary facts heavily.

It follows from the identity \eqref{eq:ode_sol} and direct computation with the chain rule that the
Jacobian $ J_{\psi^{f}_\tau}(x)$ of the flow map $\psi_\tau ^f$ (with respect
to the input $x$) is given by
\begin{equation}\label{eq:Jacobi}
J_{\psi^{f}_\tau}(x)=I+\int_0^\tau J_f(y(t))J_{y(t)}(x){\rm d} t=I+\int_0^\tau J_f(\psi^{f}_t(x))J_{\psi^{f}_t}(x){\rm d} t,
\end{equation}
where $J_f$ denotes the Jacobian of $f$ with respect to its argument.
Moreover, the solution $J_{\psi^{f}_\tau}(x)$ to \eqref{eq:Jacobi} is
given by
\begin{equation}\label{eqn:Jac-form}
J_{\psi^{f}_\tau}(x)=e^{\int_{0}^\tau J_f(\psi^{f}_t(x)){\rm d}t}.
\end{equation}

Let $\Delta_j(x^\alpha)=\min_{\alpha_1\neq\alpha_2}|x^{\alpha_1}_j-x^{\alpha_2}_j|$
be the minimal distance in the $j$th coordinate of the collection of points $\{x^\alpha\}_{
\alpha\in[n-1]^d}$. The construction of the INN $\tilde F_{\rm nn}$ below requires $\Delta_j
(x^\alpha)>0$ for some $j$ at the beginning. The mapping $\eta\in\mathcal{A}_{\mathcal{F}}$
constructed in the following lemma ensures $\Delta_{d}\big(\eta(x^\alpha)\big)=N^{-1}$. Thus, it
fulfills step (i) of the overall procedure.
\begin{lemma}\label{lem:perturb_1}
For the uniform grids $\{x^\alpha=\frac{\alpha}{n}\}_{\alpha\in[n-1]^{d}}$ of the hypercube
$K=[0,1]^{d}$, there exists an invertible mapping $\eta\in\mathcal{A}_{\mathcal{F}}$ such that
\begin{align*}
&(\eta(x))_{j}=x_{j},\quad \forall j=1,\cdots,d-1,\; x\in\mathbb{R}^{d},\quad \Delta_{d}\big(\eta(x^\alpha)\big)=N^{-1},\\
&\max_{x\in\mathbb{R}^{d}}\|J_{\eta}(x)\|_2\leq\frac{n}{n-1}\quad \mbox{and}\quad \max_{x\in\mathbb{R}^{d}}\|J_{\eta^{-1}}(x)\|_2\leq\frac{n}{n-1}.
\end{align*}
\end{lemma}
\begin{proof}
Let $f_j(x)={\rm diag}(\boldsymbol{0}_{d-1},1)\sigma({\rm diag}(\boldsymbol{0}_{d-1},x_{j}))\in\mathcal{F}$, $j=1,\cdots, d-1$,
and define the mapping
$\eta:=\psi^{f_{d-1}}_{n^{-(d-1)}}\circ\cdots\circ\psi^{f_2}_{n^{-2}}\circ \psi^{f_1}_{n^{-1}}\in\mathcal{A}_\mathcal{F}$. Then
there holds $\{\eta(x^\alpha)_d\}_{\alpha\in[n-1]^{d}}=\{i N^{-1}\}_{i=0}^{N-1}$.
In fact, by the definitions of $\psi^{f}_\tau$ and $f_j$, we have
\begin{equation}\label{eqn:perturb_d}
\eta(x^\alpha)_d=x^\alpha_d+\sum_{j=1}^{d-1}n^{-j}x^\alpha_j.
\end{equation}
Thus, $\eta$ maps the set $\{x^\alpha=\frac{\alpha}{n}\}_{\alpha\in[n-1]^{d}}$ to a set of $N$ points
with distinct values for the last component, and
$\Delta_{d}\big(\eta_{d}(x^{\alpha})\big)=N^{-1}.$
Meanwhile, for any $j=1,\cdots,d-1$, by the identity \eqref{eqn:Jac-form}, we have
\begin{align*}
J_{\psi^{f_j}_{n^{-j}}}(x)={\rm exp}\Big(\int_{0}^{n^{-j}} J_{f_j}(\psi^{f_j}_t(x)){\rm d}t\Big):={\rm exp}\big(n^{-j} A_j(x_{j})\big)= I+n^{-j}A_j(x_{j}),
\end{align*}
since the matrix $A_j(x_{j})\in\mathbb{R}^{d\times d}$ has only one nonzero
element
$[A_j(x_{j})]_{d,j}=\frac{{\rm d}\sigma (x_{j})}{{\rm d}x_j}:=a_j,$
which is also independent of $x_{d}$, and thus $A_j(x_{j})^k=0$ for any $k\geq2$. Furthermore, we have $|a_j|\leq 1$.
Also note that the components $x_{j}$, $j=1,\cdots,d-1$, stay unchanged under
any of the mappings $\psi^{f_j}_{n^{-j}}$. Thus, there holds
\begin{align*}
J_\eta(x)=&J_{\psi^{f_{d-1}}_{n^{-(d-1)}}\circ\cdots\circ \psi^{f_1}_{n^{-1}}}(x)=J_{\psi^{f_{d-1}}_{n^{-d+1}}}(x)\cdots J_{\psi^{f_{1}}_{n^{-1}}}(x)\\
=&(I+n^{-(d-1)}A_{d-1}(x_{d-1}))\cdots(I+n^{-1}A_1(x_{1}))=I+\sum_{j=1}^{d-1}n^{-j}A_{j}(x_{j}),
\end{align*}
in view of the following identity for any $i,j=1,\cdots,d-1$:
$A_{i}(x_{i})A_{j}(x_{j})=0$.
Thus, we can bound the spectral norm $\|J_{\eta}(x)\|_2$ of the Jacobian $J_\eta(x)$ by
\begin{align*}
\|J_{\eta}(x)\|_2\leq& 1+\Big\|\sum_{j=1}^{d-1}n^{-j}A_{j}(x_{j})\Big\|_2 \leq 1+ \bigg({\rm tr}\Big(\big(\sum_{j=1}^{d-1}n^{-j}A_{j}(x_{j})\big)^t\big(\sum_{j=1}^{d-1}n^{-j}A_{j}(x_{j})\big)\Big)\bigg)^\frac12\\
\leq& 1+\bigg({\sum_{j=1}^{d-1}n^{-2j}}\bigg)^\frac12\leq 1+\left(\frac{n^{-2}}{1-n^{-2}}\right)^\frac{1}{2} \leq 1+\frac{1}{n-1} =\frac{n}{n-1}.
\end{align*}
Since ${\rm det}(J_\eta(x))=1$, the mapping $\eta$ is invertible, and
\begin{align*}
\max_{x\in\mathbb{R}^{d}}\|J_{\eta^{-1}}(x)\|_2=&\max_{x\in\mathbb{R}^{d}}\|(J_{\eta}(x))^{-1}\|_2=\Big\|I-\sum_{j=1}^{d-1}n^{-j}A_{j}(x_{j})\Big\|_2
\leq \frac{n}{n-1}.
\end{align*}
This shows the last estimate and completes the proof of the lemma.
\end{proof}

\begin{remark}\label{rem:perturb_1}
The forward process of the mapping $\eta: x\to y$ in Lemma \ref{lem:perturb_1} is equivalent
to a coupling-based INN {\rm(}with one weight layer and one added layer,
i.e., identity mapping, and $n$ neurons{\rm)} and the inverse process $\eta^{-1}$ likewise
\cite{DinhSohlBengio:2017,ArdizzoneKruse:2019,TeshimaIshikawaSugiyama:2020}
\begin{align*}
  y_j = \left\{\begin{aligned}
    x_j, &\quad j=1,\ldots, d-1,\\
    x_{d}+\sum_{j=1}^{d-1}n^{-j}\sigma(x_{j}), &\quad j=d,
  \end{aligned}\right.
  \quad \mbox{and} \quad
  x_j = \left\{\begin{aligned}
    y_{j}, &\quad j=1,\cdots,d-1,\\
    y_{d}-\sum_{j=1}^{d-1}n^{-j}\sigma(x_{j}), &\quad j=d.
  \end{aligned}\right.
\end{align*}
\end{remark}

For the construction below, it is convenient to introduce the hypercube
\begin{equation*}
S_{x^{\alpha}}(s):= (x^{\alpha}_1-s,x^{\alpha}_1+s)\times\cdots\times(x^\alpha_{d}-s,x^{\alpha}_d+s),
\end{equation*}
which is a hypercube centered at $x^\alpha$, each side of length $2s$. Also we
decompose a vector $x\in\mathbb{R}^{d}$ into
$x=\begin{pmatrix}
x'\\
x_{d}
\end{pmatrix}$, where the subvector $x'\in\mathbb{R}^{d-1}$ contains the first $d-1$ components of
the vector $x$. Then we construct an NN, for any $\alpha\in[n-1]^{d}$, to transport
$(x^{\alpha})'$ to $(y^{\alpha})'(=\big(F(x^{\alpha})\big)')$, and keep the last component
$x_{d}$ for any $x\in\mathbb{R}^{d}$ and $(x)'$ for any $x\in (S_{x^{\alpha}}(\frac{1}{2N}))^c=
\mathbb{R}^{d}\setminus S_{x^{\alpha}}(\frac{1}{2N})$ unchanged. Note that the set
$\{S_{x^{\alpha}}(\frac{1}{2N})\}_{\alpha\in[n-1]^{d}}$ consists of a collection of disjoint
open hypercubes in $\mathbb{R}^{d}$.

\begin{lemma}\label{lem:trans_1}
Let $\{(x^{\alpha})_{d}\}_{\alpha\in[n-1]^{d}}$ be distinct with the smallest distance $\frac{1}{N}$.
Then for any $\alpha_1\in[n-1]^{d}$, there exists an invertible mapping $\varphi\in\mathcal{A}_\mathcal{F}$,
such that
\begin{align*}
\varphi(x^{\alpha_1})'=(y^{\alpha_1})',\quad
(\varphi(x))_{d}=&x_{d}, \; \forall x\in\mathbb{R}^{d} \quad\mbox{and}\quad
\varphi(x)=x , \; \forall x\in (S_{x^{\alpha_1}}(\tfrac{1}{2N}))^c.
\end{align*}
Furthermore, the following estimates hold
\begin{align*}
\max_{x\in\mathbb{R}^{d}}\|J_\varphi(x)\|_2\leq 1+6N\|(y^{\alpha_1}-x^{\alpha_1})'\|_2 \quad \mbox{and}
\quad \max_{x\in\mathbb{R}^{d}}\|J_{\varphi^{-1}}(x)\|_2\leq 1+6N\|(y^{\alpha_1}-x^{\alpha_1})'\|_2.
\end{align*}
\end{lemma}
\begin{proof}
First, we define an intermediate mapping
\begin{align*}
\ell(x)&=\sigma\big(-\sigma(\tfrac12
W^0x)+b^0\big)-\sigma\big(-\sigma(W^0
x)+b^0\big):=\ell^{(1)}(x) + \ell^{(2)}(x),
\end{align*}
with the weight matrix $W^0 =\begin{pmatrix}
O_{d\times (d-1)}&\boldsymbol{1}_d
\end{pmatrix}\in \mathbb{R}^{d\times d}$ and bias vector $b^0 = \boldsymbol{1}_d\in \mathbb{R}^d$.
The mapping $\ell(x)$ can be equivalently written as
$\ell(x)=\ell_0(x_{d})\boldsymbol{1}_d,$
where the function $\ell_0: \mathbb{R}\longrightarrow\mathbb{R}$ is defined by
\begin{equation}\label{eq:ell_0}
\ell_0(x)=\sigma\big(-\sigma(\tfrac12
x)+1\big)-\sigma\big(-\sigma(
x)+1\big).
\end{equation}
Note that $\ell_0(x)$ is a hat function supported on $x\in [0,2]$:
\begin{equation*}
  \ell_0(x) = \left\{\begin{array}{ll}
    0, & x\in (-\infty,0)\cup(2,+\infty),\\
    \frac{x}{2}, & x\in [0,1],\\
    1-\frac{x}{2}, & x\in [1,2].
  \end{array}\right.
\end{equation*}
Next let
$A={\rm diag}(\boldsymbol{0}_{d-1},2\Delta_{d}^{-1}(x^\alpha))\in\mathbb{R}^{d\times d}$ and
$b^{\alpha_1}=(\boldsymbol{0}_{d-1},1-2\Delta_{d}^{-1}(x^\alpha)x^{\alpha_1}_d)^t\in\mathbb{R}^d$.
Note that the minimal distance $\Delta_d(x^\alpha)=\min_{\alpha_2\neq\alpha_3}|x^{\alpha_2}_d-x^{\alpha_3}_d|=\frac{1}{N}$
is independent of $\alpha$, so is the matrix $A$. Then for any $\alpha_1\in[n-1]^{d}$, the following properties hold
\begin{align*}
(Ax^{\alpha_1}+b^{\alpha_1})_{d}=1 \quad\mbox{and}\quad
(Ax+b^{\alpha_1})_{d}=2\Delta_{d}^{-1}(x^\alpha)(x_{d}-x^{\alpha_1}_d)+1\notin (0,2), \quad \forall x\in (S_{x^{\alpha_1}}(\tfrac{1}{2N}))^c,
\end{align*}
since $|x_{d}-x^{\alpha_1}_d|\geq \frac{1}{2N}$ for any $x\in (S_{x^{\alpha_1}}(\tfrac{1}{2N}))^c$ implies $|2\Delta_{d}^{-1}
(x^\alpha)(x_{d}-x^{\alpha_1}_d)|\geq 1$. Further, it implies that for any $\alpha_1\in[n-1]^{d}$,
\begin{align*}
&\ell(Ax^{\alpha_1}+b^{\alpha_1})=\tfrac12\boldsymbol{1}_d\quad\mbox{and}\quad
\ell(Ax+b^{\alpha_1})=0,  \;\;\;\forall x\in (S_{x^{\alpha_1}}(\tfrac{1}{2N}))^c.
\end{align*}
Next, we define, for some fixed $\alpha_1\in[n-1]^{d}$, that
\begin{align*}
D_{\alpha_1}=&{\rm diag}(2(y_1^{\alpha_1}-x_1^{\alpha_1}),\cdots,2(y_{d-1}^{\alpha_1}-x_{d-1}^{\alpha_1}),0)\in\mathbb{R}^{d\times d}.
\end{align*}
Then we claim that $\varphi:=\psi_1^f$ with $f(x)=D_{\alpha_1}\ell(Ax+b^{\alpha_1})$ is the desired mapping
if $\varphi\in\mathcal{A}_\mathcal{F}$ and it is invertible. In fact, we can decompose $\ell$
into $\ell=\ell^{(1)}+\ell^{(2)}$, with $\ell^{(1)},\ell^{(2)}\in \mathcal{F}$. Note that these
two mappings actually depend only on the unchangeable value $x_{d}$ along the time. Let $f^{(i)}(x)=D_{\alpha_1}\ell^{(i)}(Ax+b^{\alpha_1})$ for $i=1,2$, then $\varphi=\psi_1^f=\psi_1^{f^{(1)}}\circ\psi_1^{f^{(2)}}\in\mathcal{A}_{\mathcal{F}}$.
By the definitions of $f$ and $\ell$, we have
\begin{align}
f(x)=&D_{\alpha_1}\ell(Ax+b^{\alpha_1})
=
2\Big[\sigma\Big(-\sigma\big[\tfrac12\big(2\Delta_{d}^{-1}(x^\alpha)(x_{d}-x_d^{\alpha_1})+1\big)\big]+1 \Big) \nonumber\\
&\quad-\sigma\Big(-\sigma\big[2\Delta_{d}^{-1}(x_\alpha) (x_{d}-x_d^{\alpha_1}) +1\big]+1\Big)\Big]\begin{pmatrix}
(y^{\alpha_1}-x^{\alpha_1})'\\
0
\end{pmatrix}
:=h_{\alpha_1}(x_{d})\begin{pmatrix}
(y^{\alpha_1}-x^{\alpha_1})'\\
0
\end{pmatrix}. \label{eq:h_1}
\end{align}
Then there holds
\begin{align*}
J_\varphi(x)=&I+\frac{{\rm d}h_{\alpha_1}(x_{d})}{{\rm d}x_{d}}\begin{pmatrix}
O_{(d-1)\times(d-1)}&(y^{\alpha_1}-x^{\alpha_1})'\\
O_{1\times(d-1)}&0
\end{pmatrix}:=I+H_{\alpha_1}(x).
\end{align*}
It follows from direct computation that
\begin{align*}
  \Big|\frac{{\rm d}h_{\alpha_1}(x_{d})}{{\rm d}x_{d}}\Big|\leq2 (\Delta_{d}^{-1}(x^\alpha)+2\Delta_{d}^{-1}
  (x^\alpha))= 6\Delta_{d}^{-1}(x^\alpha)=6N \quad \mbox{and}\quad {\rm det}(J_\varphi(x))=1.
\end{align*}
This directly implies that $\varphi\in\mathcal{A}_\mathcal{F}$ is invertible and further,
\begin{align*}
\|J_\varphi(x)\|_2\leq 1+\|H_{\alpha_1}(x)\|_2\leq 1+\sqrt{{\rm tr}\big((H_{\alpha_1}(x))^t
H_{\alpha_1}(x)\big)}\leq 1+6N\|(y^{\alpha_1}-x^{\alpha_1})'\|_2.
\end{align*}
Similarly,  we deduce
\begin{align*}
\max_{x\in\mathbb{R}^{d}}\|J_{\varphi^{-1}}(x)\|_2=\max_{x\in\mathbb{R}^{d}}\|(J_{\varphi}(x))^{-1}
\|_2=\max_{x\in\mathbb{R}^{d}}\|I-H_{\alpha_1}(x)\|_2\leq 1+6N\|(y^{\alpha_1}-x^{\alpha_1})'\|_2.
\end{align*}
This completes the proof of the lemma.
\end{proof}

\begin{remark}\label{rem:trans_1}
The forward process of the mapping $\varphi_\alpha: x\to y$ in Lemma \ref{lem:trans_1} is equivalent to the
following coupling-based NN {\rm(}with six weight layers and two added layers, i.e.,
identity mapping, and $d$ neurons on each layer{\rm)} with the mapping $h_\alpha:\mathbb{R}^d\to\mathbb{R}^d$
defined as \eqref{eq:h_1} and likewise the inverse process $\varphi^{-1}$:
\begin{align*}
\left\{\begin{aligned}
  y' &= x'+h_{\alpha}(x_{d})(y^{\alpha}-x^{\alpha})',\\
  y_{d} & = x_{d},\\
\end{aligned}\right.
\quad {\rm and}\quad \left\{\begin{aligned}
  x_{d}&= y_{d},\\
   x'&= y'-h_\alpha(x_{d})(y^{\alpha}-x^{\alpha})'.
\end{aligned}
\right.
\end{align*}
Equivalently, $y=x+D_{\alpha}\ell(Ax+b^{\alpha})$.
\end{remark}

So far we have already successfully transported the set $\{x^\alpha\}_{\alpha\in[n-1]^{d}}$
to $\{z^\alpha\}_{\alpha\in[n-1]^{d}}$ with $z^\alpha=\big((y^\alpha)'^t,\eta(x^{\alpha}_d)\big)^t$.
Following the ideas of Lemmas \ref{lem:perturb_1} and \ref{lem:trans_1}, we can further transport the
collection  $\{z^\alpha\}_{\alpha\in [n-1]^d}$ of points to $\{y^\alpha\}_{\alpha\in[n-1]^{d}}$. To
this end, we first perturb the set $\{z^\alpha\}_{\alpha\in[n-1]^{d}}$ by $\tilde{\eta}$ so that
$\Delta_j(\tilde{\eta}(z^\alpha))\geq\epsilon$ for some $j\neq d$ and $\epsilon>0$. This distinctness
at the $j$th coordinate is needed for performing the last step of the construction, and it can be
achieved as follows. If there exists some $j\neq d$ such that $\Delta_j(z^\alpha)>0$, we take
$\epsilon=\Delta_j(z^\alpha)$ and $\tilde{\eta}={\rm id}$. Otherwise, we construct the desired mapping
$\tilde{\eta}\in\mathcal{A}_\mathcal{F}$ in Lemma \ref{lem:perturb_3} below. The next lemma analyzes
the fundamental case of two points $z^1, z^2\in\mathbb{R}^{d}$.

\begin{lemma}\label{lem:perturb_2}
Let $z^1,z^2\in\mathbb{R}^{d}$, $z^1\neq z^2$ and $z^1_{1}=z^2_{1}$. Then for any $a\in\mathbb{R}$, $\epsilon:=a-z^1_1$ and $0<\Delta\leq \max_{j=1,\cdots,d}|z^1_{j}-z^2_{j}|$, there
exists an invertible mapping $\tilde\eta_0\in\mathcal{A}_\mathcal{F}$ such that
\begin{align*}
&(\tilde\eta_0(z^1))_{1}=a, \quad (\tilde\eta_0(z^1))_{j}=z^1_{j},\quad j=2,\cdots,d,\\
&\tilde\eta_0(z)=z, \; \forall z\in (S_{z^1}(\Delta))^c \quad \mbox{and}\quad \|\tilde\eta_0(z)-z\|_{\mathcal{L}^\infty(\mathbb{R}^{d})}\leq|\epsilon|.
\end{align*}
\end{lemma}
\begin{proof}
The construction of $\tilde\eta_0$  is inspired by the argument of \cite[Lemma 4.14]{LiLinShen:2019}. Since $z^1\neq z^2$, there exists an index
$j\neq1$ such that $z^1_{j}\neq z^2_{j}$, and we may assume $z^1_{2}\neq z^2_{2}$. Then for any $0<\Delta\leq|z^1_{2}-z^2_{2}|$, consider the following initial value problem
\begin{equation*}
\left\{\begin{aligned}
\frac{{\rm d}y}{{\rm d}t}&=g(
W_{22}y+e^2):=f(y),\quad t> 0,\\
y(0)&=z,
\end{aligned}\right.
\end{equation*}
with $W_{22}\in\mathbb{R}^{d\times d}$, $e^2\in \mathbb{R}^d$ and $g$ given by
\begin{align*}
  W_{22}&=\mathrm{diag}(0,\Delta^{-1},\boldsymbol{0}_{d-2})\in\mathbb{R}^{d\times d}, \quad e^2=(0,1-\Delta^{-1}z^1_2,\boldsymbol{0}_{d-2})^t\in\mathbb{R}^d,\\
   g(z)&:=\sigma \big(-\sigma(\tfrac12 W_{12}z)+e^1\big) -\sigma \big(-\sigma(W_{12}z)+e^1\big),
\end{align*}
where the matrix $W_{12}\in\mathbb{R}^{d\times d}$ has only one nonzero entry $[W_{12}]_{12}=1$ and the vector
$e^1\in \mathbb{R}^d$ has one nonzero entry $[e^1]_1=1$. By the definitions of $\Delta$ and $g$, we have
$f(z^1)=\tfrac{1}{2}e^1$ and $f(z)=0$, for any $z\in\mathbb{R}^d\mbox{ with } |z_{2}-z^1_{2}|\geq\Delta$,
since for $z\in\mathbb{R}^d$ with $|z_{2}-z^1_{2}|\geq\Delta$,
$(W_{22}z+e^2)_2 = \Delta^{-1}(z_2-z^1_2) + 1 \not\in (0,2)$,
and $g$ is supported on the interval $[0,2]$.
Thus, for any $z\in\mathbb{R}^d$ satisfying $|z_{2}-z^1_{2}|\geq\Delta$, $z^1_{j}$ for
any $j=2,\cdots,d$ stay unchanged. Moreover, since $f(z^1)=\frac{e^1}{2}$,
for any $\epsilon:=a-z^1_1\geq 0$, when $t=2\epsilon$, we have
\begin{align*}
&[\psi^f_t(z^1)]_1-z^1_{1}=[\psi^f_t(z^1)]_1-z^2_{1}=\epsilon,\\
&\|\psi^f_t(z)-z\|_{\mathcal{L}^\infty(\mathbb{R}^{d})}\leq |[\psi^f_t(z^1)]_1-z^1_{1}|=\epsilon,
\end{align*}
since by assumption $z_1^1=z_1^2$, and $|f(z)|\leq \frac12$. Similarly, for any $\epsilon=a-z^1_1 < 0$, when $t=-2\epsilon$, we have
\begin{align*}
&[\psi^{-f}_t(z^1)]_1-z^1_{1}=[\psi^{-f}_t(z^1)]_1-z^2_{1}=\epsilon,\\
&\|\psi^{-f}_t(z)-z\|_{\mathcal{L}^\infty(\mathbb{R}^{d})}\leq |[\psi^{-f}_t(z^1)]_1-z^1_{1}|=-\epsilon.
\end{align*}
Let $$\tilde\eta_0(z):=\left\{\begin{array}{ll}
\psi^f_{2\epsilon}(z),   &\mbox{if } \epsilon \geq 0,\\
\psi^{-f}_{-2\epsilon}(z),   & {\rm otherwise}.
\end{array}\right.$$
Then repeating the argument of Lemma
\ref{lem:trans_1} shows $\tilde\eta_0\in\mathcal{A}_\mathcal{F}$ by decomposing $g=g^{(1)}
+g^{(2)}$, where $g^{(1)},g^{(2)}\in \mathcal{F}$ depend only on the unchangeable value. Note that
\begin{align*}
J_{\tilde{\eta}_0}(z)
={\rm exp}\big(2\epsilon \tfrac{{\rm d}h(z_{2})}{{\rm d}z_{2}}W_{12}\big)
:={\rm exp}\big(2\epsilon A(z_{2})\big)=I+2\epsilon A(z_{2}),
\end{align*}
with $A(z_2)=\frac{{\rm d}h(z_{2})}{{\rm d}z_{2}}W_{12}$, where
\begin{align*}
h(z_{2})=\left\{\begin{array}{ll}
\frac{1}{2}\big(\Delta^{-1}(z_{2}-z^1_{2})+1\big),   &\mbox{if } \Delta^{-1}(z_{2}-z^1_{2})+1\in(0,1],\\
1-\frac{1}{2}\big(\Delta^{-1}(z_{2}-z^1_{2})+1\big),   &\mbox{if } \Delta^{-1}(z_{2}-z^1_{2})+1\in(1,2],\\
0,                  & {\rm otherwise}.
\end{array}\right.
\end{align*}
Thus, ${\rm det}(J_{\tilde{\eta}_0}(z))= 1$, and the mapping
$\tilde{\eta}$ is invertible and the inverse $(J_{\tilde{\eta}_0}(z))^{-1}$ of $J_{\tilde{\eta}_0}(z)$ is given by
$(J_{\tilde{\eta}_0}(z))^{-1}=I-2\epsilon A(z_{2})$.
Then with the choice $\tilde{\eta}_0$, the asserted properties in the lemma follow.
\end{proof}

\begin{remark}\label{rem:perturb_2}
The forward process of the mapping $\eta: x\to y$ in Lemma \ref{lem:perturb_2} is
equivalent to a coupling-based NN {\rm(}with four weight layers,
two added layers, i.e., identity mapping and $d$ neurons on each layer{\rm)}, and likewise the inverse process $\eta^{-1}$:
\begin{align*}
y_j =\left\{\begin{aligned}
   x_j+2\epsilon h(x_{2}), &\quad j=1,\\
   x_{j}, &\quad j=2,\cdots,d,
   \end{aligned}\right.
   \quad \mbox{and} \quad x_j =\left\{\begin{aligned}
   y_{j}, &\quad j=2,\cdots,d,\\
   y_j+2\epsilon h(x_{2}), &\quad j=1.
   \end{aligned}\right.
\end{align*}
\end{remark}

The next lemma generalizes the transformation in Lemma \ref{lem:perturb_2} to multiple points.
\begin{lemma}\label{lem:perturb_3}
Let the set of points $\{z^\alpha=\big((y^\alpha)'^t,x^{\alpha}_d\big)^t\}_{\alpha\in[n-1]^{d}}$ be such that $\Delta_{d}
(z^\alpha)=N^{-1}$. Then for any $\epsilon\in(0,1)$, there exists an INN $\tilde{\eta}
\in\mathcal{A}_\mathcal{F}$ and an index $j_0< d$ such that
\begin{align}\label{eq:prop}
\Delta_{j_0}\big(\tilde{\eta}(z^\alpha)\big)\geq N^{-1}\epsilon,\quad \max_{z\in\mathbb{R}^{d}}
|(\tilde\eta(z)-z)_{j_0}|<\epsilon \quad \mbox{and}\quad (\tilde{\eta}(z^\alpha))_{j}=z^{\alpha}_j, \; \forall j\neq j_0.
\end{align}
Moreover, the following estimates hold
\begin{align*}
\max_{x\in\mathbb{R}^{d}}\|J_{\tilde\eta}(x)\|_2\leq 1+2N\epsilon \quad \mbox{and}\quad \max_{x\in\mathbb{R}^{d}}\|J_{{\tilde\eta}^{-1}}(x)\|_2\leq 1+2N\epsilon.
\end{align*}
\end{lemma}
\begin{proof}
For any $\epsilon\in(0,1)$ and any fixed index $j_0<d$, there exists a sequence of disjoint closed intervals $\{I_j\}_{j=1}^J$, each of length $\epsilon$, with $\mathbb{N}\ni J\leq N$, such that the inclusion $\{z^\alpha_{j_0}\}_{\alpha\in[n-1]^{d}} \subset \bigcup_{j=1}^J I_j$ holds. Consider an interval $I_j:=[b_l,b_r]$ and assume that there are $K$ components of $\{z^\alpha_{j_0}\}_{\alpha\in[n-1]^{d}}$, denoted by $\{z^{\alpha_i}_{j_0}\}_{i=1}^K$, lie in $I_j$. Clearly we have $K\leq N$. Then there exists a set $\{a^{\alpha_i}\}_{i=1}^{K}\subset I_j$ such that
\begin{align*}
\min\big(\min_{1\leq i<j\leq K}|a^{\alpha_i}-a^{\alpha_j}|,\min_{1\leq i\leq K}|a^{\alpha_i}-b_r|\big)= K^{-1}\epsilon\geq N^{-1}\epsilon.
\end{align*}
By Lemma \ref{lem:perturb_2}, for any index $1\leq i\leq K$, we can construct an invertible mapping $\tilde\eta_{\alpha_i}\in\mathcal{A}_\mathcal{F}$ with $\Delta=
(2N)^{-1}$ such that
\begin{align*}
(\tilde\eta_{\alpha_i}(z^{\alpha_i}))_{j_0}=a^{\alpha_i}, \;\;\|\tilde\eta_{\alpha_i}(z)-z\|_{\mathcal{L}^\infty(\mathbb{R}^{d})}= |a^{\alpha_i}-z^{\alpha_i}_{j_0}|\leq\epsilon,
\end{align*}
and $\tilde\eta_{\alpha_i}(z)=z$ for any $z\in (S_{z^{\alpha_i}}(\frac{1}{2N}))^c$.
Similarly, there exist a set $\{a^{\alpha_i}\}_{i=1}^{N}\subset \bigcup_{j=1}^J I_j$ and a sequence of invertible mapping $\{\tilde\eta_{\alpha_i}\}_{i=1}^N$ such that
\begin{align*}
\min_{1\leq i<j\leq N}|a^{\alpha_i}-a^{\alpha_j}|\geq N^{-1}\epsilon,\;\;
(\tilde\eta_{\alpha_i}(z^{\alpha_i}))_{j_0}=a^{\alpha_i}, \;\;\|\tilde\eta_{\alpha_i}(z)-z\|_{\mathcal{L}^\infty(\mathbb{R}^{d})}= |a^{\alpha_i}-z^{\alpha_i}_{j_0}|\leq\epsilon,
\end{align*}
and $\tilde\eta_{\alpha_i}(z)=z$ for any $z\in (S_{z^{\alpha_i}}(\frac{1}{2N}))^c$.
Furthermore, we have
$J_{\tilde{\eta}_{\alpha_i}}(z)=I+2\epsilon A_{\alpha_i}(z_{d})$,
where the matrix $A_{\alpha_i}(z_{d})\in\mathbb{R}^{d\times d}$ has only one nonzero entry at $[A_{\alpha_i}(z_{d})
]_{j_0,d}$ with $|[A_{\alpha_i}(z_{d})]_{j_0,d}|\leq \frac12 \Delta^{-1}$.
Next we define
$\tilde\eta:= \tilde\eta_{\alpha_N}\circ\cdots\circ\tilde\eta_{\alpha_1}\in\mathcal{A}_\mathcal{F}$.
It can be shown that $\tilde\eta$ satisfies the properties in \eqref{eq:prop}, since
\begin{equation*}
\Delta_{j_0}\big(\tilde{\eta}(z^\alpha)\big)=\min_{1\leq i<j\leq N}|a^{\alpha_i}-a^{\alpha_j}|\geq N^{-1}\epsilon.
\end{equation*}
Further, since the supports of $ J_{\tilde \eta_{\alpha_i}}$ are disjoint from each other, there hold
\begin{align*}
\|J_{\tilde\eta}(x)\|_2\leq & \max_{i=1,\cdots,N} \|J_{\tilde\eta_{\alpha_i}}(x)\|_2
\leq\max_{i=1,\cdots,N}\|I+2\epsilon A_{\alpha_i}(z_{d})\|_2\leq  1+2\epsilon \max_{i=1,\cdots,N}\| A_{\alpha_i}(z_{d})\|_2\\
\leq& 1+2\epsilon \max_{i=1,\cdots,N}\sqrt{{\rm tr }( A_{\alpha_i}(z_{d}))^t(A_{\alpha_i}(z_{d})\big)}
\leq 1+\epsilon \Delta^{-1}=1+2N\epsilon,\\
\|(J_{\tilde\eta}(x))^{-1}\|_2\leq&\max_{i=1,\cdots,N}\|I-2\epsilon A_{\alpha_i}(z_{d})\|_2\leq1+2N\epsilon.
\end{align*}
This completes the proof of the lemma.
\end{proof}

Under the condition $\Delta_{j_0}\big(\tilde{\eta}(z^\alpha)\big)\geq N^{-1}\epsilon$,
we can now construct an INN $\tilde\varphi$ for any $\alpha_1\in[n-1]^{d}$, to transport $z_d^{\alpha_1}$ to $y_d^{\alpha_1}=(\tilde\varphi(z^{\alpha_1}))_{d}$, and keep $(z)'$ for any $z\in\mathbb{R}^{d}$ and $z_{d}$
for any $z\in (S_{z^{\alpha_1}}(\frac{\epsilon}{2N}))^c$ unchanged.

\begin{lemma}\label{lem:trans_2}
Let $\{z^{\alpha}_{j_0}\}_{\alpha\in[n-1]^{d}}$ be distinct with the smallest distance equal or larger than $\Delta_{j_0}=\frac{
\epsilon}{N}$ for some $j_0\neq d$ and $\epsilon>0$. Then for any $\alpha_1\in[n-1]^{d}$, there
exists an INN $\tilde\varphi\in\mathcal{A}_\mathcal{F}$ such that
\begin{align*}
&(\tilde\varphi(z^{\alpha_1}))_{d}=y_d^{\alpha_1},
\quad (\tilde\varphi(z))_{j}=z_{j}, \; \forall z\in\mathbb{R}^{d},\; j\neq d \quad\mbox{and}\quad
\tilde\varphi(z)=z , \; \forall z\in (S_{z^{\alpha_1}}(\tfrac{\epsilon}{2N}))^c.
\end{align*}
Furthermore, the following estimates hold
\begin{align*}
\max_{z\in\mathbb{R}^{d}}\|J_{\tilde\varphi}(z)\|_2\leq 1+6N\epsilon^{-1}|y_d^{\alpha_1}-z_d^{\alpha_1}| \quad \mbox{and}\quad \max_{x\in\mathbb{R}^{d}}\|J_{{\tilde\varphi}^{-1}}(x)\|_2\leq 1+6N\epsilon^{-1}|y_d^{\alpha_1}-z_d^{\alpha_1}|.
\end{align*}
\end{lemma}
\begin{proof}
We define an intermediate map $
\ell(z)=\sigma\big(-\sigma\big(\tfrac12Pz\big)+e^1\big)-\sigma\big(-\sigma\big(Pz\big)+e^1\big)$,
where $P\in\mathbb{R}^{d\times d}$ has only one nonzero entry $[P]_{d,j_0}=1$ and $e^1=\boldsymbol{1}_d\in\mathbb{R}^d$.
Let $A\in\mathbb{R}^{d\times d}$ with only one nonzero entry $[A]_{j_0,j_0}=2\Delta_{j_0}^{-1}$ and $b^{\alpha_1}\in
\mathbb{R}^{d}$ with only one nonzero entry $b^{\alpha_1}_{j_0}=1-2\Delta_{j_0}^{-1}z_{j_0}^{\alpha_1}$.
Then there hold
\begin{align*}
&(Az^{\alpha_1}+b^{\alpha_1})_{j_0}=1, \quad\; \;\forall \alpha_1\in[n-1]^{d},\\
&(Az+b^{\alpha_1})_{j_0}\notin (0,2), \quad \forall z\in (S_{z^{\alpha_1}}(\tfrac{\epsilon}{2N}))^c.
\end{align*}
These two properties together imply
\begin{align}
[\ell(Az^{\alpha_1}+b^{\alpha_1})]&=(\boldsymbol{0}_{d-1},\tfrac12)^t, \quad \forall \alpha_1\in[n-1]^{d},\\
\ell(Az+b^{\alpha_1})&=0, \,\qquad\qquad\qquad\forall z\in (S_{z^{\alpha_1}}(\tfrac{\epsilon}{2N}))^c.
\end{align}
For a fixed index $\alpha_1\in[n-1]^{d}$, define
$D_{\alpha_1}={\rm diag}(\boldsymbol{0}_{d-1},2(y_d^{\alpha_1}-z_d^{\alpha_1}))\in\mathbb{R}^{d\times d}.$
Then $\tilde\varphi:=\psi_1^f\in\mathcal{A}_{\mathcal{F}}$ with $f(z)=D_{\alpha_1}\ell(Az+b_{\alpha_1})$ is the desired mapping.
Indeed, by the definition of $f$, we have
\begin{align*}
f(x)=&
2\Big[\sigma\Big(-\sigma\big[\tfrac12\big(2\Delta_{j_0}^{-1}(z_{j_0}-z_{j_0}^{\alpha_1})+1\big)\big]+1 \Big)\\
&\quad-\sigma\Big(-\sigma\big[2\Delta_{j_0}^{-1} (z_{j_0}-z_{j_0}^{\alpha_1}) +1\big]+1\Big)\Big]\begin{pmatrix}\boldsymbol{0}_{d-1}\\
y_d^{\alpha_1}-z_d^{\alpha_1}
\end{pmatrix}
:=h_{\alpha_1}(z_{j_0})\begin{pmatrix}\boldsymbol{0}_{d-1}\\
y_d^{\alpha_1}-z_d^{\alpha_1}
\end{pmatrix}.
\end{align*}
Then there holds $J_{\tilde\varphi}(z)=I+H_{\alpha_1}(z)$,
where the matrix $H_{\alpha_1}\in\mathbb{R}^{d\times d}$ has only one nonzero entry
\begin{align*}
[H_{\alpha_1}]_{d,j_0}=(y_d^{\alpha_1}-z_d^{\alpha_1})\frac{{\rm d}h_{\alpha_1}(z_{j_0})}{{\rm d}z_{j_0}},\quad \mbox{with }
\Big|\frac{{\rm d}h_{\alpha_1}(z_{j_0})}{{\rm d}z_{j_0}}\Big|\leq 6\Delta_{j_0}^{-1}=\frac{6N}{\epsilon}.
\end{align*}
Consequently,
\begin{align*}
\|J_{\tilde\varphi}(z)\|_2&\leq 1+\|H_{\alpha_1}(z)\|_2\leq 1+6N\epsilon^{-1}|y_d^{\alpha_1}-z_d^{\alpha_1}|,\\
\max_{z\in\mathbb{R}^{d}}\|J_{{\tilde\varphi}^{-1}}(z)\|_2&=\max_{z\in\mathbb{R}^{d}}\|(J_{\tilde\varphi}(z))^{-1}
\|_2=\max_{z\in\mathbb{R}^{d}}\|I-H_{\alpha_1}(z)\|_2\leq 1+6N\epsilon^{-1}|y_d^{\alpha_1}-z_d^{\alpha_1}|.
\end{align*}
This completes the proof of the lemma.
\end{proof}
\begin{remark}\label{rem:trans_2}
Similar to Remark \ref{rem:trans_1}, the forward process of the mapping $\tilde{\varphi}_\alpha$ in Lemma \ref{lem:trans_2}
amounts to a coupling-based INN with six weight layers and two added layers, and $d$ neurons on each layer.
\end{remark}

Last we construct an INN $\tilde{F}_{\rm nn}$ to approximate the map $F$. This is the main result of
this section.
\begin{theorem}\label{thm:F_NN}
Given a finite collection of evaluations $\{(x^{\alpha} , y^{\alpha})\}_{\alpha\in[n-1]^{d}}$, for any tolerance $\epsilon>0$, there exists a bi-Lipschitz continuous invertible mapping $\tilde{F}_{\rm nn}\in\mathcal{A}_\mathcal{F}$,
such that $\|\tilde{F}_{\rm nn}(x^{\alpha})-y^{\alpha}\|_2<\epsilon$ for any $\alpha\in[n-1]^{d}$, and the Lipschitz
constants of $\tilde{F}_{\rm nn}$ and $\tilde{F}_{\rm nn}^{-1}$ are bounded by
\begin{align*}
\Lip{\tilde{F}_{\rm nn}}\leq\frac{n}{n-1}\big(1+6Nc\big)(1+2N\epsilon)\big(1+6N\epsilon^{-1}c\big), \\
\Lip{\tilde{F}_{\rm nn}^{-1}}\leq\frac{n}{n-1}\big(1+6Nc\big)(1+2N\epsilon)\big(1+6N\epsilon^{-1}c\big),
\end{align*}
with $c:=\max_{i=1,\cdots N}\|y^{\alpha_i}-x^{\alpha_i}\|_2+n^{-1}$.
\end{theorem}
\begin{proof}
By Lemma \ref{lem:perturb_1}, there exists an INN $\eta\in \mathcal{A}_\mathcal{F}$ such
that $\Delta_{d}\big(\eta(x^\alpha)\big)=N^{-1}$ for all $\alpha\in [n-1]^d$ and
\begin{align*}
 \max_{x\in\mathbb{R}^{d}}\|J_{\eta}(x)\|_2\leq\frac{n}{n-1}\quad \mbox{and}\quad \max_{x\in\mathbb{R}^{d}}\|J_{\eta^{-1}}(x)\|_2\leq\frac{n}{n-1}.
\end{align*}
Then, by Lemma \ref{lem:trans_1}, we can construct a composite INN $\varphi^N:=\varphi_{\alpha_N}\circ\cdots\circ
\varphi_{\alpha_1}$, with $\big(\varphi^j(\eta(x^{\alpha_i}))\big)'=\big(y^{\alpha_i}\big)'$ for any $1\leq i\leq j\leq N$.
Note that each $\varphi_{\alpha_i}$, $i=1,\cdots,N$, is the identity mapping in  $(S_{\eta(x^{\alpha_i})}
(\frac{1}{2N}))^c$, and $\{S_{\eta(x^{\alpha_i})}(\frac{1}{2N})\}_{i=1}^N$ are disjoint from each other. Thus, we have
\begin{align*}
\max_{x\in\mathbb{R}^{d}}\|J_{\varphi^N}(x)\|_2\leq\max_{i=1,\cdots N,x\in\mathbb{R}^{d}} \|J_{\varphi_{\alpha_i}}(x)\|_2\leq 1+6N\max_{i=1,\cdots N}\|(y^{\alpha_i}-x^{\alpha_i})'\|_2.
\end{align*}
Similarly, there holds
\begin{equation*}
  \max_{x\in\mathbb{R}^{d}}\|J_{(\varphi^N)^{-1}}(x)\|_2\leq 1+6N\max_{i=1,\cdots N}\|(y^{\alpha_i}-x^{\alpha_i})'\|_2.
\end{equation*}
Then by Lemma \ref{lem:perturb_3}, we can construct $\tilde\eta\in\mathcal{A}_\mathcal{F}$ such that
\begin{align*}
\Delta_{j_0}\big(\tilde{\eta}\circ \varphi^N\circ \eta(x^\alpha)\big)\geq N^{-1}\epsilon,\quad \max_{x\in\mathbb{R}^{d}}\|J_{\tilde\eta}(x)\|_2\leq 1+2N\epsilon \quad \mbox{and}\quad \max_{x\in\mathbb{R}^{d}}
\|J_{{\tilde\eta}^{-1}}(x)\|_2\leq 1+2N\epsilon.
\end{align*}
Then repeating the argument with  Lemma \ref{lem:trans_2}, there exists an INN $\tilde{\varphi}^N:=
\tilde{\varphi}_{\alpha_N}\circ\cdots\circ\tilde{\varphi}_{\alpha_1}\in\mathcal{A}_\mathcal{F}$ such that
$(\tilde{\varphi}^N\circ\tilde{\eta}\circ \varphi^N \circ \eta(x^{\alpha_i})))_{d}=(\tilde\varphi^j(z^{\alpha_i}))_{d}=y^{\alpha_i}_d$ for any $1\leq i\leq j\leq N$ and
\begin{align*}
   \max_{z\in\mathbb{R}^{d}}\|J_{\tilde\varphi}(z)\|_2\leq& 1+6N\epsilon^{-1}\max_{i=1,\cdots N}
   |y_d^{\alpha_i}-z_d^{\alpha_i}|,\\
   \max_{z\in\mathbb{R}^{d}}\|J_{{\tilde\varphi}^{-1}}(z)\|_2\leq
   & 1+6N\epsilon^{-1}\max_{i=1,\cdots N}|y_d^{\alpha_i}-x_d^{\alpha_i}|.
\end{align*}
By \eqref{eqn:perturb_d}, we can bound $|y_d^{\alpha_i}-z_d^{\alpha_i}|$ by
\begin{equation*}
|y_d^{\alpha_i}-z_d^{\alpha_i}|\leq|y_d^{\alpha_i}-x_d^{\alpha_i}|+|z^{\alpha_i}_d-x^{\alpha_i}_d|\leq |y_d^{\alpha_i}-x_d^{\alpha_i}|+\sum_{j=1}^{d-1}n^{-j}x^{\alpha_i}_j\leq|y_d^{\alpha_i}-x_d^{\alpha_i}|+n^{-1}.
\end{equation*}
Finally, let $\tilde{F}_{\rm nn}=\tilde{\varphi}^N\circ\tilde{\eta}\circ \varphi^N \circ \eta\in
\mathcal{A}_\mathcal{F}$. This is the desired mapping. With $c:=\max_{i=1,\cdots N}\|y^{\alpha_i}-x^{\alpha_i}\|_2+n^{-1}$, the preceding estimates together imply
\begin{align*}
\Lip{\tilde{F}_{\rm nn}}=&\max_{z\in\mathbb{R}^{d}}\|J_{\tilde{F}_{\rm nn}}(z)\|
\leq \frac{n}{n-1}\big(1+6Nc\big)(1+2N\epsilon)\big(1+6N\epsilon^{-1}c\big),\\
\Lip{{\tilde{F}_{\rm nn}}^{-1}}=&\max_{z\in\mathbb{R}^{d}}\|J_{{\tilde{F}_{\rm nn}}^{-1}}(z)\|
\leq \frac{n}{n-1}\big(1+6Nc\big)(1+2N\epsilon)\big(1+6N\epsilon^{-1}c\big).
\end{align*}
This completes the proof of the theorem.
\end{proof}

\begin{remark}\label{rem:lip_F_NN}
By the construction in Theorem \ref{thm:F_NN} and the analysis in Remarks
\ref{rem:perturb_1}-\ref{rem:trans_2}, the INN $\tilde{F}_{\rm nn}$ is
a coupling-based INN with $1+16N$ weight layers and $1+6N$
added layers, and $d$ neurons on each layer.
\end{remark}

\subsection{Error estimation}
Now we bound the error of the INN approximation
$F_{\rm nn}:=\tilde{F}_{\nn}\circ H^r_{\nn}$, where $\tilde{F}_{\nn}$ (cf. Theorem \ref{thm:F_NN}) has $1+16N$ weight layers, $1+6N$
added layers and each layer having $d$ neurons, and $H^r_{\nn}(x):=(h^r_{\nn}(x_{1}),\cdots,h^r_{\nn}
(x_{d}))$ (cf. Section \ref{subsubsec:H^r}) with a two-layer NN $h^r_{\nn}$ of $2(n+1)$ neurons.
\begin{theorem}
For any $\epsilon>0$, define the bi-Lipschitz INNs
$F_{\nn}:=\tilde{F}_{\nn}\circ H^r_{\nn}, F_{\nn}^{-1}:=H^r_{\nn}\circ \tilde{F}_{\nn}^{-1}\in\mathcal{A}_\mathcal{F}$,
where the INNs $\tilde{F}_{\nn}$, $\tilde{F}_{\nn}^{-1}$ and $H^r_{\nn}$ are defined in Theorem \ref{thm:F_NN}. Then the following error estimates hold
\begin{align*}
\|F_{\nn}-F\|^2_{\mathcal{L}^2(K)}
\leq&2\big((\Lip{\tilde{F}_{\nn}}+\Lip{F})(1-r)\sqrt{d}n^{-1}+\epsilon\big)^2+2\Lip{F}^2(2r-1)^2dn^{-2}\\
&+2\big(\tfrac12(\Lip{\tilde{F}_{\nn}}+\Lip{F})\sqrt{d}n^{-1}+\epsilon\big)^2(1-r^d),\\
\|F^{-1}_{\nn}-F^{-1}\|^2_{\mathcal{L}^2(F(K))}\leq&
2\Lip{F}^{d}\Big[\Big(2+\big(\Lip{{\tilde{F}_{\nn}}^{-1}}+\Lip{F^{-1}}\big)\Lip{F}(1-r)\Big)\sqrt{d}n^{-1}+\Lip{F^{-1}}\epsilon\Big]^2\\
&+2\Lip{F}^d (2r-1)^2 d n^{-2}+8\Lip{F}^{d}(1-r^{d}).
\end{align*}
\end{theorem}
\begin{proof}
By the triangle inequality, we can decompose the total error $\|F_{\nn}-F\|^2_{\mathcal{L}^2(K)}$
between the INN $F_{\nn}$ and the bi-Lipschitz map $F$ into two parts
\begin{align*}
\|F_{\nn}-F\|^2_{\mathcal{L}^2(K)}\leq&2 \|\tilde{F}_{\nn}\circ H^r_{\nn}-F\circ H^r\|^2_{\mathcal{L}^2(K)}+2\|F\circ H^r- F\|^2_{\mathcal{L}^2(K)}\\
=&2\|\tilde{F}_{\nn}\circ H^r-F\circ H^r\|^2_{\mathcal{L}^2(K)}+2\|F\circ H^r- F\|^2_{\mathcal{L}^2(K)},
\end{align*}
since $H^r=H^r_{\nn}$, cf. the construction in Section \ref{subsubsec:H^r}. Then with
the union of the (disjoint) hypercubes
\begin{equation*}
K^r:=\bigcup_{\alpha\in[n-1]^{d}}[x^\alpha_{1},x^\alpha_{1}+\tfrac{r}{n}]\times\cdots\times[x^\alpha_{d},x^\alpha_{d}+\tfrac{r}{n}],
\end{equation*}
and using the argument of Lemma \ref{lem:H^r}, since $|K\setminus K^r|=1-r^d$, we deduce
\begin{align}
\|\tilde{F}_{\nn}\circ H^r-F\circ H^r\|^2_{\mathcal{L}^2(K)}
=&\|\tilde{F}_{\nn}\circ H^r-F\circ H^r\|^2_{\mathcal{L}^2(K^r)}+\|\tilde{F}_{\nn}\circ H^r-F\circ H^r\|^2_{\mathcal{L}^2(K\setminus K^r)}\nonumber\\
\leq&|K^r|\|\tilde{F}_{\nn}\circ H^r-F\circ H^r\|^2_{\mathcal{L}^\infty(K^r)}+|K\setminus K^r|\|\tilde{F}_{\nn}-F\|^2_{\mathcal{L}^\infty(K)}\nonumber\\
\leq&\|\tilde{F}_{\nn}\circ H^r-F\circ H^r\|^2_{\mathcal{L}^\infty(K^r)}+(1-r^d)\|\tilde{F}_{\nn}-F\|^2_{\mathcal{L}^\infty(K)}.\label{ineq:key}
\end{align}
For any index $\alpha\in [n-1]^{d}$, let $\xi^\alpha:=\tilde{F}_{\nn}(x^\alpha)-F(x^\alpha)$, i.e., the approximation error at the point $x^\alpha$. Then there holds $\|\xi^\alpha\|_2\leq \epsilon$ by the construction
of $\tilde{F}_{\nn}$, cf. Lemmas \ref{lem:trans_1} and \ref{lem:perturb_3}. Note that
\begin{align*}
&\|\tilde{F}_{\nn}\circ H^r-F\circ H^r\|_{\mathcal{L}^\infty(K^r)}\\
=&\|\tilde{F}_{\nn}\circ(H^r-\lim_{r\to 1^-}H^r)+(\tilde{F}_{\nn}-F)\circ\lim_{r\to 1^-}H^r+F\circ(\lim_{r\to 1^-}H^r- H^r)\|_{\mathcal{L}^\infty(K^r)}\\
\leq&\|\tilde{F}_{\nn}\circ(H^r-\lim_{r\to 1^-}H^r)\|_{\mathcal{L}^\infty(K^r)}+\|(\tilde{F}_{\nn}-F)\circ\lim_{r\to 1^-}H^r\|_{\mathcal{L}^\infty(K^r)}\\
 &+\|F\circ(\lim_{r\to 1^-}H^r- H^r)\|_{\mathcal{L}^\infty(K^r)}\\
\leq&\omega_{\tilde{F}_{\nn}}((1-r)\sqrt{d}n^{-1})+\max_{\alpha\in[n-1]^{d}}\|\xi^\alpha\|_2+\omega_{F}((1-r)\sqrt{d}n^{-1})\\
\leq& \Lip{\tilde{F}_{\nn}}(1-r)\sqrt{d}n^{-1}+\epsilon+\Lip{F}(1-r)\sqrt{d}n^{-1}.
\end{align*}
Meanwhile, we have the following estimate
\begin{align*}
\|\tilde{F}_{\nn}-F\|_{\mathcal{L}^\infty(K)}=&\min_{\alpha\in[n-1]^{d}}\|\tilde{F}_{\nn}(x)-\tilde{F}_{\nn}(x^\alpha) +\tilde{F}_{\nn}(x^\alpha)-F(x^\alpha)+F(x^\alpha)-F(x)\|_{\mathcal{L}^\infty(K)}\\
\leq&\min_{\alpha\in[n-1]^{d}}\big(\|\tilde{F}_{\nn}(x)-\tilde{F}_{\nn}(x^\alpha)\|_{\mathcal{L}^\infty(K)}+\|\xi^\alpha\|_2 +\|F(x^\alpha)-F(x)\|_{\mathcal{L}^\infty(K)}\big)\\
\leq&\omega_{\tilde{F}_{\nn}}(\sqrt{d}(2n)^{-1})+\epsilon+\omega_{F}(\sqrt{d}(2n)^{-1})
\leq \Lip{\tilde{F}_{\nn}}\sqrt{d}(2n)^{-1}+\epsilon+\Lip{F}\sqrt{d}(2n)^{-1}.
\end{align*}
Consequently, we obtain
\begin{align*}
\|\tilde{F}_{\nn}\circ H^r-F\circ H^r\|^2_{\mathcal{L}^2(K)}
\leq&\big((\Lip{\tilde{F}_{\nn}}+\Lip{F})(1-r)\sqrt{d}n^{-1}+\epsilon\big)^2+\big(\tfrac12(\Lip{\tilde{F}_{\nn}}+\Lip{F})\sqrt{d}n^{-1}+\epsilon\big)^2(1-r^d).
\end{align*}
Further, by Lemma \ref{lem:H^r}, we have
\begin{align*}
\|F\circ H^r-F\|^2_{\mathcal{L}^2(K)}\leq \Lip{F}^2(2r-1)^2dn^{-2}.
\end{align*}
Similarly, we derive
\begin{align*}
\|F_{\nn}^{-1}-F^{-1}\|^2_{\mathcal{L}^2(F(K))}
\leq&2\|H^r\circ{\tilde{F}_{\nn}}^{-1} -H^r\circ F^{-1}\|^2_{\mathcal{L}^2(F(K))}+2\|H^r\circ F^{-1}- F^{-1}\|^2_{\mathcal{L}^2(F(K))}.
\end{align*}
It follows from the estimate
\begin{align*}
\|{\tilde{F}_{\nn}}^{-1}\big(\tilde{F}_{\nn}(x^\alpha)\big)-F^{-1}\big(\tilde{F}_{\nn}(x^\alpha)\big)\|_2 =\|x^{\alpha}-F^{-1}\big(F(x^\alpha)+\xi^\alpha\big)\|_2\leq \Lip{F^{-1}}\epsilon
\end{align*}
that
\begin{align*}
\|{\tilde{F}_{\nn}}^{-1}-F^{-1}\|_{\mathcal{L}^\infty(F(K^r))}\leq&\omega_{{\tilde{F}_{\nn}}^{-1}}(\Lip{F}(1-r)\sqrt{d}n^{-1})
+\omega_{F^{-1}}(\Lip{F}(1-r)\sqrt{d}n^{-1})+\Lip{F^{-1}}\epsilon\\
\leq&\big(\Lip{{\tilde{F}_{\nn}}^{-1}}+\Lip{F^{-1}}\big)\Lip{F}(1-r)\sqrt{d}n^{-1}+\Lip{F^{-1}}\epsilon.
\end{align*}
By the definition of $H^r$, we have
\begin{align*}
&\omega_{H^r}(\|{\tilde{F}_{\nn}}^{-1}-F^{-1}\|_{\mathcal{L}^\infty(F(K^r))})\leq \|{\tilde{F}_{\nn}}^{-1}-F^{-1}\|_{\mathcal{L}^\infty(F(K^r))}+2\sqrt{d}n^{-1}.
\end{align*}
Thus, we obtain
\begin{align*}
&\|H^r\circ{\tilde{F}_{\nn}}^{-1} -H^r\circ F^{-1}\|^2_{\mathcal{L}^2(F(K))}\\
\leq&|F(K^r)|\omega_{H^r}(\|{\tilde{F}_{\nn}}^{-1}-F^{-1}\|_{\mathcal{L}^\infty(F(K^r))})^2+|F(K\setminus K^r)|\|H^r\circ{\tilde{F}_{\nn}}^{-1} -H^r\circ F^{-1}\|_{\mathcal{L}^\infty(F(K\setminus K^r))}^2\\
\leq&\Lip{F}^{d}\big(\|{\tilde{F}_{\nn}}^{-1}-F^{-1}\|_{\mathcal{L}^\infty(F(K^r))}+2\sqrt{d}n^{-1}\big)^2
+(1-r^{d})\Lip{F}^{d}2^2\\
\leq&\Lip{F}^{d}\Big(\big(\Lip{{\tilde{F}_{\nn}}^{-1}}+\Lip{F^{-1}}\big)\Lip{F}(1-r)\sqrt{d}n^{-1}
+2\sqrt{d}n^{-1}+\Lip{F^{-1}}\epsilon\Big)^2+4(1-r^{d})\Lip{F}^{d}\\
\leq&\Lip{F}^{d}\Big[\Big(2+\big(\Lip{{\tilde{F}_{\nn}}^{-1}}+\Lip{F^{-1}}\big)\Lip{F}(1-r)\Big)\sqrt{d}n^{-1}
+\Lip{F^{-1}}\epsilon\Big]^2+4(1-r^{d})\Lip{F}^{d}.
\end{align*}
By Lemma \ref{lem:H^r}, there holds
\begin{align*}
\|H^r\circ F^{-1}- F^{-1}\|^2_{\mathcal{L}^2(F(K))}\leq \Lip{F}^d (2r-1)^2 d n^{-2}.
\end{align*}
Combining the preceding estimates completes the proof of the theorem.
\end{proof}

The above estimates are derived in the $\mathcal{L}^2$ space, and similar results can be
derived for the $\mathcal{L}^p$ spaces, $1\leq p<\infty$. The next result provides quantitative estimates for the INN approximations
$F_{\nn}$ and $F_{\nn}^{-1}$.
\begin{corollary}\label{rem:cNN}
By choosing $\epsilon=
c_\epsilon n^{-1}$ with some $c_\epsilon>0$ and
\begin{align*}
r\geq c^r_{\nn}:=\max\big(\big(1-(\Lip{\tilde{F}_{\nn}}+\Lip{F})^{-2}\big)^{\frac1d},1-\big(\Lip{{\tilde{F}_{\nn}}^{-1}} +\Lip{F^{-1}}\big)^{-1},(1-n^{-2})^{\frac1d} \big),
\end{align*}
the following error estimates hold
\begin{align*}
  \|F_{\nn}-F\|^2_{\mathcal{L}^2(K)}&\leq2\big[(3+\Lip{F}^2)d+3c_\epsilon^2\big]n^{-2},\\
\|F^{-1}_{\nn}-F^{-1}\|^2_{\mathcal{L}^2(F(K))}&\leq 2\Lip{F}^{d}\big[\big(2(2+\Lip{F})^2+1\big) d +2\Lip{F^{-1}}c_\epsilon+6\big]n^{-2}.
\end{align*}
\end{corollary}
\begin{proof}
By the choice of $\epsilon$ and $r$, we have
$\Lip{F}(1-r)\leq(\Lip{\tilde{F}_{\nn}}+\Lip{F})(1-r)\leq(\Lip{{\tilde{F}_{\nn}}}+\Lip{F}\big)(1-r^d)^{\frac12}\leq1$,
$\big(\Lip{{\tilde{F}_{\nn}}^{-1}}+\Lip{F^{-1}}\big)(1-r)\leq 1$ and $1-r^d\leq n^{-2}$.
Then it follows from Theorem \ref{thm:F_NN} that
\begin{align*}
\|F_{\nn}-F\|^2_{\mathcal{L}^2(K)}
\leq&2\big((\Lip{\tilde{F}_{\nn}}+\Lip{F})(1-r)\sqrt{d}+c_\epsilon\big)^2n^{-2}+2\Lip{F}^2(2r-1)^2dn^{-2}\\
&+2\big(\tfrac12(\Lip{\tilde{F}_{\nn}}+\Lip{F})\sqrt{d}+c_\epsilon\big)^2(1-r^d)n^{-2}\\
\leq&2\big[\big(\sqrt{d}+c_\epsilon\big)^2+\Lip{F}^2d+\big(\tfrac12\sqrt{d}+c_\epsilon n^{-1}\big)^2\big]n^{-2}\\
\leq&2\big[2d+2c_\epsilon^2+\Lip{F}^2d+\tfrac12 d+2c_\epsilon^2n^{-2}\big]n^{-2}
\leq2\big[(3+\Lip{F}^2)d+3c_\epsilon^2\big]n^{-2}.
\end{align*}
Similarly, we have
\begin{align*}
\|F^{-1}_{\nn}-F^{-1}\|^2_{\mathcal{L}^2(F(K))}\leq&
2\Lip{F}^{d}\big(\big(2+(\Lip{{\tilde{F}_{\nn}}^{-1}}+\Lip{F^{-1}})\Lip{F}(1-r)\big)\sqrt{d}+\Lip{F^{-1}}c_\epsilon\big)^2n^{-2}\\
&+2\Lip{F}^d (2r-1)^2 d n^{-2}+8\Lip{F}^{d}(1-r^{d})\\
\leq&2\Lip{F}^{d}\big[\big(\big(2+\Lip{F}\big)\sqrt{d}+\Lip{F^{-1}}c_\epsilon\big)^2+ d +4\big]n^{-2}\\
\leq&2\Lip{F}^{d}\big[\big(2(2+\Lip{F})^2+1\big) d +2\Lip{F^{-1}c_\epsilon}+6\big]n^{-2}.
\end{align*}
This completes the proof of the corollary.
\end{proof}

\begin{remark}
The approximation error decays at a rate $O(N^{-\frac{1}{d}})$, and thus it suffers from the
usual curse of dimensionality. This result agrees with the fact that for Lipschitz maps, standard fully
connected NNs also suffer from this curse similarly \cite[Section 8.7]{DeVoreHaninPetrova:2021}
{\rm(}up to a logarithmic factor{\rm)}. However, it is still unclear whether this is also the
lower bound for coupling based INNs. In order to overcome the curse, additional structural
information, e.g., manifold assumption or low complexity, is needed \cite[Section 8.10]{DeVoreHaninPetrova:2021}.
\end{remark}

\section{INN operator approximation}\label{sec:overall}

Now we consider the approximation of the following nonlinear operator equation using INNs:
\begin{equation}\label{eqn:ivp}
F^{\dag}(x)=y,
\end{equation}
where $F^{\dag}: \mathcal{D}(F^{\dag})\to\mathcal{I}(F^{\dag})$ is a bi-Lipschitz continuous mapping with its domain
$\mathcal{D}(F^{\dag})\subset \mathcal{X}$ and image space $\mathcal{I}(F^{\dag}) \subset\mathcal{Y}$, and
$\mathcal{X}$ and $\mathcal{Y}$ are infinite-dimensional separable Hilbert spaces with inner products $\langle
\cdot,\cdot\rangle$ and norms $\|\cdot\|$, respectively. The model \eqref{eqn:ivp} can describe many important scientific
problems. We aim at constructing a DNN that simultaneously approximates the forward process
$F^\dag:x\to y$ and the inverse process $(F^{\dag })^{-1}:y\to x$ by training the DNN on a finite collection of paired observations
$\{(x^n, y^n)\}_{n=1}^N$ (i.e., the training data). Throughout, we assume that the evaluation points $x^n$s are drawn independent
and identically distributed (i.i.d.) with respect to an (unknown) probability measure $\mu$ supported on $\mathcal{X}$,
and $y^n$s are i.i.d. with respect to the push-forward measure $F^{\dag}_{\sharp}\mu$.

\subsection{Principal component analysis}
We first recall model reduction by principal component analysis (PCA). Following the construction in
\cite{BhattacharyaHosseiniStuart:2021}, we project any $x\in\mathcal{X}$ and $y\in\mathcal{Y}$ into finite-dimensional
spaces by principal component analysis (PCA) on $\mathcal{X}$ and $\mathcal{Y}$, respectively. We approximate the identity
mappings $I_\mathcal{X}: \mathcal{X}\to\mathcal{X}$ and $I_\mathcal{Y}: \mathcal{Y}\to\mathcal{Y}$ by the composition of
two maps, known as the encoder and decoder, respectively, denoted by $G_{\mathcal{U}}: \mathcal{U}\to \mathbb{R}^{
d_{\mathcal{U}}}$ and $G^*_{\mathcal{U}}: \mathbb{R}^{d_{\mathcal{U}}}\to \mathcal{U}$ where $\mathcal{U}=\mathcal{X}$
or $\mathcal{Y}$, which have a finite-dimensional range and domain, such that
$G^*_{\mathcal{U}}\circ G_{\mathcal{U}}\approx I_\mathcal{U}$.
Specifically, let $\nu$ denote a probability measure supported on the space $\mathcal{U}$, and we make the assumption of
a finite fourth moment: $\E_{u\sim\nu}[\|u\|^4]<\infty$. We denote by $\{u^i\}_{i=1}^N$ a finite collection of $N$ samples
drawn independent and identically distributed (i.i.d.) from the probability measure $\nu$ that will be used as the training
data on which PCA is based. We consider the empirical, non-centered covariance operator $C_N^{\mathcal{U}}:=\frac1N
\sum_{i=1}^N u^i \otimes u^i$, where $\otimes$ denotes the outer product. $C_N^\mathcal{U}$ is a non-negative, self-adjoint,
trace-class operator on the Hilbert space $\mathcal{U}$, of rank at most $N$ with eigenvectors $\{\phi_\mathcal{U}^{N,i}
\}_{i=1}^N$ and its corresponding eigenvalues $\lambda^\mathcal{U}_{N,1}\geq\lambda^\mathcal{U}_{N,2}\geq\cdots\geq
\lambda^\mathcal{U}_{N,N}\geq 0$ in a nonincreasing order. We define that, with some truncation level $0<d_{\mathcal{U}}\leq N$,
\begin{align*}
\mbox{(PCA encoder)}   \quad\;\; G_{\mathcal{U}}(u):=&\big(\langle u, \phi_\mathcal{U}^{N,1} \rangle,\cdots,\langle u,
\phi_\mathcal{U}^{N,d_{\mathcal{U}}} \rangle\big)^t\in\mathbb{R}^{d_\mathcal{U}}, \quad\forall u\in\mathcal{U},\\
\mbox{(PCA decoder)}   \quad G^*_{\mathcal{U}}(u^*):=&\sum_{i=1}^{d_\mathcal{U}} u^*_{i}\phi_\mathcal{U}^{N,i},
\quad\forall u^*=(u^*_1,\cdots,u^*_{d_\mathcal{U}})^t\in \mathbb{R}^{d_\mathcal{U}}.
\end{align*}
Below we denote the associated projection operator via PCA by $\mathcal{R}^\mathcal{U}_N$, which clearly depends on
the random samples $\{u^i\}_{i=1}^N$ used to compute the empirical covariance $C_N^\mathcal{U}$.

With the PCA on the spaces $\mathcal{X}$ and $\mathcal{Y}$, we can now describe the reduced model. Specifically, we
transform the evaluations $\{(x^{\dag,i}, y^{\dag,i})\}_{i=1}^N$ to $\{x^i:=G_{\mathcal{X}}(x^{\dag,i}), y_i:=
G_{\mathcal{Y}}(y^{\dag,i})\}_{i=1}^N$, and finally return them back to $\{\tilde{x}^i:=G^*_{\mathcal{X}}(x^i),
\tilde{y}^i:=G^*_{\mathcal{Y}}(y^i)\}_{i=1}^N$. Throughout, we take $\{(\tilde{x}^i, \tilde{y}^i)\}_{i=1}^N\subset
\{(x^{\dag,i}, y^{\dag,i})\}_{i=1}^N$ which implies $F^\dag(\tilde{x})=\tilde{y}$, and define the following mapping
\begin{align}\label{eqn:reducedivp}
  F:=G_{\mathcal{Y}}\circ F^\dag\circ G^{*}_{\mathcal{X}}: G_{\mathcal{X}}\big(\mathcal{D}(F^{\dag})\big)\subset
  \mathbb{R}^{d_\mathcal{X}}\to G_{\mathcal{Y}}\big(\mathcal{I}(F^{\dag})\big)\subset\mathbb{R}^{d_\mathcal{Y}}.
\end{align}
Then for any $x^i\in G_{\mathcal{X}}(\mathcal{D}(F^{\dag}))$, $i=1,\ldots,N$, there holds
$\tilde{x}^i=G^*_{\mathcal{X}}(x^i)$ and $y^i=G_{\mathcal{Y}}(y^{\dag,i})=G_{\mathcal{Y}}(\tilde{y}^i)$.
Since $F^\dag(\tilde{x}^i)=\tilde{y}^i$ , we have
$F(x^i)=G_{\mathcal{Y}}\circ F^\dag\circ G^{*}_{\mathcal{X}}(x^i)
=G_{\mathcal{Y}}\circ F^\dag(\tilde{x}^i)=G_{\mathcal{Y}}(\tilde{y}^i)=y^i$.
Further, we have the inverse mapping $F^{-1}$ of $F$ given by
\begin{align*}
F^{-1}=&G_{\mathcal{X}}\circ (F^\dag)^{-1}\circ G^{*}_{\mathcal{Y}}:G_{\mathcal{Y}}\big(\mathcal{I}(F^{\dag})\big)\to G_{\mathcal{X}}\big(\mathcal{D}(F^{\dag})\big),
\end{align*}
since for any $x\in G_{\mathcal{X}}\big(\mathcal{D}(F^{\dag})\big)$, there holds
$F^{-1}\circ F(x)=G_{\mathcal{X}}\circ (F^\dag)^{-1}\circ G^{*}_{\mathcal{Y}}\circ G_{\mathcal{Y}}\circ F^\dag\circ G^{*}_{\mathcal{X}}(x)=x$.

\subsection{Proposed approach for approximating bi-Lipschitz maps}

Now we construct an approximation of the bi-Lipschitz map $F^\dag:\mathcal{X}\to \mathcal{Y}$ that approximates
both $F^\dag$ and its inverse $(F^\dag)^{-1}$ simultaneously. It proceeds in two steps:
\begin{itemize}
  \item[Step 1]  approximate the true
mapping $F^\dag$ by the reduced one $G^{*}_{\mathcal{Y}}\circ F\circ G_{\mathcal{X}}$ with the PCA encoder $G_{\mathcal{X}}$
and decoder $G^{*}_{\mathcal{Y}}$, with $F$ given in \eqref{eqn:reducedivp}, and also approximate the inverse mapping $(F^\dag)^{-1}$ by
$G^{*}_{\mathcal{X}}\circ F^{-1}\circ G_{\mathcal{Y}}$ with the PCA encoder $G_{\mathcal{Y}}$ and
the decoder $G^{*}_{\mathcal{X}}$.
  \item[Step 2] approximate the forward map $G^{*}_{\mathcal{Y}}
\circ F\circ G_{\mathcal{X}}$ and its inverse $G^{*}_{\mathcal{X}}\circ F^{-1}\circ G_{\mathcal{Y}}$ (with $F$ and $F^{-1}$ defined on $\mathbb{R}^d$) by $G^{*}_{\mathcal{Y}} \circ F_{\nn}\circ G_{\mathcal{X}}$ and $G^{*}_{\mathcal{X}}\circ F^{-1}_{\nn}\circ G_{\mathcal{Y}}$ in
Section \ref{sec:INN}, respectively, where $F_{\nn}$ is a coupling based INN.
\end{itemize}

\begin{figure}
\centering
\begin{tikzpicture}[font=\sf]
\node[draw,circle,label=below:$\mathbb{R}^d$] (n1) at ({0},-1){};
\node[draw,circle,label=above:$\mathbb{R}^d$] (n2) at (0,1){};
\node[draw,circle,label=above:$\mathcal{X}$] (s2) at (-4,1){};
\node[draw,circle,label=above:$\mathcal{X}$] (s3) at (4,1){};
\node[draw,circle,label=below:$\mathcal{Y}$] (s1) at (-4,-1){};
\node[draw,circle,label=below:$\mathcal{Y}$] (s4) at (4,-1){};

\draw[blue,thick,-latex] (-0.1,-.8) -- (-0.1,.8);
\draw[red,thick,-latex] (0.1,.8) -- (0.1,-.8);
\draw[red,thick,-latex] (s2) -- (n2);
\draw[blue,thick,-latex] (n2) -- (s3);
\draw[blue,thick,-latex] (s4) -- (n1);
\draw[red,thick,-latex] (n1) -- (s1);
\draw[thick,-latex] (s2) -- (s1);
\node at (-4.3,0) {$F^{\dag}$};
\node at (4.6,0) {$(F^\dag)^{-1}$};
\draw[thick,-latex] (s4) -- (s3);
\node at (-2,1.2) {\color{red}Encoder $G_\mathcal{X}$};
\node at (2,1.2) {\color{blue}Decoder $G^*_\mathcal{X}$};
\node at (2,-1.3) {\color{blue}Encoder $G_\mathcal{Y}$};
\node at (-2,-1.3) {\color{red}Decoder $G^*_\mathcal{Y}$};
\node at (-.5,0) {\color{blue}$F_{\rm nn}^{-1}$};
\node at (0.5,0) {\color{red}$F_{\rm nn}$};

\draw[dashed,thick,-latex] (-3.8,1.1) .. controls ($(n2)+(0.0cm,1.5cm)$) .. (3.8,1.1);
\node at (-0,2.5) {$I_{\mathcal{X}}\approx G_\mathcal{X}^*G_\mathcal{X}$};
\draw[dashed,thick,-latex] (3.8,-1.1) .. controls ($(n1)-(0.0cm,1.5cm)$) .. (-3.8,-1.1);
\node at (-0,-2.5) {$I_{\mathcal{Y}}\approx G_\mathcal{Y}^*G_\mathcal{Y}$};
\end{tikzpicture}
\caption{A schematic diagram for the proposed simultaneous approximations to
$F^\dag$ {\color{blue}and $(F^\dag)^{-1}$ based on the INN $F_{\nn}$} and its inverse $F_{\nn}^{-1}$.\label{fig:INN-diag}}
\end{figure}
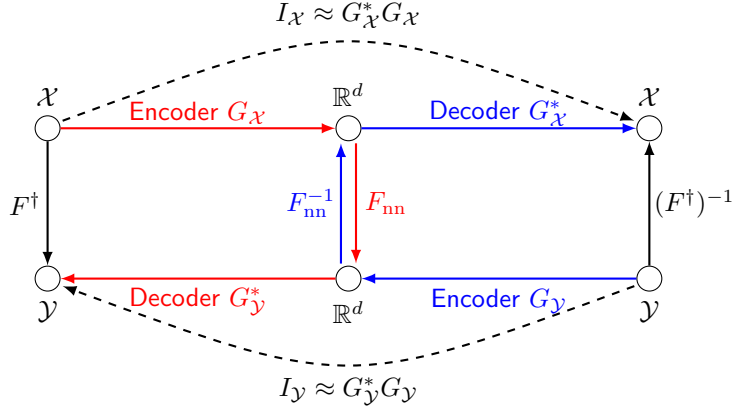

The diagram in Fig. \ref{fig:INN-diag} illustrates the complete procedure of constructing and inferring the INN approximation.
Given the training dataset $\{x^i,y^i\}_{i=1}^N$, we first generate the pairs of encoder and decoder on the
space $\mathcal{X}$ and $\mathcal{Y}$,  i.e., $(G_\mathcal{X},G_\mathcal{X}^*)$ and $(G_\mathcal{Y},G_\mathcal{Y}^*)$,
and then learn the coupling-based NN $F_{\nn}$ by minimizing the following
empirical loss
\begin{equation*}
  \mathcal{L}(\theta) = \sum_{i=1}^N\left(\frac{c_0}{2}\| G^{*}_{\mathcal{X}}\circ F_{\rm nn}^{-1}\circ G_{\mathcal{Y}}(y^i)-x^i\|_\mathcal{X}^2  + \frac12 \|G^{*}_{\mathcal{Y}}\circ F_{\nn}\circ G_{\mathcal{X}}(x^i)-y^i\|^2_\mathcal{Y}\right),
\end{equation*}
where $c_0>0$ is a scalar that controls the balance between the two terms in the loss, and $\theta$ denotes the vector of trainable
parameters in the INN $F_{\rm nn}$ approximating the reduced map. The optimization problem is typically minimized by
gradient type methods, e.g., Adam \cite{KingmaBa:2015}.

The obtained approximation involves several different sources of errors: the truncation / approximation
error due to the use of PCA encoder / decoder, the approximation error due to the INN approximation
to the bi-Lipschitz map, generalization error due to using a finite-number of samples (instead of the underlying
distribution), and the optimization error due to the optimizer (which fails to find a global minimizer
of the empirical loss $\mathcal{L}(\theta)$, since the loss landscape is highly complex). Our analysis below
focuses on the first two components. The key tool is to analyze the approximation error of the INN
approximation in Section \ref{sec:INN}.

\subsection{Error analysis of INN operator approximation}
Now we give the error analysis of the approximation for both forward and inverse processes,
which involve both PCA model reduction and INN approximation of bi-Lipschitz map on finite-dimensional spaces.
We prove that there exists an INN that achieves small
approximation error in Theorem \ref{thm:full}.

\subsubsection{Error estimate of PCA model reduction}

Now we discuss the error due to the PCA. There are two sources of errors: one is due to the truncation, and
the other due to random sampling (in approximating the covariance operator by an empirical one). We denote
the covariance operator of the infinite-dimensional input data space $\mathcal{X}$ and output data space
$\mathcal{Y}$ by $C^\mathcal{U}:=\E_{u\sim \nu}[u\otimes u]$, with $\mathcal{U}=\mathcal{X}$ or $\mathcal{Y}$,
respectively. By the i.i.d. assumption on the samples $\{u^i\}_{i=1}^N$, we have $C^\mathcal{U}=
\E_{\{u^i\}\sim \nu}[C^\mathcal{U}_N]$, where $\E[\cdot]$ denotes taking expectation with respect to the
samples from the probability measure $\nu$. We also define the eigenvectors of $C^\mathcal{U}$ by $\{
\phi_\mathcal{U}^{i}\}_{i=1}^\infty$ and its corresponding nonnegative eigenvalues $\lambda^\mathcal{U}_{1}
\geq\lambda^\mathcal{U}_{2}\geq\cdots\geq \cdots$ in a nonincreasing order. Then the PCA encoder and decoder
for the infinite-dimensional space $\mathcal{U}$ are defined respectively by
\begin{align*}
\mbox{(PCA encoder)}   \quad\;\; G_{\mathcal{U}}(u):=&\big(\langle u, \phi_\mathcal{U}^{1} \rangle,
\cdots,\langle u, \phi_\mathcal{U}^{d_{\mathcal{U}}} \rangle\big)^t\in \mathbb{R}^{d_\mathcal{U}} , \quad\forall u\in\mathcal{U},\\
\mbox{(PCA decoder)}   \quad G^{*}_{\mathcal{U}}(u^*):=&\sum_{i=1}^{d_\mathcal{U}} u^*_{i}\phi^\mathcal{U}_{i},
\quad\forall u^*=(u^*_1,\cdots,u^*_{d_\mathcal{U}})^t\in \mathbb{R}^{d_\mathcal{U}}.
\end{align*}
We define $d:=d_\mathcal{X}=d_\mathcal{Y}\geq 2$ by choosing suitable truncation indices $d_\mathcal{X}$ and
$d_\mathcal{Y}$. 
Let $\tilde u = \mathcal{R}_N^\mathcal{U}u^\dag=G^{*}_{\mathcal{U}}\circ G_{\mathcal{U}}u^\dag$ ($\mathcal{R}_N^\mathcal{U}$ is the
orthogonal projection into the PCA subspace).

Now, we bound the mean-squared errors $\E_{\{x^{\dag,i}\}\sim \nu}\big[\E_{x^\dag\sim \nu}[ \|x^{\dag}-\tilde{x}\|^2]\big]$
and $\E_{\{y^{\dag,i}\}\sim {F^{\dag}_{\sharp}\nu}}\big[\E_{y^\dag\sim F^\dag_\sharp\nu}[ \|y^{\dag}-\tilde{y}\|^2]\big]$
caused by computing the PCA encoder $G_{\mathcal{U}}$ and decoder $G^{*}_{\mathcal{U}}$ from the finite data
$\{(x^{\dag,i}, y^{\dag,i})\}_{i=1}^N$ and truncating the order of the reduced model, on the spaces
$\mathcal{X}$ and $\mathcal{Y}$, separately.
\begin{lemma}\label{lem:PCAerr}
Given a finite collection of evaluations $\{u^{\dag,i}\}_{i=1}^N$, which are i.i.d. with respect to a probability
measure $\nu$ supported on $\mathcal{U}$, then with $c_\nu=\Big(\E_{u\sim\nu}[{\rm tr}\big((u\otimes u-C^{\mathcal{U}}
)^t (u\otimes u-C^{\mathcal{U}})\big)]\Big)^\frac12$, there holds
\begin{align*}
\E_{\{u^{\dag,i}\}\sim \nu}\big[\E_{u^\dag\sim \nu}[ \|u^{\dag}-\tilde{u}\|^2]\big]\leq c_\nu\sqrt{d}N^{-\frac12}+\sum_{j=d+1}^{\infty}\lambda^\mathcal{U}_j.
\end{align*}
\end{lemma}
\begin{proof}
The whole proof is shown in \cite[Theorem 3.4]{BhattacharyaHosseiniStuart:2021} that
\begin{align*}
\E_{\{u^{\dag,i}\}\sim \nu}\big[\E_{u^\dag\sim \nu}[ \|u^{\dag}-\tilde{u}\|^2]\big]&\leq c_\nu\sqrt{d}N^{-\frac12}
+\E_{u^\dag\sim \nu}[ \|u^{\dag}-G^{*}_{\mathcal{U}} \circ G_{\mathcal{U}}(u^\dag)\|^2]=c_\nu\sqrt{d}N^{-\frac12}+\sum_{j=d+1}^{\infty}\lambda^\mathcal{U}_j.
\end{align*}
This completes the proof of the lemma.
\end{proof}

Now we bound the error between the approximate mapping $G^{*}_{\mathcal{Y}}
\circ F\circ G_{\mathcal{X}}$ and the exact one $F^\dag$, and the analogue for the
inverse $(\mathcal{F}^{\dag})^{-1}$. In the statement, $L_F$
denotes the Lipschitz constant of the map $F$, and the constants $c_\mu$ and
$c_{F^\dag \mu}$ are defined in Lemma \ref{lem:PCAerr}.
\begin{theorem}\label{thm:PCAerr}
For the operator $F$ defined in \eqref{eqn:reducedivp}, the following error estimates hold
\begin{align*}
&\E_{\{x^{\dag,i}\}\sim \mu}\big[\E_{x^\dag\sim \mu}[ \|G^{*}_{\mathcal{Y}}\circ F\circ G_{\mathcal{X}}(x^\dag)-F^\dag(x^\dag)\|^2]\big]\\
\leq&2(\Lip{F^\dag}^2c_\mu +c_{F^{\dag}_{\sharp}\mu})\sqrt{d}N^{-\frac12}+2\Big(\Lip{F^\dag}^2\sum_{j=d+1}^{\infty}\lambda^\mathcal{X}_j+\sum_{j=d+1}^{\infty}\lambda^\mathcal{Y}_j\Big),\\
&\E_{\{y^\dag_i\}\sim F^{\dag}_{\sharp}\mu}\big[\E_{y^\dag\sim F^{\dag}_{\sharp}\mu}[ \|G^{*}_{\mathcal{X}}\circ F^{-1}\circ G_{\mathcal{Y}}(y^\dag)-(F^\dag)^{-1}(y^\dag)\|^2]\big]\\
\leq&2(\Lip{F^{\dag-1}}^2c_{F^{\dag}_{\sharp}\mu} +c_{\mu})\sqrt{d}N^{-\frac12}+2\Big(\Lip{F^{\dag-1}}^2\sum_{j=d+1}^{\infty}\lambda^\mathcal{Y}_j+\sum_{j=d+1}^{\infty}\lambda^\mathcal{X}_j\Big).
\end{align*}
\end{theorem}
\begin{proof}
By the definition of the operator $F$, the triangle inequality and Cauchy-Schwarz inequality, with $\tilde x=G^{*}_{\mathcal{X}}\circ G_{\mathcal{X}}(x^\dag)$ and $\tilde y =G^{*}_{\mathcal{Y}}\circ  G_{\mathcal{Y}}\circ y^\dag$, we have
\begin{align*}
 &\|G^{*}_{\mathcal{Y}}\circ F\circ G_{\mathcal{X}}(x^\dag)-F^\dag(x^\dag)\|^2
=\|G^{*}_{\mathcal{Y}}\circ G_{\mathcal{Y}}\circ F^\dag\circ G^{*}_{\mathcal{X}}\circ G_{\mathcal{X}}(x^\dag)-F^\dag(x^\dag)\|^2\\
\leq&2\|G^{*}_{\mathcal{Y}}\circ G_{\mathcal{Y}}\circ F^\dag(\tilde{x})-G^{*}_{\mathcal{Y}}\circ G_{\mathcal{Y}}\circ F^\dag(x^\dag) \|^2+2\|G^{*}_{\mathcal{Y}}\circ  G_{\mathcal{Y}}\circ y^\dag-y^\dag\|^2\\
\leq&2\Lip{G^{*}_{\mathcal{Y}}\circ G_{\mathcal{Y}}\circ F^\dag}^2\| \tilde{x}-x^\dag \|^2+2\|\tilde{y}-y^\dag\|^2
\leq2\Lip{F^\dag}^2\| \tilde{x}-x^\dag \|^2+2\|\tilde{y}-y^\dag\|^2,
\end{align*}
since $G^{*}_{\mathcal{Y}}\circ G_{\mathcal{Y}}$ is an orthogonal projection.
Using the error estimate in Lemma \ref{lem:PCAerr}, we derive
\begin{align*}
&\E_{\{x^{\dag,i}\}\sim \mu}\big[\E_{x^\dag\sim \mu}[ \|G^{*}_{\mathcal{Y}}\circ F\circ G_{\mathcal{X}}(x^\dag)-F^\dag(x^\dag)\|^2]\big]\\
\leq& 2\Lip{F^\dag}^2\Big(c_\mu \sqrt{d}N^{-\frac12}+\sum_{j=d+1}^{\infty}\lambda^\mathcal{X}_j\Big)+2\Big(c_{F^{\dag}_{\sharp}\mu}\sqrt{d}N^{-\frac12} +\sum_{j=d+1}^{\infty}\lambda^\mathcal{Y}_j\Big).
\end{align*}
Similarly, there holds
\begin{align*}
&\|G^{*}_{\mathcal{X}}\circ F^{-1}\circ G_{\mathcal{Y}}(y^\dag)-F^{\dag-1}(y^\dag)\|^2
=\|G^{*}_{\mathcal{X}}\circ G_{\mathcal{X}}\circ F^{\dag-1}\circ G^{*}_{\mathcal{Y}}\circ G_{\mathcal{Y}}(y^\dag)-(F^\dag)^{-1}(y^\dag)\|^2\\
\leq&2\Lip{G^{*}_{\mathcal{X}}\circ G_{\mathcal{X}}\circ F^{\dag-1}}^2\|\tilde{y}-y^\dag\|^2+2\| \tilde{x}-x^\dag\|^2
\leq2\Lip{F^{\dag-1}}^2 \|\tilde{y}-y^\dag\|^2+2\| \tilde{x}-x^\dag\|^2,
\end{align*}
and hence
\begin{align*}
&\E_{\{y^{\dag,i}\}\sim F^{\dag}_{\sharp}\mu}\big[\E_{y^\dag\sim F^{\dag}_{\sharp}\mu}[ \|G^{*}_{\mathcal{X}}\circ F^{-1}\circ G_{\mathcal{Y}}(y^\dag)-F^{\dag-1}(y^\dag)\|^2]\big]\\
\leq&2\Lip{F^{\dag-1}}^2\Big(c_{F^{\dag}_{\sharp}\mu}\sqrt{d}N^{-\frac12}+\sum_{j=d+1}^{\infty}\lambda^\mathcal{Y}_j\Big)+2\Big(c_\mu \sqrt{d}N^{-\frac12}+\sum_{j=d+1}^{\infty}\lambda^\mathcal{X}_j\Big),
\end{align*}
which directly completes the proof of the theorem.
\end{proof}

\begin{remark}
In the discussion, we have ignored the possible presence of data noise in $y^{\dag,i}$. Should there be any
noise $\{\xi_i\}_{i=1}^N$ in the observation, it will impact the accuracy of the reduced model. The mean
squared errors indicate that the bound depends on the spectral decay of the covariance operators
$C^\mathcal{X}$ and $C^\mathcal{Y}$ via the remainders $\sum_{j=d+1}^{\infty}\lambda^\mathcal{X}_j$ and $\sum_{j=d+1}^{\infty}
\lambda^\mathcal{Y}_j$ and the number of samples $N$ {\rm(}at a rate $N^{-\frac12}${\rm)}. In practice, the
covariance operator $C^\mathcal{X}$ is often represented by an integral operator, and then its spectral
decay can be characterized by the smoothness of the associated kernel \cite{GriebelLi:2018}.
\end{remark}

\subsubsection{Error estimate of INN operator approximation}
Now we derive the full error estimate of the INN approximation of the forward and inverse mappings. The
following theorem is the main result of this section.

\begin{theorem}\label{thm:full}
Let $F_{\nn}:=\tilde{F}_{\nn}\circ H^r_{\nn}$ and $F^{-1}_{\nn}=H^r_{\nn}\circ\tilde{F}^{-1}_{\nn}$ when $\epsilon=
c_\epsilon n^{-1}$ with some $c_\epsilon>0$ and $r\geq c^r_{\nn}$, where $H^r_{\nn}(x):=(h^r_{\nn}(x_1),\cdots,h^r_{\nn}(x_{d}))$
with a two-layer NN $h^r_{\nn}$ of $2(N^{\frac1d}+1)$ neurons, and $\tilde{F}_{\nn}$ has $1+16N$ weight
layers, $1+6N$ added layers and $N^{\frac1d}$ neurons.
Then the expected errors of the whole systems are bounded by
\begin{align*}
&\E_{\{x^\dag_i\}\sim \mu}\big[\E_{x^\dag\sim \mu}[ \|G^{*}_{\mathcal{Y}}\circ F_{\nn}\circ G_{\mathcal{X}}(x^\dag)-F^\dag(x^\dag)\|^2]\big]\\
\leq&2c_{\nn}N^{-\frac{2}{d}}+4(\Lip{F^\dag}^2c_\mu +c_{F^{\dag}_{\sharp}\mu})\sqrt{d}N^{-\frac12}+4(\Lip{F^\dag}^2\sum_{j=d+1}^{\infty}\lambda^\mathcal{X}_j+\sum_{j=d+1}^{\infty}\lambda^\mathcal{Y}_j),\\
&\E_{\{y^\dag_i\}\sim F^{\dag}_{\sharp}\mu}\big[\E_{y^\dag\sim F^{\dag}_{\sharp}\mu}[ \|G^{*}_{\mathcal{X}}\circ F^{-1}_{\nn}\circ G_{\mathcal{Y}}(y^\dag)-(F^\dag)^{-1}(y^\dag)\|^2]\big]\\
\leq& 2c_{\nn}'N^{-\frac{2}{d}}+4(\Lip{F^{\dag-1}}^2c_{F^{\dag}_{\sharp}\mu} +c_{\mu})\sqrt{d}N^{-\frac12}+4\Big(\Lip{F^{\dag-1}}^2\sum_{j=d+1}^{\infty}\lambda^\mathcal{Y}_j+\sum_{j=d+1}^{\infty}\lambda^\mathcal{X}_j\Big),
\end{align*}
with $c_{\nn}=2\big[(3+\Lip{F}^2)d+3c_\epsilon^2\big]$ and $c_{\nn}'=2\Lip{F}^{d}[\big(2(2+\Lip{F})^2+1\big) d +2\Lip{F^{-1}}c_\epsilon+6]$.
\end{theorem}
\begin{proof}
By the triangle inequality, we have
\begin{align*}
 \|G^{*}_{\mathcal{Y}}\circ F_{\nn}\circ G_{\mathcal{X}}(x^\dag)-F^\dag(x^\dag)\|^2
\leq&2 \|G^{*}_{\mathcal{Y}}\circ F_{\nn}\circ G_{\mathcal{X}}(x^\dag)-G^{*}_{\mathcal{Y}}\circ F\circ G_{\mathcal{X}}(x^\dag)\|^2\\
 &+2\|G^{*}_{\mathcal{Y}}\circ F\circ G_{\mathcal{X}}(x^\dag)-F^\dag(x^\dag)\|^2:={\rm I}+{\rm II}.
\end{align*}
It suffices to bound the two terms separately. For the term ${\rm I}$, we have
\begin{align*}
  \E_{\{x^\dag_i\}\sim \mu}\big[\E_{x^\dag\sim \mu}[{\rm I}]\big]
  \leq&2\E_{\{x_i\}\sim {G_{\mathcal{X}}}_{\sharp}\mu}\big[\E_{x\sim {G_{\mathcal{X}}}_{\sharp}\mu}[ \|G^{*}_{\mathcal{Y}}\circ F_{\nn} (x)-G^{*}_{\mathcal{Y}}\circ F(x)\|^2]\big]\\
\leq&2\E_{\{x_i\}\sim {G_{\mathcal{X}}}_{\sharp}\mu}\big[\|G^{*}_{\mathcal{Y}}\|^2 \|F_{\nn}- F\|_{\mathcal{L}^2(K)}^2\big]\leq2\|F_{\nn}- F\|_{\mathcal{L}^2(K)}^2.
\end{align*}
In view of Remark \ref{rem:cNN}, there holds
\begin{align*}
\E_{\{x^\dag_i\}\sim \mu}\big[\E_{x^\dag\sim \mu}[{\rm I}]\big]
\leq2c_{\nn }n^{-2}\leq2c_{\nn}N^{-\frac{2}{d}}.
\end{align*}
Meanwhile, by Theorem \ref{thm:PCAerr}, we have
\begin{align*}
  \E_{\{x^\dag_i\}\sim \mu}\big[\E_{x^\dag\sim \mu}[{\rm II}]\big] \leq 4(\Lip{F^\dag}^2c_\mu +c_{F^{\dag}_{\sharp}\mu})\sqrt{d}N^{-\frac12}+4\Big(\Lip{F^\dag}^2\sum_{j=d+1}^{\infty}\lambda^\mathcal{X}_j +\sum_{j=d+1}^{\infty}\lambda^\mathcal{Y}_j\Big).
\end{align*}
Combining the last two estimates yields the first assertion. Similarly, there holds
\begin{align*}
 \|G^{*}_{\mathcal{X}}\circ F^{-1}_{\nn}\circ G_{\mathcal{Y}}(y^\dag)-(F^\dag)^{-1}(y^\dag)\|^2
\leq&2 \|G^{*}_{\mathcal{X}}\circ F^{-1}_{\nn}\circ G_{\mathcal{Y}}(y^\dag)-F^{-1}(y^\dag)\|^2\\
  &+2\|G^{*}_{\mathcal{X}}\circ F^{-1}\circ G_{\mathcal{Y}}(y^\dag)-(F^\dag)^{-1}(y^\dag)\|^2 =:{\rm III} + {\rm IV}.
\end{align*}
By Remark \ref{rem:cNN} and Theorem \ref{thm:PCAerr}, we derive
\begin{align*}
  \E_{\{y^\dag_i\}\sim F^{\dag}_{\sharp}\mu}\big[\E_{y^\dag\sim F^{\dag}_{\sharp}\mu}[{\rm III}]\big]
  & \leq 2\|F^{-1}_{\nn}-F^{-1}\|^2_{\mathcal{L}^2(F(K))} \leq 2c_{\nn}'N^{-\frac{2}{d}},\\
  \E_{\{y^\dag_i\}\sim F^{\dag}_{\sharp}\mu}\big[\E_{y^\dag\sim F^{\dag}_{\sharp}\mu}[ {\rm IV}]\big]
  & \leq 4\big(\Lip{F^{\dag-1}}^2c_{F^{\dag}_{\sharp}\mu} +c_{\mu}\big)\sqrt{d}N^{-\frac12}
   +4\Big(\Lip{F^{\dag-1}}^2\sum_{j=d+1}^{\infty}\lambda^\mathcal{Y}_j +\sum_{j=d+1}^{\infty}\lambda^\mathcal{X}_j\Big).
\end{align*}
Combining the preceding estimates completes the proof of the theorem.
\end{proof}

\begin{remark}
Theorem \ref{thm:full} gives only the existence of an INN that approximates the forward
and inverse maps simultaneously to a certain tolerance. It does not address the important
issue of realizing the approximation in practice via optimizing the empirical loss $\mathcal{L}
(\theta)$, for which one has to analyze also the statistical error and optimization error. The
former would indicate how many samples are needed in order to achieve the tolerance, and
the latter is notoriously challenging due to the complex landscape of the associated
optimization problem. We leave the study of these important errors to future works.
\end{remark}
\section{Numerical experiments} \label{sec:numer}
In this section, we present preliminary numerical results to showcase the feasibility
of the proposed INN based approach to simultaneously approximate the
forward and inverse processes of a bi-Lipschitz map on infinite-dimensional spaces.
Consider the following second-order elliptic PDE on a smooth bounded domain $D$
(with a boundary $\partial D$):
\begin{align}\label{eqn:pde}
\left\{
\begin{aligned}
     -\nabla\cdot (u(x)\nabla y(x))&=f, \quad  \mbox{in }D, \\
     y&=0, \quad \mbox{on } \partial D,
\end{aligned}
\right.
\end{align}
where $u\in L^\infty(\Omega)$ is the diffusion coefficient in the PDE, and $f\in L^2(D)$ is the given
source term. The diffusion coefficient $u$ is assumed to satisfy the standard ellipticity and boundedness
assumption so that for any given $u$, there exists a unique solution $y\equiv y(u)\in H_0^1(\Omega)$.
The map $F$ of interest is defined by $F:u\mapsto y$. Note that this setting has been extensively studied in
the literature \cite{CohenDeVore:2015}, and it often serves as a model problem for parameter
identification of partial differential equations.

In the numerical experiments below, we take a unit square $D=(0,1)^2$, and fix $f\equiv1$. To represent
the diffusion coefficient $u$, we employ truncated Karhunen-Lo\`{e}ve expansion
\cite{GhanemSpanos:1991,BabuskaTempone:2004}. We simulate $M= 10000$ samples of the
coefficient $u$, generated by the following Fourier representation
\begin{equation*}
u(x,\xi)=2+\sum_{i,j=1}^{20}\frac{\xi_{i,j}}{i^3+j^3}\cos(i\pi x_1)\cos(j\pi x_2),\quad x=[x_1,x_2]\in D,
\end{equation*}
where $\xi=[\xi_{1,1},\xi_{1,2},\cdots,\xi_{1,20},\xi_{2,1},\cdots,\xi_{20,20}]\in \mathbb{R}^{20\times 20}$ and each $\xi_{i,j}$
follows the standard normal distribution $N(0,1)$ (for each entry independently). For each realization $u$, we
discretize the boundary value problem \eqref{eqn:pde} with a mesh size $h=\frac{1}{50}$
using the standard piecewise linear Galerkin finite element method, implemented in FEniCS \cite{LoggWells:2012}, to
obtain the corresponding output $y_h\in \mathbb{R}^{51\times 51}$.
For the ease of exposition, we apply a one-to-one transformation $T$ on $u$, given by $T(u)=\xi\in \mathbb{R}^{20\times 20}$,
to transform the infinite-dimensional inputs to a finite-dimensional Euclidean space without additional discretization.
Then we reduce the dimensionalities of both inputs $\xi$ and outputs $y_h$ by PCA (using the implementation in \texttt{scikit-learn},
available from \url{https://github.com/scikit-learn/scikit-learn/releases/tag/1.2.2}) with a reduced dimensionality $d=10$ (which capture more than 99\% of the spectral energy) and obtain the truncated inputs
$\hat{u}=[\hat{u}_1,\cdots,\hat{u}_{10}]\in\mathbb{R}^{10}$ and truncated outputs $\hat{y}=[\hat{y}_1,\cdots,\hat{y}_{10}]
\in\mathbb{R}^{10}$. We also record the $10$ largest singular values of the empirical covariance to form normalized, in terms of the $\ell^1$-norm, weight vectors $w_u\in\mathbb{R}^{10}$ and $w_y\in\mathbb{R}^{10}$. This normalization amounts to working in suitable Sobolev
space scales so that the map is nearly bi-Lipschitz, in accordance with the known stability estimates \cite{JinZhou:2021}.

\begin{figure}[hbt!]
\centering
\begin{tabular}{cc}
\includegraphics[width=0.8\textwidth,trim={0cm 0cm 0cm 0cm}, clip]{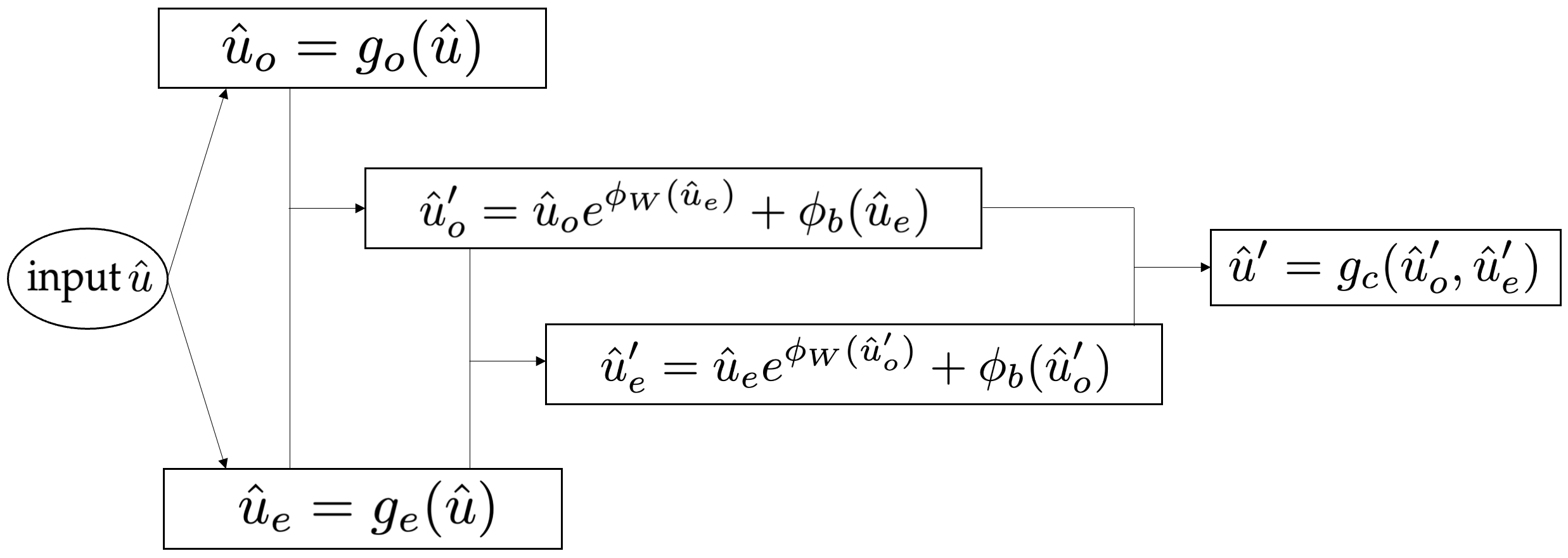} \\
(a) affine coupling block for the forward map\\
\includegraphics[width=0.8\textwidth,trim={0cm 0cm 0cm 0cm}, clip]{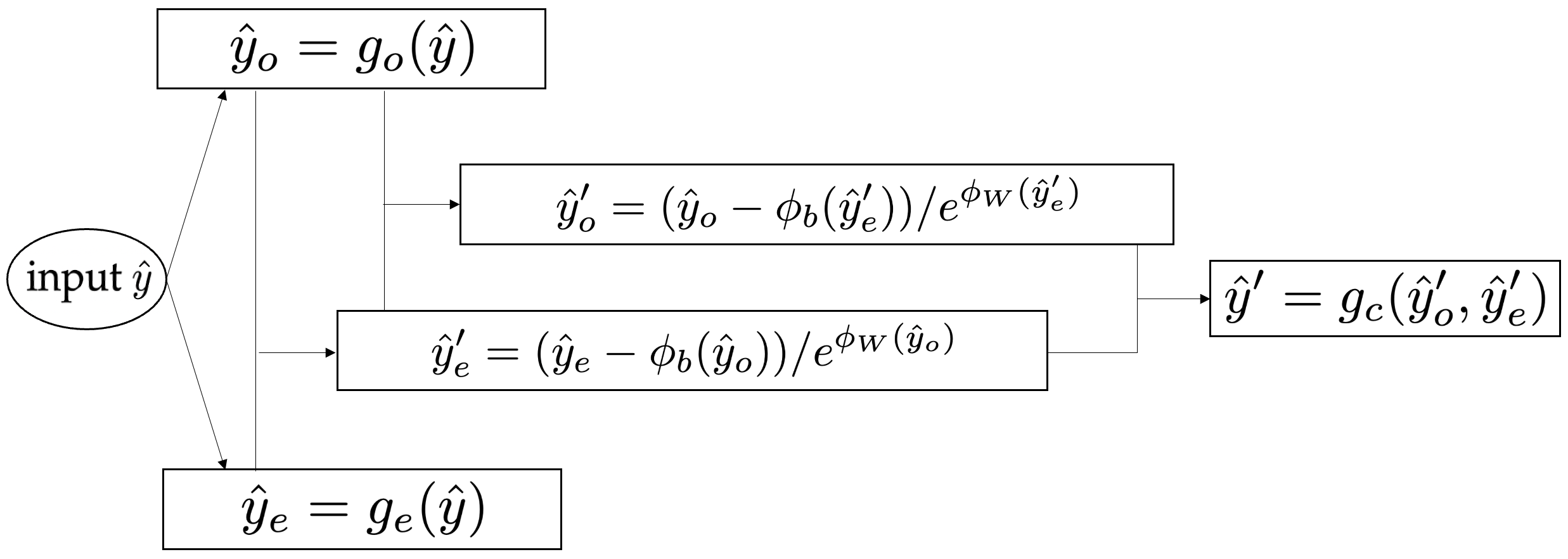}\\
(b) affine coupling block for the inverse map
\end{tabular}
\caption{The forward affine coupling block with input $\hat{u}$ and output $\hat{u}'$ (top) and the inverse affine
coupling block with input $\hat{y}$ and output $\hat{y}'$ (bottom). Both $\phi_W$ and $\phi_b$ are four-layer full
connected NNs with ReLU.\label{fig:couplingNN_forward}}
\end{figure}

In our experiments, we use an identical NN architecture with 3 blocks. One block of the forward process
and inverse process is sketched in Fig. \ref{fig:couplingNN_forward}, which are commonly used in an affine coupling
layer \cite{DinhSohlBengio:2017,ArdizzoneKruse:2019}. In the figure, we define the maps $g_{o}$ (odd part), $g_{e}$
(even part) and $g_{c}$ (combination) as
\begin{align*}
g_{o}(u) &= [u_1,u_3,\cdots,u_9], \quad \forall u\in\mathbb{R}^{10},\\
g_{e}(u)&= [u_2,u_4,\cdots,u_{10}], \quad \forall u\in\mathbb{R}^{10},\\
g_{c}(u',u'') &= [u'_1,u''_1,u'_2,u''_2,\cdots,u'_5,u''_5], \quad \forall u',u''\in\mathbb{R}^{5}.
\end{align*}

Let $\mathcal{T}_r$ be a training dataset of size $N_t=|\mathcal{T}_r|$ (taken to be $100$, $500$, and $1000$),
and $\mathcal{T}_e$ be the test data set of size $|\mathcal{T}_e|=10000$.
We train the INN for the forward and inverse processes as
constructed above simultaneously, denoted by $\Phi_{\nn}$ and $\Phi_{\nn}^{-1}$,
respectively, on the data in the training set $\mathcal{T}_r$.
We define $\hat{y}_{\nn}=\Phi_{\nn }(\hat{u})$ and $\hat{u}_{\nn}=\Phi_{\nn}^{-1}(\hat{y})$.
The loss function $\ell$ that we employ for the training is defined by
\begin{equation}\label{eqn:loss}
   \ell(\hat{u},\hat{y})=\frac{c_0}{2}\sum_{(\hat{u},\hat{y})\in \mathcal{T}_r}\|(\hat{u}-\hat{u}_{\nn})\odot
   w_u\|^2+\frac{1}{2}\sum_{(\hat{u},\hat{y})\in \mathcal{T}_r}\|(\hat{y}-\hat{y}_{\nn})\odot w_y\|^2,
\end{equation}
where $c_0$ is a penalty parameter to be specified according to the relative
errors of forward and inverse processes, and $\odot$ denotes the Hadamard product between two vectors.
In the experiment, we fix $c_0$ at $0.001$ for problem \eqref{eqn:pde}, since
the magnitude of the loss of the inverse process $\sum_{(\hat{u},\hat{y})\in \mathcal{T}_r}\|(\hat{u}-\hat{u}_{\nn})\odot w_u\|^2$
is more than 100 times that of the forward process $\sum_{(\hat{u},\hat{y})\in \mathcal{T}_r}\|
(\hat{y}-\hat{y}_{\nn})\odot w_y\|^2$, and the Lipschitz constant of the inverse process is much bigger than
that of the forward process. The
optimization problem \eqref{eqn:loss} is minimized by the Adam algorithm \cite{KingmaBa:2015} with
a learning rate $r_0=0.001$.

First, we compare our results to that obtained by two different 5-layer fully-connected NN (FNN)
for forward and inverse processes of problem \eqref{eqn:pde}, respectively. Likewise we run Adam \cite{KingmaBa:2015}
for at most 500000 steps to minimize the resulting optimization problem. The numerical results are given in Table
\ref{tab:err}, where $e_{g}$ and $e_{a}$ denote respectively the smallest relative generalization
error (on the test dataset $\mathcal{T}_e$) and the relative approximation error
(of $N_t$ training data points) of the NN that trained on the training dataset $\mathcal{T}_r$ along
the iteration trajectory. The numerical results show that the proposed INN approach can achieve both approximation
 and generalization errors comparable with that by the FNN for all sizes of training dataset for the
forward process; and can also achieve the generalization accuracy largely comparable with that by the FNN
for the inverse process, albeit FNN tends to do a slightly better job. Note that a carefully tuned $c_0$ is needed for different problem settings.
In practice, we only have access to a very limited amount of training data, for which both INN and FNN can achieve
acceptable accuracy. However, INN allows training forward and inverse processes simultaneously
using one single NN.

\begin{table}[htp!]
  \centering
  \begin{threeparttable}
  \caption{The comparison between the proposed INN and FNN (for forward and inverse processes separately).\label{tab:err}}
    \begin{tabular}{cccccccccccc}
    \toprule
    \multicolumn{2}{c}{Neural Network}&
    \multicolumn{2}{c}{INN}&\multicolumn{1}{c}{FNN (forward)}&\multicolumn{1}{c}{FNN (inverse)}\\
    \cmidrule(lr){3-4} \cmidrule(lr){5-5} \cmidrule(lr){6-6}
    $|\mathcal{T}_r|$ &error& forward&inverse&forward&inverse\\
    \midrule
    $100$& $e_{a}$ &1.00e-3 & 1.67e-2 & 1.40e-3 &3.40e-3\\
       & $e_{ g}$ & 1.47e-2 & 2.90e-2 & 1.95e-2 &1.10e-2\\
    \hline
    $500$& $e_{ a}$ &3.10e-3 & 1.83e-2& 8.00e-4 & 3.50e-3\\
       & $e_{ g}$ &9.90e-3  &3.54e-2 & 7.72e-3 &9.00e-3\\
    \hline
    $1000$& $e_{ a}$ &9.00e-4 & 2.95e-2& 9.00e-4&3.70e-3\\
       & $e_{ g}$ &6.99e-3 & 4.17e-2& 6.38e-3&7.70e-3\\
    \bottomrule
    \end{tabular}
    \end{threeparttable}
\end{table}

\begin{figure}[hbt!]
\centering
\setlength{\tabcolsep}{4pt}
\begin{tabular}{cc}
\includegraphics[width=0.48\textwidth]{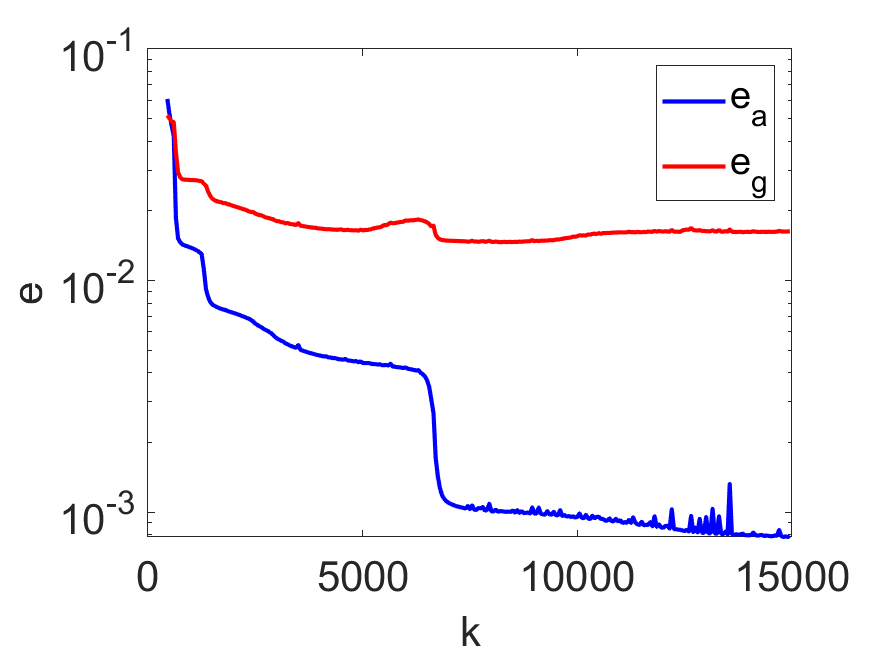}& \includegraphics[width=0.48\textwidth]{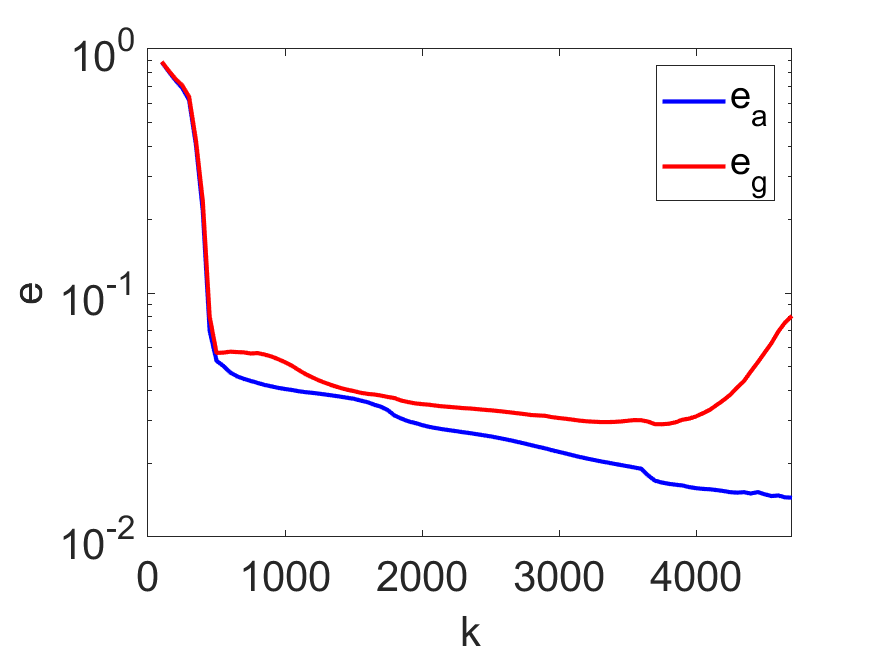}\\
\includegraphics[width=0.48\textwidth]{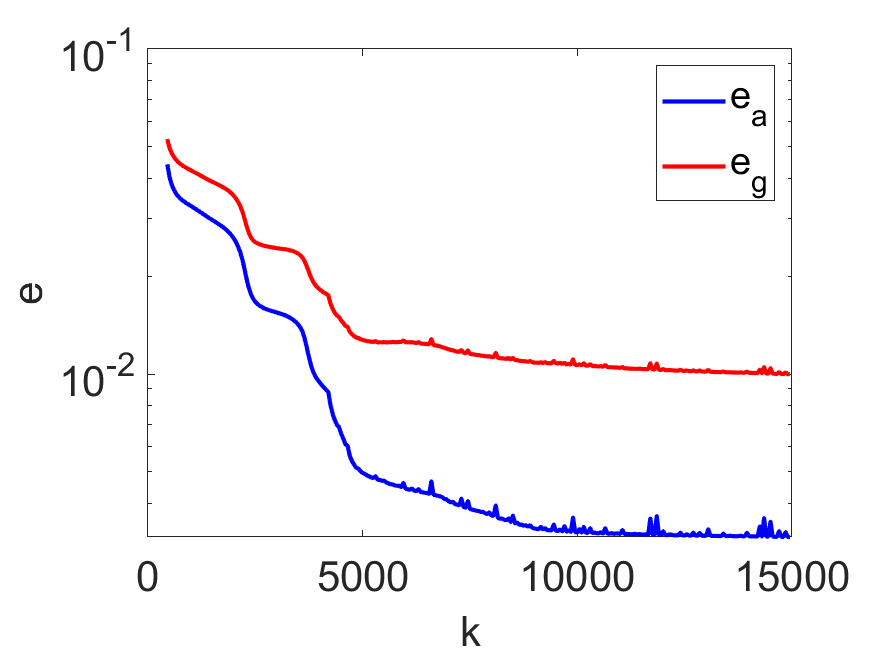}&\includegraphics[width=0.48\textwidth]{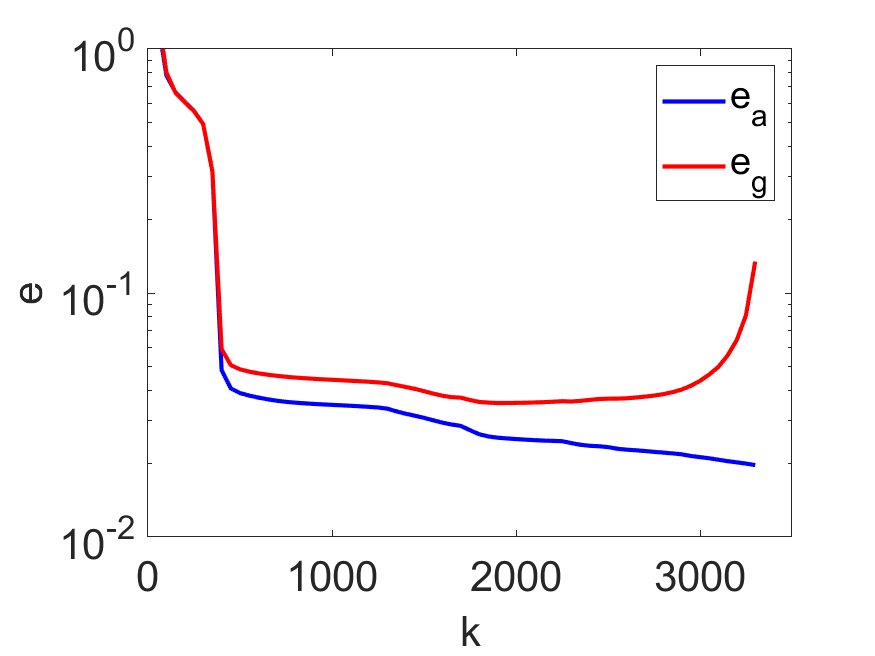}\\
\includegraphics[width=0.48\textwidth]{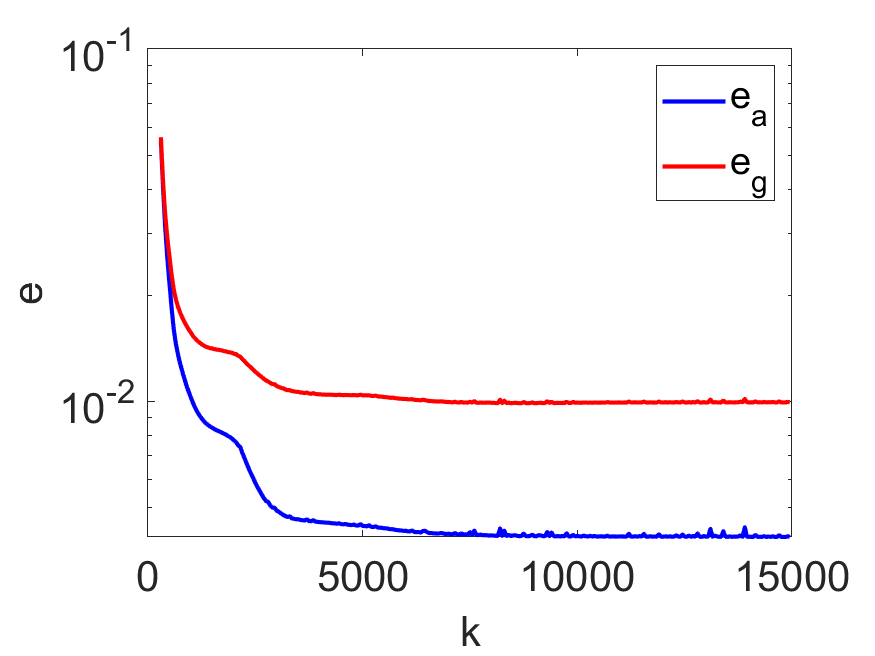}&\includegraphics[width=0.48\textwidth]{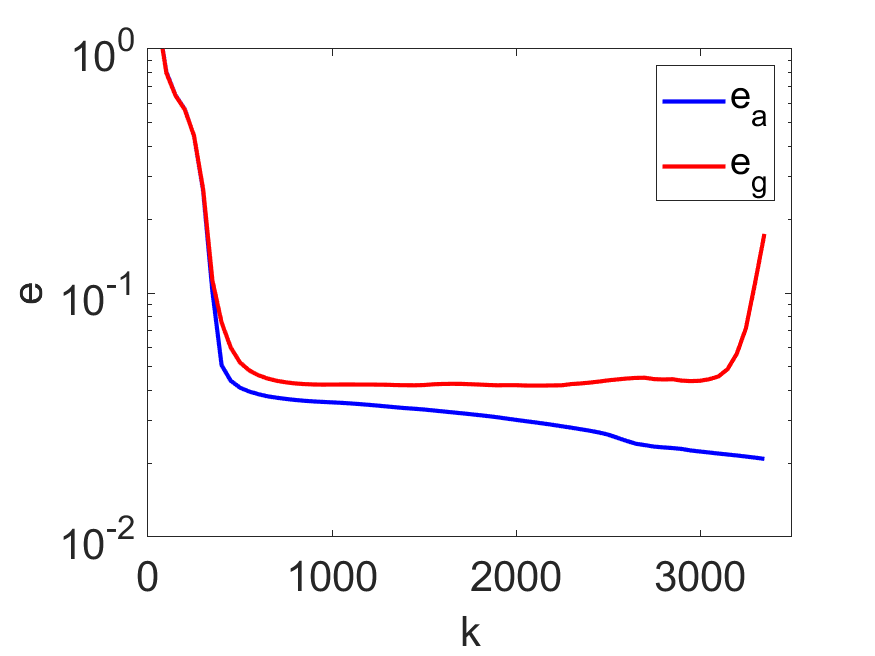}\\
{(a) forward process}& { (b) inverse process}
\end{tabular}
\caption{The relative error $e$ of the proposed INN of the forward process and inverse processes, respectively, versus
the iteration number $k$ of ADAM, when the size $N_t$ of the training dataset varies from 100 (top), 500 (middle) and
1000 (bottom).\label{fig:convergence}}
\end{figure}

Now, we examine more closely the convergence behaviour of the Adam iterates for the training with 100, 500 and 1000
training data, see Fig. \ref{fig:convergence}. We observe that both relative approximation and generalization errors
for the forward and inverse processes decay rapidly at first several iterations and then reach a steady state. This
behavior is commonly observed for many problems. Additionally it is observed that the generalization error $e_g$
exhibits a semi-convergence phenomenon, especially for the inverse process: it first decreases steadily up to a finite number of
iterations, and then starts to increase (occasionally in a dramatic manner) as the iteration further proceeds. Hence,
suitable regularization might be necessary to overcome the phenomenon, e.g., early stopping or suitable explicit
regularization. In contrast, the approximation error remains fairly stable throughout, indicating a stable optimization
process. Compared with the inverse process, the forward process tends to be more stable numerically, since the forward map $F$ enjoys a far more favorable stability estimate and the inverse map $F^{-1}$ enjoys
only conditional (H\"{o}lder type) stability estimates (see \cite{JinZhou:2021} and references therein for details).

\begin{figure}[hbt!]
\centering
\setlength{\tabcolsep}{4pt}
\begin{tabular}{cc}
\includegraphics[width=0.48\textwidth]{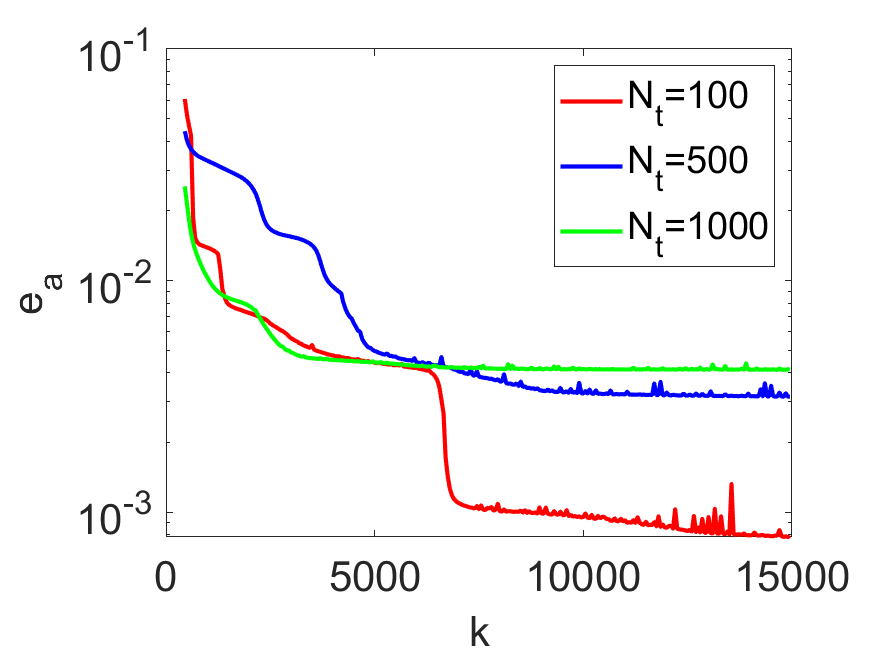}& \includegraphics[width=0.48\textwidth]{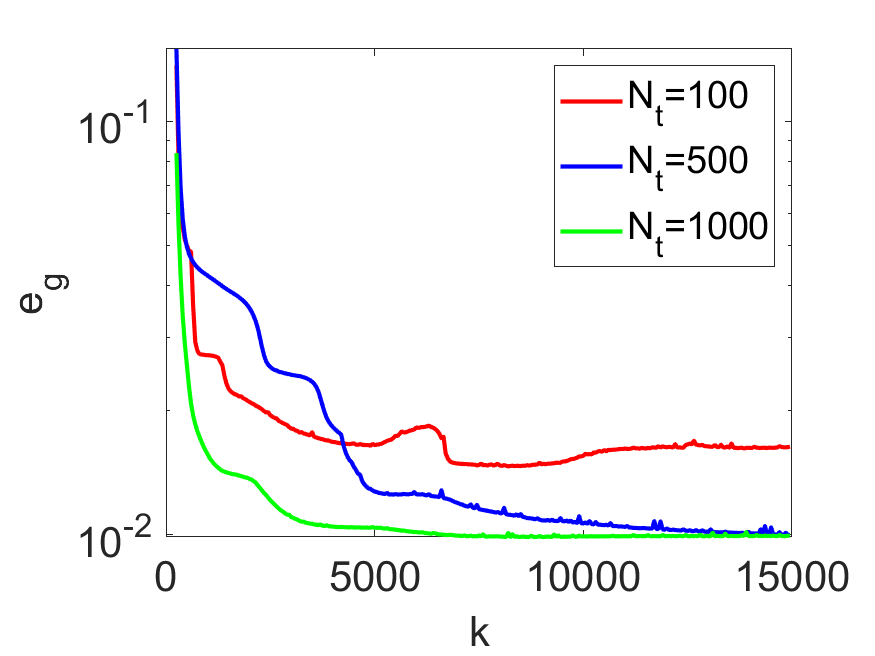}\\
(a) {approximation error $e_a$}& (b) {generalization error $e_g$}
\end{tabular}
\caption{The evolution of the approximation error $e_a$ (a) and generalization error $e_g$ (b) versus the Adam iteration
number $k$ of the proposed INN of the forward process using training dataset of three different sizes $N_t$ (100, 500, 1000).\label{fig:trainingsize}}
\end{figure}

\begin{figure}[hbt!]
\centering
\setlength{\tabcolsep}{4pt}
\begin{tabular}{cc}
\includegraphics[width=0.48\textwidth]{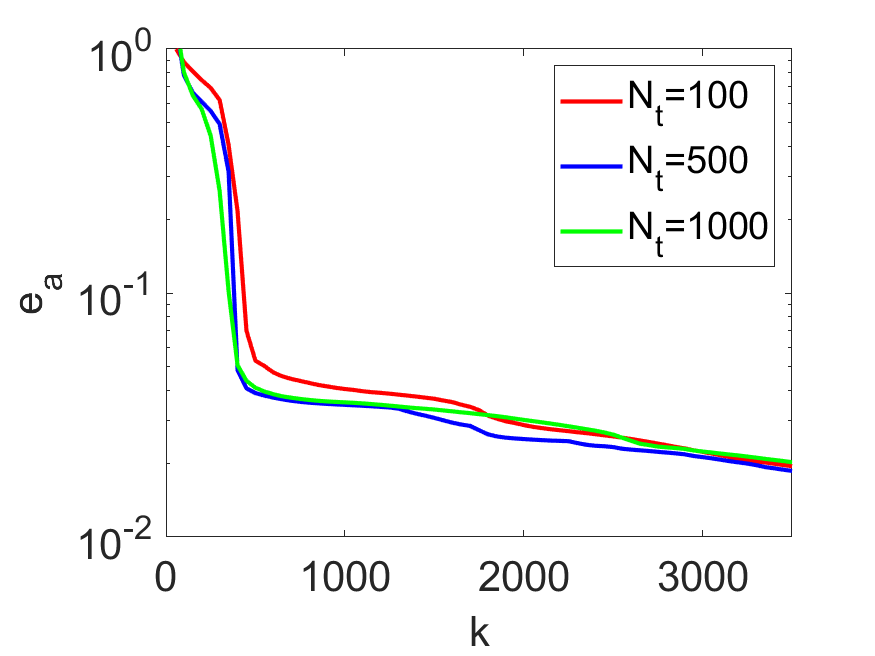}& \includegraphics[width=0.48\textwidth]{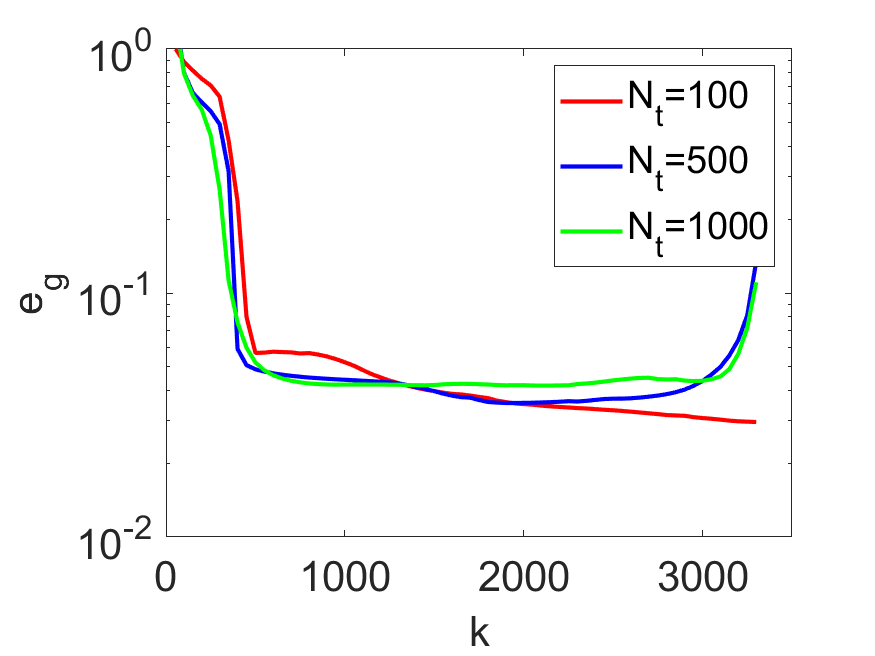}\\
(a) {approximation error $e_a$}& (b) {generalization error $e_g$}
\end{tabular}
\caption{The evolution of the approximation error $e_a$ (a) and generalization error $e_g$ (b) versus the Adam iteration number $k$ of the proposed
INN of the inverse process using training dataset of three different sizes $N_t$ (100, 500, 1000).\label{fig:trainingsizeinv}}
\end{figure}

Finally, the results in Figs. \ref{fig:trainingsize} and \ref{fig:trainingsizeinv} indicate that the training data size
$N_t=500$ is sufficient for problem \eqref{eqn:pde} to achieve reasonable accuracy. A smaller training set yields lower
approximation errors, but performs worse on the generalization; Larger training datasets lead to comparable
performance in terms of both approximation and generalization. This observation agrees well with the intuition that within the
supervised learning paradigm, often one needs a lot of training data in order to achieve good performance of the trained NNs.

\bibliographystyle{abbrv}
\bibliography{NN}

\appendix

\section{An alternative construction of INN $\tilde F_{\nn}$}

In this appendix, we provide an alternative construction an INN $\tilde F_{\rm nn}$, based on
lifting the map $F$ into $\mathbb{R}^{\bd}$, with $\bd= 2d+2$.
Throughout, we are given a finite collection of evaluations $\{(x^{\alpha} , y^{\alpha})\}_{\alpha
\in[n-1]^{d}}$ on a uniform grid. Then the construction proceeds in the following six steps:
\begin{itemize}
  \item[(i)] start from the input $x^{\alpha}$ and lift it to a higher-dimensional channel space by copying its last component and then insert $d+1$ zeros between the original vector and the new component, i.e., transforming $x^{\alpha}$ into $\hat{x}^\alpha=R_L(x^{\alpha})=\big((x^{\alpha})^t,\boldsymbol{0}_{d+1},x^{\alpha}_d\big)^t\in\mathbb{R}^{\bd}$;
  \item[(ii)] make a small perturbation $\eta$ to $\{\hat{x}^{\alpha}\}_{\alpha\in[n-1]^{d}}$ such that
  $\{\eta(\hat{x}^{\alpha})_{\bd}\}_{\alpha\in[n-1]^{d}}$ are distinct, i.e., the last component of all
  the mapped points $\{\eta(\hat{x}^{\alpha})\}_{\alpha\in[n-1]^{d}}$ are distinct;
  \item[(iii)] construct a mapping $\varphi^N$ that satisfies that, for any $\alpha\in[n-1]^{d}$, $(\varphi^N\circ\eta(\hat{x}^{\alpha}))_{j}=y^\alpha_{j}$ for $j=1,\cdots,d$ and $(\varphi^N\circ \eta(\hat{x}^{\alpha}))_{j}=\eta(\hat{x}^{\alpha})_{j}$ for any $j=d+1,\ldots, \bd$, i.e., keeping the last $(d+2)$ components of all the points unchanged;
  \item[(iv)] apply a mapping $\varphi_c$ on $\varphi^N\circ\eta(\hat{x}^{\alpha})$ to transform it into $\varphi_c\big(\varphi^N\circ\eta(\hat{x}^{\alpha})\big)=\big(\big((\varphi^N\circ\eta(\hat{x}^{\alpha}))'^1\big)^t,0,
      \big((\varphi^N\circ\eta(\hat{x}^{\alpha}))'^1\big)^t,\big(\varphi^N\circ\eta(\hat{x}^{\alpha}))_{\bd}\big)^t\in \mathbb{R}^{\bd}$, where $(\varphi^N\circ\eta(\hat{x}^{\alpha}))'^1$ is the subvector that contains the first $d$ components in the vector $\varphi^N\circ\eta(\hat{x}^{\alpha})$;
  \item[(v)] derive some $\tilde{\varphi}^N$ so that $[\tilde{\varphi}^N\circ \varphi_c\circ \varphi^N\circ\eta (\hat{x}^{\alpha})]_{d+1}=[\tilde{\varphi}^N\circ \varphi_c\circ \varphi^N\circ\eta (\hat{x}^{\alpha})]_{\bd}=0$ for any $\alpha\in[n-1]^{d}$ and $[\tilde{\varphi}^N\circ \varphi_c\circ \varphi^N\circ\eta (\hat{x}^{\alpha})]_{j}=[\varphi_c\circ\varphi^N\circ\eta(\hat{x}^{\alpha})]_{j}$ for any $j=1,\cdots,d, d+2,\cdots,\bd-1$;
  \item[(vi)] project it back to $\mathbb{R}^d$ by discarding the last $(d+2)$ components, using a mapping $R_P$.
\end{itemize}

First, we apply lifting a map $R_L: \mathbb{R}^d\longrightarrow\{\hat x\in\mathbb{R}^{\bd}:\;\hat x_{\bd}= \hat x_d,\;x_j=0, j=d+1,\ldots, \bd-1\}$
to lift the inputs $\{x^\alpha\}_{\alpha\in[n-1]^d}$ into $\{\hat{x}^\alpha\}_{
\alpha\in[n-1]^d}=\{R_L(x^{\alpha})\}_{\alpha\in[n-1]^d}=\big((x^{\alpha})^t,\boldsymbol{0}_{d+1},x^{\alpha}_d\big)^t\in\mathbb{R}^{\bd}\}_{
\alpha\in[n-1]^d}$. The map $R_L: \mathbb{R}^d\longrightarrow\{\hat x\in\mathbb{R}^{\bd}:\;\hat x_d=\hat x_{\bd},\;\hat x_j=0, j= d+1,\ldots,\bd-1\}$ is given by
\begin{equation*}
R_L(x)=V_{L}x,\quad \mbox{with } V_{R_L}=\begin{pmatrix}
    I_d\\
    \tilde{V}_{L}
\end{pmatrix}\in\mathbb{R}^{\bd\times d},
\end{equation*}
where $\tilde V_L\in \mathbb{R}^{(d+2)\times d}$ has only one nonzero, $[\tilde V_L]_{d+2,d}=1$;
and its inverse $R_L^{-1}$ is given by
$R_L^{-1}(\hat x)=V_{L^{-1}}\hat x$, with $V_{L^{-1}}=\begin{pmatrix}
    I_d & O_{d \times (d+2)}
\end{pmatrix}\in\mathbb{R}^{d\times \bd}$.
For any vector $\hat x\in\mathbb{R}^{\bd}$, we denote by $\hat x'^1$ the subvector of the first $d$
components of $\hat x$, and $\hat x'^2$ the subvector from the $d+2$ till $\bd-1$th components
of $\hat x$. The projection mapping $R_P: \{\hat x\in\mathbb{R}^{\bd}:\;\hat x_{d+1}=\hat x_{\bd}=0,\;\hat x'^1=\hat x'^2\}\longrightarrow\mathbb{R}^d$ at step (vi) is given by $R_P(x)=V_{P}x$ with $V_P=V_{L^{-1}}$;
and its inverse $R_P^{-1}$ is given by
\begin{equation*}
R_P^{-1}(x)=V_{P^{-1}}x,\quad \mbox{with }V_{P^{-1}}=\begin{pmatrix}
    I_d\\
    O_{1 \times d}\\
    I_d\\
    O_{1 \times d}
\end{pmatrix}\in\mathbb{R}^{\bd \times d}.
\end{equation*}
The key of the construction is to ensure that each of steps (ii)-(v) can be realized by a coupling-based INN, which is given below. We employ three control families $\mathcal{F}_0$, $\mathcal{F}_1$ and $\mathcal{F}_2$ of one-layer, three-layer and four-layer fully connected NN architectures defined by
\begin{align*}
\mathcal{F}_0:=&\big\{V\cdot+b^{(1)}\;|\; V\in\mathbb{R}^{\bd\times\bd},\; b^{(1)}\in\mathbb{R}^{\bd}\big\},\\
\mathcal{F}_1:=&\big\{V\sigma\big(W^{(2)}\sigma(W^{(1)}\cdot+b^{(1)})+b^{(2)}\big)\;|\; V, W^{(1)}, W^{(2)}\in\mathbb{R}^{\bd\times\bd},\; b^{(1)},b^{(2)}\in\mathbb{R}^{\bd}\big\},\\ \mathcal{F}_2:=&\big\{V\sigma\big(W^{(3)}\big(\sigma\big(W^{(2)}\sigma(W^{(1)}\cdot+b^{(1)})+b^{(2)}\big)+b^{(3)}\big)\big)\;| V, W^{(i)}\in\mathbb{R}^{\bd\times \bd}, b^{(i)}\in\mathbb{R}^{\bd}, i=1,2,3\big\}.
\end{align*}
Also we define the attainable set $\mathcal{A}_{{\mathcal{F}}_i}$ corresponding to $\mathcal{F}_i$.

For the uniform grids $\{x^\alpha=\frac{\alpha}{n}\}_{\alpha\in[n-1]^{d}}$ of the hypercube
$K=[0,1]^{d}$, and let $\hat{x}^\alpha=R_L(x^{\alpha})=\big((x^{\alpha})^t,\boldsymbol{0}_{d+1},x^{\alpha}_d\big)^t\in\mathbb{R}^{\bd}$. Let
$\Delta_j(\hat{x}^\alpha)=\min_{\alpha_1\neq\alpha_2}|\hat{x}^{\alpha_1}_j-\hat{x}^{\alpha_2}_j|$
be the minimal distance in the $j$th coordinate of the collection of points $\{\hat{x}^\alpha\}_{
\alpha\in[n-1]^d}$. In the construction of the INN $\tilde F_{\rm nn}$ below, we need $\Delta_{\bd}
(\hat{x}^\alpha)>0$ at the beginning. The mapping $\eta\in\mathcal{A}_{\mathcal{F}_1}$
in the next lemma ensures $\Delta_{\bd}\big(\eta(\hat{x}^\alpha)\big)=N^{-1}$,
thereby fulfilling step (ii) of the overall procedure.

\begin{lemma}\label{lem:perturb_1-v2}
There exists an invertible mapping $\eta\in\mathcal{A}_{\mathcal{F}_1}$ such that
\begin{align*}
&(\eta(\hat x))_{j}=\hat x_{j},\quad j=1,\ldots,\bd-1,\; \hat x\in\mathbb{R}^{\bd},\quad \Delta_{\bd}\big(\eta(\hat{x}^\alpha)\big)=N^{-1},\\
&\max_{\hat x\in\mathbb{R}^{\bd}}\|J_{\eta}(\hat x)\|_2\leq\frac{n}{n-1}\quad \mbox{and}\quad \max_{\hat x\in\mathbb{R}^{\bd}}\|J_{\eta^{-1}}(\hat x)\|_2\leq\frac{n}{n-1}.
\end{align*}
\end{lemma}
\begin{proof}
The proof is similar to Lemma \ref{lem:perturb_1}, but with
$f_j(\hat x)={\rm diag}(\boldsymbol{0}_{2d+1},1)\sigma({\rm diag}(\boldsymbol{0}_{2d+1},\hat x_{j}))\in\mathcal{F}_1$, $j=1,\cdots, d-1$.
Then the map
$\eta:=\psi^{f_{d-1}}_{n^{-(d-1)}}\circ\cdots\circ\psi^{f_2}_{n^{-2}}\circ \psi^{f_1}_{n^{-1}}\in\mathcal{A}_{\mathcal{F}_1}$ satisfies
the desired assertion.
\end{proof}

\begin{remark}\label{rem:perturb_1-v2}
The forward process of the invertible mapping $\eta: \hat x\to \hat y$ in Lemma \ref{lem:perturb_1-v2} is equivalent
to a coupling-based INN {\rm(}with one weight layer and one added layer,
i.e., identity mapping, and $\bd$ neurons on each layer{\rm)}.
\end{remark}

Now we introduce the hypercube centered at $\hat{x}^\alpha$, each side of length $2s$:
\begin{equation*}
S_{\hat{x}^{\alpha}}(s):= (\hat{x}^{\alpha}_1-s,\hat{x}^{\alpha}_1+s)\times\cdots\times(\hat{x}^\alpha_{\bd}-s,\hat{x}^{\alpha}_{\bd}+s).
\end{equation*}
Then we construct an NN, for any $\alpha\in[n-1]^{d}$, to transport
$x^{\alpha}=(\hat{x}^{\alpha})'^1$ to $y^{\alpha}(=F(x^{\alpha}))$, and keep the last $d+2$ components
of $\hat x$ for any $\hat x\in\mathbb{R}^{\bd}$ and $\hat x$ for any $\hat x\in (S_{\hat{x}^{\alpha}}(\frac{1}{2N}))^c=
\mathbb{R}^{\bd}\setminus S_{\hat{x}^{\alpha}}(\frac{1}{2N})$ unchanged.

\begin{lemma}\label{lem:trans_1-v2}
Let $\{\hat{x}^{\alpha}_{\bd}\}_{\alpha\in[n-1]^{d}}$ be distinct with the smallest distance $N^{-1}$.
Then for any $\alpha_1\in[n-1]^{d}$, there exists an invertible mapping $\varphi\in\mathcal{A}_{\mathcal{F}_1}$,
such that
\begin{align*}
\varphi(\hat{x}^{\alpha_1})'^1=y^{\alpha_1},\quad
(\varphi(\hat x))_{j}=&\hat x_{j}, \; j=1,\ldots, d+1,\; \hat x\in\mathbb{R}^{\bd}, \quad\mbox{and}\quad
\varphi(\hat x)=\hat x , \; \forall \hat x\in (S_{\hat{x}^{\alpha_1}}(\tfrac{1}{2N}))^c.
\end{align*}
Furthermore, the following estimates hold
\begin{align*}
\max_{\hat x\in\mathbb{R}^{\bd}}\|J_\varphi(\hat x)\|_2\leq 1+6N\|y^{\alpha_1}-x^{\alpha_1}\|_2 \quad \mbox{and}
\quad \max_{\hat x\in\mathbb{R}^{\bd}}\|J_{\varphi^{-1}}(\hat x)\|_2\leq 1+6N\|y^{\alpha_1}-x^{\alpha_1}\|_2.
\end{align*}
\end{lemma}
\begin{proof}
The proof of the lemma is similar to Lemma \ref{lem:trans_1}. We define an intermediate mapping
\begin{align*}
\ell(\hat x)&=\sigma\big(-\sigma(\tfrac12
W^0\hat x)+b^0\big)-\sigma\big(-\sigma(W^0
\hat x)+b^0\big):=\ell^{(1)}(\hat x) + \ell^{(2)}(\hat x),
\end{align*}
with the weight matrix $W^0=\begin{pmatrix}O_{\bd \times(\bd-1)} &\boldsymbol{1}_{\bd}
\end{pmatrix} \in \mathbb{R}^{\bd\times\bd}$ and bias $b^0= \boldsymbol{1}_{\bd}\in\mathbb{R}^{\bd}$.
Using $\ell_0$ in \eqref{eq:ell_0}, the mapping $\ell(x)$ can be written as
$\ell(x)=\ell_0(1,\cdots,1)^t$. Let
\begin{align*}
A={\rm diag}(\boldsymbol{0}_{\bd-1},2\Delta_{\bd}^{-1}(\hat{x}^\alpha))\in\mathbb{R}^{\bd\times\bd} \quad \mbox{and}\quad
b^{\alpha_1}=(\boldsymbol{0}_{\bd-1},1-2\Delta_{\bd}^{-1}(\hat{x}^\alpha)\hat{x}^{\alpha_1}_{\bd})^t\in\mathbb{R}^{\bd}.
\end{align*}
Then repeating the argument in Lemma \ref{lem:trans_1}, for any $\alpha_1\in[n-1]^{d}$, the following identities hold
\begin{align*}
(A\hat{x}^{\alpha_1}+b^{\alpha_1})_{\bd}=1 \quad\mbox{and}\quad
(A\hat x+b^{\alpha_1})_{\bd}=2\Delta_{\bd}^{-1}(\hat{x}^\alpha)(\hat x_{\bd}-\hat{x}^{\alpha_1}_{\bd})\notin (0,2), \quad \forall \hat x\in (S_{\hat{x}^{\alpha_1}}(\tfrac{1}{2N}))^c.
\end{align*}
Next, we define, for any $\alpha_1\in[n-1]^{d}$, that $D_{\alpha_1}={\rm diag}(2(y_1^{\alpha_1}-\hat{x}_1^{\alpha_1}),\cdots,2(y_{d}^{\alpha_1}-\hat{x}_{d}^{\alpha_1}),\boldsymbol{0}_{d+2})\in\mathbb{R}^{\bd\times \bd}$.
Then we claim that $\varphi:=\psi_1^f$ with $f(\hat x)=D_{\alpha_1}\ell(A\hat x+b^{\alpha_1})$ is the desired mapping
if $\varphi\in\mathcal{A}_{\mathcal{F}_1}$ and it is invertible. The argument in Lemma \ref{lem:trans_1} implies
\begin{align*}
f(\hat x)=&h_{\alpha_1}(\hat x_{\bd})\begin{pmatrix}
y^{\alpha_1}-(\hat{x}^{\alpha_1})'^1\\
\boldsymbol{0}_{d+2}
\end{pmatrix}=h_{\alpha_1}(\hat x_{\bd})\begin{pmatrix}
y^{\alpha_1}-x^{\alpha_1}\\
\boldsymbol{0}_{d+2}
\end{pmatrix},
\end{align*}
with $h_{\alpha_1}$ defined in \eqref{eq:h_1}.
Then there holds
\begin{align*}
J_\varphi(\hat x)=&I+\frac{{\rm d}h_{\alpha_1}(\hat x_{\bd})}{{\rm d}x_{\bd}}\begin{pmatrix}
O_{d\times(2d+1)}&y^{\alpha_1}-x^{\alpha_1}\\
O_{(d+2)\times (2d+1)}&\boldsymbol{0}_{d+2}
\end{pmatrix}:=I+H_{\alpha_1}(\hat x).
\end{align*}
It follows from direct computation that
\begin{align*}
  \Big|\frac{{\rm d}h_{\alpha_1}(\hat x_{\bd})}{{\rm d}\hat x_{\bd}}\Big|\leq2 (\Delta_{\bd}^{-1}(\hat{x}^\alpha)+2\Delta_{\bd}^{-1}
  (\hat{x}^\alpha))= 6\Delta_{\bd}^{-1}(\hat{x}^\alpha)=6N \quad \mbox{and}\quad {\rm det}(J_\varphi(\hat x))=1.
\end{align*}
The bound on $J_\varphi(\hat x)$ and $J_{\varphi^{-1}}(\hat x) $ follows identically as Lemma \ref{lem:trans_1}.
\end{proof}

\begin{remark}\label{rem:trans_1-v2}
The forward process of the mapping $\varphi_\alpha: x\to y$ in Lemma \ref{lem:trans_1-v2} is equivalent to a coupling-based INN {\rm(}with six weight layers and two added layers {\rm(}i.e.,
identity mapping{\rm)}, and $\bd$ neurons on each layer{\rm)} with the mapping $h_\alpha:\mathbb{R}\to\mathbb{R}$
defined as \eqref{eq:h_1} and likewise the inverse process $\varphi^{-1}$:
\begin{align*}
\left\{\begin{aligned}
  y'^1&= x'^1+h_{\alpha}(x_{\bd})(y^{\alpha}-x^{\alpha}),\\
 y_{j}&= x_{j}, \;j= d+1,\ldots,\bd,
\end{aligned}\right.
 \quad {\rm and}\quad
 \left\{\begin{aligned}
   x_{j}&= y_{j}, \; j= d+1,\ldots,\bd, \\
   x'^1&= y'^1-h_\alpha(x_{\bd})(y^{\alpha}-x^{\alpha}).
 \end{aligned}\right.
\end{align*}
\end{remark}

Having successfully transported the set $\{\hat{x}^\alpha\}_{\alpha\in[n-1]^{d}}$
to $\{z^\alpha\}_{\alpha\in[n-1]^{d}}$ with $z^\alpha=\big((y^\alpha)^t,\boldsymbol{0}_{d+1},\eta(\hat{x}^{\alpha})_{\bd}\big)^t$,
following the idea of Lemma \ref{lem:trans_1-v2}, we can transport $\{z^\alpha\}_{\alpha\in [n-1]^d}$ of points to $\{((y^\alpha)^t,0,(y^\alpha)^t,0)^t\}_{\alpha\in[n-1]^{d}}$. To
this end, we first transform $z^\alpha$ into $\hat{z}^\alpha=\varphi_c(z^\alpha)=\big(\big((z^\alpha)'^1\big)^t,0,\big((z^\alpha)'^1\big)^t,z^\alpha_{\bd}\big)^t\in\mathbb{R}^{\bd}$ using a mapping $\varphi_c: \mathbb{R}^{\bd}\longrightarrow\mathbb{R}^{\bd}$.
\begin{lemma}\label{lem:copy_1}
For any $\alpha\in[n-1]^{d}$, let $z^\alpha=\big((y^\alpha)^t,\boldsymbol{0}_{d+1},\eta(\hat{x}^{\alpha})_{\tilde{d}}\big)^t$. Then, there exists an invertible mapping $\varphi_c\in\mathcal{A}_{\mathcal{F}_0}$,
such that
\begin{align}\label{eq:copy}
\varphi_c(z^{\alpha})'^1=\varphi_c(z^{\alpha})'^2=y^{\alpha},\quad
(\varphi_c(z))_{j}=&z_{j}, \;  j=1,\ldots, d+1 \mbox{ or }j=\bd, \mbox{ and } z\in\mathbb{R}^{\bd}.
\end{align}
Furthermore, the following estimates hold
\begin{align*}
\max_{z\in\mathbb{R}^{\bd}}\|J_{\varphi_c}(z)\|_2\leq 1+\sqrt{d} \quad \mbox{and}
\quad \max_{z\in\mathbb{R}^{\bd}}\|J_{{\varphi_c}^{-1}}(z)\|_2\leq 1+\sqrt{d}.
\end{align*}
\end{lemma}
\begin{proof}
Note that the linear mapping $\varphi_c(z)=z+f(z):=z+Wz$, with $$W=\begin{pmatrix}
O_{(d+1)\times(d+1)}&O_{(d+1)\times(d+1)}\\
\tilde{W}& O_{(d+1)\times(d+1)}
\end{pmatrix} \quad\mbox{with}\quad \tilde{W}={\rm diag}(\boldsymbol {1}_{d},0)\in \mathbb{R}^{(d+1)\times (d+1)},$$
satisfies \eqref{eq:copy} and $\varphi_c=\psi_1^f\in\mathcal{A}_{\mathcal{F}_0}$. The mapping $\varphi_c$ is invertible with $\varphi_c^{-1}=\psi_1^{-f}$ and
$J_{\varphi_c}(z)=I+W$ and $J_{{\varphi_c}^{-1}}(z)=I-W$.
Thus,
${\rm det}(J_\varphi(z))=1$ and $\max\big(\max_{z\in\mathbb{R}^{\bd}}\|J_{\varphi_c}(z)\|_2\;,\max_{z\in\mathbb{R}^{\bd}}\|J_{{\varphi_c}^{-1}}(z)\|_2\big)\leq 1+\sqrt{d}$.
This completes the proof of the lemma.
\end{proof}
\begin{remark}\label{rem:copy_1}
The forward process of the mapping $\varphi_c: z\to y$ in Lemma \ref{lem:copy_1} is equivalent to a
coupling-based NN {\rm(}with one weight layer and one added layer {\rm(}i.e.,
identity mapping{\rm)}, and $\bd$ neurons on each layer{\rm)} and likewise the inverse process $\varphi_c^{-1}$:
\begin{align*}
y_j  = \left\{\begin{aligned}
      z_j+z_{j-(d+1)}, &\quad j=d+2,\ldots, 2d+1,\\
      z_{j}, &\quad\mbox{otherwise},
      \end{aligned}\right.\quad \mbox{and}\quad
        z_j=\left\{\begin{aligned}
    y_j-y_{j-(d+1)}, &\quad j=d+2,\ldots, 2d+1,\\
     y_{j}, &\quad\mbox{otherwise}.
  \end{aligned}\right.
\end{align*}
\end{remark}

Now we introduce a hypercube in $\mathbb{R}^d$, centered at $y^\alpha$, each side of length $2s$:
$S_{y^{\alpha}}(s):= (y^{\alpha}_1-s,y^{\alpha}_1+s)\times\cdots\times(y^\alpha_{d}-s,y^{\alpha}_{d}+s)$.
Since $F$ is bi-Lipschitz continuous, we have
\begin{equation*}
   \|y^{\alpha_1}-y^{\alpha_2}\|_2\geq \Lip{F^{-1}}^{-1}\|F^{-1}(y^{\alpha_1})-F^{-1}(y^{\alpha_2})\|_2=\Lip{F^{-1}}^{-1}\|x^{\alpha_1}-x^{\alpha_2}\|_2,\quad \forall \alpha_1 \neq\alpha_2.
\end{equation*}
For any fixed $\alpha_1\neq\alpha_2$, there must exist an index $j\leq d$ such that
\begin{equation*}
   |y_{j}^{\alpha_1}-y_{j}^{\alpha_2}|\geq d^{-\frac12} \Lip{F^{-1}}^{-1}\|x^{\alpha_1}-x^{\alpha_2}\|_2\geq d^{-\frac12}\Lip{F^{-1}}^{-1}n^{-1}:=\Delta.
\end{equation*}
Thus, the set
$\{S_{y^{\alpha}}(\frac{\Delta}{2})\}_{\alpha\in[n-1]^{d}}$ consists of a collection of disjoint
open hypercubes in $\mathbb{R}^{d}$.
Then we construct an NN applying to $z\in\mathbb{R}^{\tilde{d}}$, for any $\alpha\in[n-1]^{d}$, to transport
$\hat{z}^{\alpha}=\varphi_c(z^{\alpha})$ to $\hat{y}^\alpha:=\big((y^{\alpha})^t,0,(y^{\alpha})^t,0\big)^t$, i.e., keeping $z_j$ for any $j=1,\cdots,\tilde{d}-1$ and $z\in\mathbb{R}^{\tilde{d}}$, and $z_{\tilde{d}}$ for any $z$ such that $z'^1=z'^2\in (S_{y^{\alpha}}(\frac{\Delta}{2}))^c=
\mathbb{R}^{d}\setminus S_{y^{\alpha}}(\frac{\Delta}{2})$ unchanged.

\begin{lemma}\label{lem:trans_2-v2}
Let $\hat{z}^\alpha=\big((z^\alpha)'^t,0,(z^\alpha)^t\big)^t=\big((y^\alpha)^t,0,(y^\alpha)^t,
\eta(\hat{x}^{\alpha})_{\tilde{d}}\big)^t\in\mathbb{R}^{\bd}$ for any $\alpha\in[n-1]^{d}$.
Then for any fixed $\alpha_1\in[n-1]^{d}$, there exists an invertible mapping $\tilde{\varphi}
\in\mathcal{A}_{\mathcal{F}_2}$, such that, for any $i=1,2$,
\begin{align*}
(\tilde{\varphi}(\hat{z}^{\alpha_1}))_{i(d+1)}=0,\quad
\tilde{\varphi}(z)'^i=z'^i&, \; \forall z\in\mathbb{R}^{\bd} \quad\mbox{and}\quad
\tilde{\varphi}(z)=z , \; \forall z \mbox{ with } z'^1=z'^2\in (S_{y^{\alpha}}(\tfrac{\Delta}{2}))^c.
\end{align*}
Furthermore, the following estimates hold
\begin{align*}
\max_{z\in\mathbb{R}^{\bd}}\|J_{\tilde\varphi}(z)\|_2\leq 1+6\Lip{F^{-1}}^{-1}n^{-1} \quad \mbox{and}
\quad \max_{z\in\mathbb{R}^{\bd}}\|J_{{\tilde\varphi}^{-1}}(z)\|_2\leq 1+6\Lip{F^{-1}}^{-1}n^{-1}.
\end{align*}
\end{lemma}
\begin{proof}
Similar to the proof of Lemma \ref{lem:trans_1}, we define an intermediate mapping
\begin{align*}
\ell(z)&=\sigma\big(-\sigma(W^0
z)+b^0\big),\quad \mbox{with }
  W^0 =
\begin{pmatrix}
\tfrac12 I_{d+1}&0\\
0& I_{d+1}
\end{pmatrix}\in \mathbb{R}^{\bd\times\bd}, b^0 = \boldsymbol{1}_{\bd}\in\mathbb{R}^{\bd}.
\end{align*}
The mapping $\ell(z)$ can be equivalently written as
$$\ell(z)=(\ell_1(z_1),\cdots,\ell_1(z_{d+1}),-\ell_2(z_{d+2}),\cdots,-\ell_2(z_{\bd}))^t,$$
with $\ell_1$ and $\ell_2$ defined in \eqref{eq:ell_0} and $\ell_0=\ell_1+\ell_2$.
Next we define
$A=2\Delta^{-1}{\rm diag}(\boldsymbol {1}_{d},0,\boldsymbol{1}_{d},0)\in\mathbb{R}^{\bd\times\bd}$, and
$b^{\alpha_1}=(1-2\Delta^{-1}\hat{z}^{\alpha_1}_{1},\cdots,1-2\Delta^{-1}\hat{z}^{\alpha_1}_{d},1,1-2\Delta^{-1}\hat{z}^{\alpha_1}_{d+2},
\cdots,1-2\Delta^{-1}\hat{z}^{\alpha_1}_{2d+1},1)^t\in\mathbb{R}^{\bd}$.
Then for any $\alpha_1\in[n-1]^{d}$ and $j=1,\cdots,d, d+2, \cdots, 2d+1$, the following identities hold
\begin{align*}
&(Az+b^{\alpha_1})_{d+1}=(Az+b^{\alpha_1})_{\bd}=1, \quad \forall z\in\mathbb{R}^{\bd},\\
&(A\hat{z}^{\alpha_1}+b^{\alpha_1})_{j}=1 \quad\mbox{and}\quad
(Az+b^{\alpha_1})_{j}=2\Delta^{-1}(z_{j}-\hat{z}^{\alpha_1}_{j})\notin (0,2), \quad \forall z \mbox{ with } z'^1=z'^2\in (S_{y^{\alpha}}(\tfrac{\Delta}{2}))^c.
\end{align*}
Thus, for any $\alpha_1\in[n-1]^{d}$, we have
\begin{align*}
&\ell_1((A\hat{z}^{\alpha_1}+b^{\alpha_1})_j)+\ell_2((A\hat{z}^{\alpha_1}+b^{\alpha_1})_{j+d+1})=\ell_0((A\hat{z}^{\alpha_1}+b^{\alpha_1})_j)=\tfrac12, \quad \forall j=1,\cdots, d,\\
&\ell_1((Az+b^{\alpha_1})_{d+1})=\tfrac12,\;\quad \ell_2((Az+b^{\alpha_1})_{\bd})=0, \quad \forall z\in \mathbb{R}^{\bd},
\end{align*}
and there is at least one zero in the set $\{\ell_1((Az+b^{\alpha_1})_j)+\ell_2((Az+b^{\alpha_1})_{j+d+1})\}_{j=1}^{d+1}$ for any $z$ such that $z'^1=z'^2\in (S_{y^{\alpha}}(\frac{\Delta}{2}))^c$, which imply
\begin{align*}
s(\ell(Az+b^{\alpha_1}))<\tfrac12 (d-1)-\tfrac12(d-1)=0 \quad \mbox{and}\quad
s(\ell(A\hat{z}^{\alpha_1}+b^{\alpha_1}))=\tfrac12 d-\tfrac12 (d-1)=\tfrac12,
\end{align*}
where the function $s$ is defined by
\begin{align*}
s(z) = \sum_{j=1}^{d}z_j-\sum_{j=d+2}^{2d+1}z_j-(d-1)z_{d+1}.
\end{align*}
Note that the function $s$ can be rewritten as
\begin{align*}
s(z)=W^1 z, \quad \mbox{with } W^1=\begin{pmatrix}
O_{(\bd-1)\times \bd}\\
w
\end{pmatrix},
\end{align*}
with $w=\begin{pmatrix} \boldsymbol{1}_{d}&-(d-1)&-\boldsymbol{1}_d&0\end{pmatrix}\in\mathbb{R}^{1\times\bd}$.
Next, for fixed $\alpha_1\in[n-1]^{d}$, define
$D_{\alpha_1}={\rm diag}(\boldsymbol{0}_{\bd-1},-2\hat{z}^{\alpha_1}_{\bd})\in\mathbb{R}^{\bd\times\bd}$.
Then $\tilde\varphi:=\psi_1^f$ with $f(x)=D_{\alpha_1}\sigma\big(s\big(\ell(Ax+b^{\alpha_1})\big)\big)$ is the desired mapping
if $\tilde\varphi\in\mathcal{A}_{\mathcal{F}_2}$ and it is invertible.
By the definitions of $f$, $s$, $\ell$, $\ell_1$, $\ell_2$ and $\ell_0=\ell_1+\ell_2$, we have $\tilde\varphi\in\mathcal{A}_{\mathcal{F}_2}$, and there holds
\begin{align}
f(x)=&D_{\alpha_1}\sigma\Big(s\big(\ell(Ax+b^{\alpha_1})\big)\Big)=D_{\alpha_1}\sigma\Big(W^1\sigma\big(-\sigma(W^0
(Ax+b^{\alpha_1}))+b^0\big)\Big)\nonumber\\
=&D_{\alpha_1}\sigma(W^1v)
=-2\hat{z}^{\alpha_1}_{\bd}\begin{pmatrix}
\boldsymbol{0}_{\bd-1}\\
h_{\alpha_1}(z)
\end{pmatrix},
 \label{eq:h_2}
\end{align}
with $v=(
\ell_1(2\Delta^{-1}(z_{1}-\hat{z}^{\alpha_1}_{1})+1)\
\ldots \ \ell_1(2\Delta^{-1}(z_{d}-\hat{z}^{\alpha_1}_{d})+1)\
\tfrac 12\
-\ell_2(2\Delta^{-1}(z_{d+2}-\hat{z}^{\alpha_1}_{d+2})+1)\
\ldots\ -\ell_2(2\Delta^{-1}(z_{2d+1}-\hat{z}^{\alpha_1}_{2d+1})+1)\
0
)^t\in\mathbb{R}^{\bd}$
and $h_{\alpha_1}(z)=\sigma\big(\sum_{j=1}^d\big(\ell_1(2\Delta^{-1}(z_{j}-\hat{z}^{\alpha_1}_{j})+1) +\ell_2(2\Delta^{-1}(z_{j+d+1}-\hat{z}^{\alpha_1}_{j+d+1})+1)\big)-\tfrac12(d-1)\big)$.
Then direct computation gives
\begin{align*}
J_{\tilde\varphi}(z)=&I-2\hat{z}^{\alpha_1}_{\bd}\begin{pmatrix}
O_{(\bd-1)\times(\bd-1)}&\boldsymbol{0}_{\bd-1}\\
\tilde h&0
\end{pmatrix}:=I+H_{\alpha_1}(z),
\end{align*}
with $\tilde h =(\frac{{\rm d}h_{\alpha_1}(z)}{{\rm d}z_{1}}\ \cdots \ \frac{{\rm d}h_{\alpha_1}(z)}{{\rm d}z_{2d+1}}) \in\mathbb{R}^{\bd-1}$.
It follows from direct computation that
\begin{align*}
&\Big|\frac{{\rm d}h_{\alpha_1}(z)}{{\rm d}z_j}\Big|\leq \left\{\begin{aligned} \Delta^{-1},& \quad j=1,\cdots,d,\\
 2\Delta^{-1}, &\quad j=d+2,\cdots,2d+1,\end{aligned}\right.
\end{align*}
$|\frac{{\rm d}h_{\alpha_1}(z)}{{\rm d}z_{d+1}}|=0$ and ${\rm det}(J_{\tilde\varphi}(z))=1.$
This directly implies that $\tilde\varphi\in\mathcal{A}_{\mathcal{F}_2}$ is invertible and furthermore,
\begin{align*}
\|J_{\tilde\varphi}(z)\|_2\leq& 1+\|H_{\alpha_1}(z)\|_2\leq 1+\sqrt{{\rm tr}\big((H_{\alpha_1}(z))^t
H_{\alpha_1}(z)\big)}\leq 1+2|\hat{z}^{\alpha_1}_{\bd}|\sqrt{5d}\Delta^{-1}\\
\leq& 1+6\Lip{F^{-1}}^{-1}|\eta(\hat{x}^{\alpha_1})_{d+1}|n^{-1}\leq 1+6\Lip{F^{-1}}^{-1}n^{-1},\\
\max_{z\in\mathbb{R}^{\bd}}\|J_{{\tilde\varphi}^{-1}}(z)\|_2=&\max_{z\in\mathbb{R}^{\bd}}\|(J_{\tilde\varphi}(z))^{-1}
\|_2=\max_{z\in\mathbb{R}^{\bd}}\|I-H_{\alpha_1}(z)\|_2\leq 1+6\Lip{F^{-1}}^{-1}n^{-1}.
\end{align*}
This completes the proof of the lemma.
\end{proof}

\begin{remark}\label{rem:trans_2-v2}
The forward process of the mapping $\tilde\varphi_\alpha: z\to y$ in Lemma \ref{lem:trans_2-v2} is equivalent to a
coupling-based NN {\rm(}with eight weight layers and two added layers {\rm(}i.e.,
identity mapping{\rm)}, and $\bd$ neurons on each layer{\rm)} with the mapping $h_\alpha:\mathbb{R}^{\bd}\to\mathbb{R}$
defined as \eqref{eq:h_2} which is independent of $z_{\bd}$:
$y_{\bd}= z_{\bd}-2\hat{z}^{\alpha_1}_{\bd}h_{\alpha}(z)=z_{\bd}-2\eta(\hat{x}^{\alpha_1})_{d+1}h_{\alpha}(z)$ and $y'= z'$.
The inverse process ${\tilde\varphi}^{-1}$ can be written as
$  z'=y'$ and $z_{\bd}= y_{\bd}+2\eta(\hat{x}^{\alpha_1})_{d+1}h_{\alpha}(z).$
\end{remark}

Last we construct an alternative INN $\tilde{F}_{\rm nn}$ to approximate the map $F$.
\begin{theorem}\label{thm:F_NN-v2}
Given a finite collection of evaluations, there exists a bi-Lipschitz continuous invertible mapping $\tilde{F}_{\rm nn}\in R_P\circ\mathcal{A}_{\mathcal{F}_2}\circ \mathcal{A}_{\mathcal{F}_0} \circ \mathcal{A}_{\mathcal{F}_1}\circ R_L$,
such that $y^{\alpha}=\tilde{F}_{\rm nn}(x^{\alpha})$ for any $\alpha\in[n-1]^{d}$, and the Lipschitz
constants of $\tilde{F}_{\rm nn}$ and $\tilde{F}_{\rm nn}^{-1}$ are bounded by
\begin{align*}
\Lip{\tilde{F}_{\rm nn}}\leq\frac{4n}{n-1}\sqrt{d}\big(1+6Nc\big)(1+6\Lip{F^{-1}}^{-1}n^{-1})\quad\mbox{and}\quad
\Lip{\tilde{F}_{\rm nn}^{-1}}\leq\frac{4n}{n-1}\sqrt{d}\big(1+6Nc\big)(1+6\Lip{F^{-1}}^{-1}n^{-1}),
\end{align*}
with $c:=\max_{i=1,\cdots N}\|y^{\alpha_i}-x^{\alpha_i}\|_2$.
\end{theorem}
\begin{proof}
First, for any $\alpha\in[n-1]^{d}$, we transform $x^\alpha$ into $\hat{x}^\alpha=R_L(x^{\alpha})=\big((x^{\alpha})^t,\boldsymbol{0}_{d+1},x^{\alpha}_d\big)^t\in\mathbb{R}^{\bd}$.
By Lemma \ref{lem:perturb_1-v2}, there exists an INN $\eta$ such
that $\Delta_{\bd}\big(\eta(\hat{x}^\alpha)\big)=N^{-1}$ for all $\alpha\in [n-1]^d$ and
\begin{align*}
 \max_{x\in\mathbb{R}^{\bd}}\|J_{\eta}(x)\|_2\leq\frac{n}{n-1}\quad \mbox{and}\quad \max_{x\in\mathbb{R}^{\bd}}\|J_{\eta^{-1}}(x)\|_2\leq\frac{n}{n-1}.
\end{align*}
Then, by Lemma \ref{lem:trans_1-v2}, we can construct a composite INN $\varphi^N:=\varphi_{\alpha_N}\circ\cdots\circ
\varphi_{\alpha_1}\in \mathcal{A}_{\mathcal{F}_1}$, where $\big(\varphi^j(\eta(\hat{x}^{\alpha_i}))\big)'^1=y^{\alpha_i}$ for any $1\leq i\leq j\leq N$.
Since each $\varphi_{\alpha_i}$, $i=1,\cdots,N$, is the identity map in $(S_{\eta(\hat{x}^{\alpha_i})}
(\frac{1}{2N}))^c$, and $\{S_{\eta(\hat{x}^{\alpha_i})}(\frac{1}{2N})\}_{i=1}^N$ are disjoint from each other, we have
\begin{align*}
\max_{\hat x\in\mathbb{R}^{\bd}}\|J_{\varphi^N}(\hat x)\|_2&\leq\max_{i=1,\cdots N,\hat x\in\mathbb{R}^{\bd}} \|J_{\varphi_{\alpha_i}}(\hat x)\|_2\leq 1+6N\max_{i=1,\cdots N}\|y^{\alpha_i}-x^{\alpha_i}\|_2,\\
  \max_{\hat x\in\mathbb{R}^{\bd}}\|J_{(\varphi^N)^{-1}}(\hat x)\|_2&\leq 1+6N\max_{i=1,\cdots N}\|y^{\alpha_i}-x^{\alpha_i}\|_2.
\end{align*}
For any $\alpha\in [n-1]^d$, we have $z^\alpha=\varphi^N \circ \eta(\hat{x}^{\alpha})=\big((y^\alpha)^t,\boldsymbol{0}_{d+1},\eta(\hat{x}^{\alpha})_{\bd}\big)^t$ and $\hat{z}^\alpha=\varphi_c(z^\alpha)=\big((y^\alpha)^t,0,(y^\alpha)^t,\eta(\hat{x}^{\alpha})_{\bd}\big)^t$ with the mapping $\varphi_c\in\mathcal{A}_{\mathcal{F}_0}$ from Lemma \ref{lem:copy_1} and
\begin{align*}
\max_{z\in\mathbb{R}^{\bd}}\|J_{\varphi_c}(z)\|_2\leq 1+\sqrt{d}\leq2\sqrt{d} \quad \mbox{and}
\quad \max_{z\in\mathbb{R}^{\bd}}\|J_{{\varphi_c}^{-1}}(z)\|_2\leq 1+\sqrt{d}\leq 2\sqrt{d}.
\end{align*}
Then repeating the argument with Lemma \ref{lem:trans_2-v2}, there exists an INN $\tilde{\varphi}^N:=
\tilde{\varphi}_{\alpha_N}\circ\cdots\circ\tilde{\varphi}_{\alpha_1}\in\mathcal{A}_{\mathcal{F}_2}$ such that
$({\tilde\varphi}^j\circ \varphi_c\circ\varphi^N\circ\eta(\hat{x}^{\alpha_i}))_{\bd}=({\tilde\varphi}^j( \hat{z}^{\alpha_i}))_{\bd}=0$ for any $1\leq i\leq j\leq N$ and
\begin{align*}
   \max_{z\in\mathbb{R}^{\bd}}\|J_{\tilde\varphi}(z)\|_2\leq 1+6\Lip{F^{-1}}^{-1}n^{-1}\quad \mbox{and}\quad
   \max_{z\in\mathbb{R}^{\bd}}\|J_{{\tilde\varphi}^{-1}}(z)\|_2\leq
 1+6\Lip{F^{-1}}^{-1}n^{-1}.
\end{align*}
Finally, let $\tilde{F}_{\rm nn}=R_P\circ\tilde{\varphi}^N\circ \varphi_c\circ \varphi^N \circ \eta\circ R_L\in
R_P\circ\mathcal{A}_{\mathcal{F}_2}\circ \mathcal{A}_{\mathcal{F}_0}\circ \mathcal{A}_{\mathcal{F}_1}\circ R_L$. This is the desired mapping.
Moreover, direct computation shows
$\Lip{R_L}\leq 2$, $\Lip{R_P}\leq 1$, $\Lip{R_L^{-1}}\leq 1$, and $ \Lip{R_P^{-1}}\leq 2$.
Then, the desired estimates on $\Lip{\tilde{F}_{\rm nn}}$ and $\Lip{{\tilde{F}_{\rm nn}}^{-1}}$ follow.
\end{proof}

\begin{remark}
By the construction in Theorem \ref{thm:F_NN-v2} and the analysis in Remarks \ref{rem:perturb_1-v2}-\ref{rem:trans_2-v2}, the INN $\tilde{F}_{\rm nn}$  is actually an invertible lifting layer $R_L$ combined with a coupling-based INN with $2+14N$ weight layers and $2+4N$
added layers {\rm(}at most $\bd$ neurons at each layer{\rm)}, and then followed by an invertible projection layer $R_P$.
\end{remark}

\end{document}